\documentclass[12pt]{article}
\usepackage[utf8]{inputenc}
\usepackage{amsmath, amssymb, amsthm, bbm, mathrsfs}
\usepackage{geometry}
\usepackage{natbib}
\usepackage{verbatim}
\usepackage{multirow}
\usepackage{booktabs, array, threeparttable}
\usepackage{threeparttable}
\usepackage{rotating}
\usepackage{graphicx}
\usepackage{subfigure}
\usepackage{float}
\usepackage{caption}
\usepackage{amsfonts}
\usepackage{xr}
\usepackage{color}
\usepackage{enumerate}
\usepackage{epstopdf}
\usepackage{threeparttable}
\usepackage{siunitx}
\usepackage{multirow}
\usepackage[colorlinks=true, allcolors=blue]{hyperref}
\usepackage{makecell}
\allowdisplaybreaks
\usepackage[linesnumbered,ruled,vlined]{algorithm2e}
\usepackage{geometry}
\graphicspath{{./figs/}}

\newcommand{\bm}{\boldsymbol}

\def\cum{\mathrm {cum}}
\def\cp{\mathop{\longrightarrow}\limits^{\mathcal{P}}}

\def\cov{\mathrm {cov}}
\def\var{\mathrm {Var}}
\def\Pr{\mathbb{P}}
\def\tr{\mathrm{tr}}

\def\bms{{\bm\Sigma}}
\def\bmss{{\Sigma}}

\def\E{\mathbb {E}}

\def\bH{\mathbb{H}}
\def\bR{\mathbb{R}}

\def\done{$\hfill\blacksquare$}
\def\mA{{\mathcal A}}
\def\mB{{\mathcal B}}

\def\mF{{\mathcal F}}

\def\X{{\bm  x}}
\def\Y{{\bm y}}
\def\Z{{\bm z}}

\newtheorem{assumption}{Assumption}
\newtheorem{theorem}{Theorem}
\newtheorem{proposition}{Proposition}
\newtheorem{remark}{Remark}
\newtheorem{lemma}{Lemma}

\title{Simultaneous Detection and Localization of Mean and Covariance Changes in High Dimensions}
\date{}
\author{Junfeng Cui$^1$, Guangming Pan$^2$, Guanghui Wang$^3$ and Changliang Zou$^3$\\
$^1$ {\normalsize\it School of Mathematical Sciences, Shenzhen University} \\
$^2$ {\normalsize\it School of Physical and Mathematical Sciences, Nanyang Technological University}\\
$^3$ {\normalsize\it School of Statistics and Data Science, Nankai University}
}

\begin{document}

\maketitle

\begin{abstract}
\baselineskip 20pt
Existing methods for high-dimensional changepoint detection and localization typically focus on changes in either the mean vector or the covariance matrix separately. This separation reduces detection power and localization accuracy when both parameters change simultaneously. We propose a simple yet powerful method that jointly monitors shifts in both the mean and covariance structures. Under mild conditions, the test statistics for detecting these shifts jointly converge in distribution to a bivariate standard normal distribution, revealing their asymptotic independence. This independence enables the combination of the individual p-values using Fisher's method, and the development of an adaptive p-value-based estimator for the changepoint. Theoretical analysis and extensive simulations demonstrate the superior performance of our method in terms of both detection power and localization accuracy.
\end{abstract}


\section{Introduction}\label{section_introduction}

High-dimensional data, arising from fields such as finance, genomics, and industrial engineering, often exhibit structural breaks, or changepoints. Detecting these changes, particularly in the mean vector and covariance matrix, is crucial for accurate statistical inference and reliable downstream analyses. Formally, let $\{\X_i=(x_{i1},\ldots,x_{ip})^{\top}\}_{i=1}^{n}$ be independent observations from $\bR^p$, with means $\bm\mu_{i}=\E(\X_i)$ and covariances $\bms_{i}=\cov(\X_{i})$. The goal is to test whether there is a changepoint in the mean or covariance, i.e.,
\begin{gather*}
	\bH_{0}:\   (\bm\mu_{1},\bms_1)=\cdots=(\bm\mu_{n},\bms_n) \quad \text{versus}\\
	\bH_{1}:\ \exists\ \tau^{\ast}\in\{1,\ldots,n-1\}, \ (\bm\mu_{1},\bms_1)=\cdots=(\bm\mu_{\tau^{\ast}},\bms_{\tau^{\ast}}) \neq(\bm\mu_{\tau^{\ast}+1},\bms_{\tau^{\ast}+1})=\cdots=(\bm\mu_{n},\bms_n).
\end{gather*}

In low-dimensional settings, likelihood-based methods efficiently address this problem \citep{csorgo1997limit, chen2012parametric}. However, for high-dimensional scenarios, current approaches generally handle mean or covariance shifts separately, assuming the other parameter remains unchanged \citep{aue2009break,zhang2010detecting,horvath2012change,jirak2015uniform,wang2018change,avanesov2018change,enikeeva2019high,wang2019multiple,zhong2019homogeneity,liu2020unified,liu2021minimax,yu2021finite,wang2022inference,zhang2022adaptive,wang2023computationally,li2024inference,dornemann2024detecting}. Recent methods have attempted to address both shifts indirectly---for instance, by vectorizing covariance matrices to convert covariance shifts into mean shifts \citep{wang2022inference}, or by employing a unified U-statistic framework that accommodates multiple parameter types \citep{liu2020unified}. Nevertheless, these methods still treat mean and covariance shifts separately, rather than jointly. To our knowledge, no existing method explicitly detects changepoints in both the mean and covariance simultaneously in high dimensions.

Nonparametric tests for general distributional shifts \citep[e.g.,][]{matteson2014nonparametric, chu2019asymptotic, londschien2023random} detect departures beyond the first two moments but may be less powerful for mean-covariance changes, since they are not designed for that specific structure.

\subsection{Limitations of Existing Approaches}\label{subsection_powerloss}

A straightforward approach to simultaneous high-dimensional changepoint detection might apply separate tests designed for mean or covariance shifts. However, this approach has two primary drawbacks. First, covariance-specific tests exhibit limited sensitivity to mean shifts, and vice versa. Second, when both parameters change, focusing on just one parameter reduces test power due to interference from shifts in the other parameter. A natural remedy is to combine separate tests via a Bonferroni correction, which controls Type-I error rates but is often overly conservative, sacrificing detection power.

\begin{figure}[htbp]
	\centering{
		\includegraphics[]{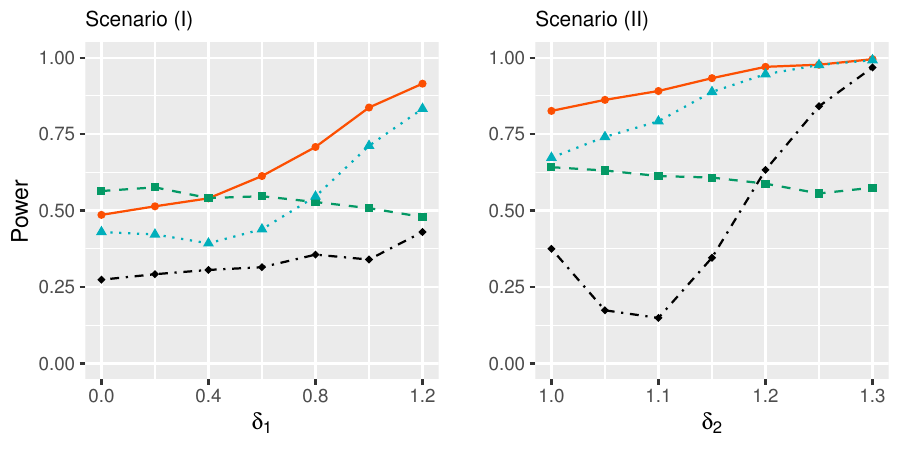}
	}
	\caption{\label{fig-intro}Empirical power comparison under scenarios (I) and (II). Shown are the individual mean or covariance changepoint test from \cite{wang2022inference} (dashed line with square markers), their Bonferroni combination (dotted line with triangle markers), the nonparametric test from \cite{chu2019asymptotic} (dash-dot line with diamond markers), and the proposed test (solid line with circle markers).}
\end{figure}

Consider a simple scenario with $n=200$, $p=100$, and a changepoint $\tau^{\ast}=100$. Samples before and after the changepoint are generated from distributions $N(\bm 0,\bms_1)$ and $N(\bm\mu,\bms_2)$, respectively, where $\bm\mu=\delta_1(p^{-1/2},\ldots,p^{-1/2})^{\top}$, $\bms_1=(0.3^{|i-j|})$, and $\bms_2=\delta_2(0.5^{|i-j|})$. We explore two scenarios: (I) varying $\delta_1$ while holding $\delta_2$ fixed, and (II) fixing $\delta_1$ while varying $\delta_2$. We compare several methods, including individual mean and covariance changepoint tests proposed by \cite{wang2022inference}, their Bonferroni combination, and a graph-based nonparametric changepoint test from \cite{chu2019asymptotic}. All methods preserve the nominal 5\% Type-I error (with the Bonferroni combination being slightly conservative), which is consistent with the corresponding theoretical results; see Table~\ref{table-S-H0} in Section~\ref{section_simulation}. Figure~\ref{fig-intro} presents empirical power comparisons, showing that the nonparametric test loses power in cases of modest change. Moreover, the individual mean or covariance test exhibits reduced sensitivity when the complementary parameter also shifts. While the Bonferroni combination improves upon individual tests, there remains considerable room for enhancing detection power.

\subsection{Contributions of this Paper}

This paper develops a simple yet effective procedure for simultaneously detecting and localizing the changepoint in both the mean vector and covariance matrix. Our contributions include:
\begin{itemize}
	\item \textit{Joint test statistic}: We construct separate $L_2$-norm-based test statistics for mean and covariance shifts, which, under mild assumptions, jointly converge in distribution to a bivariate standard normal distribution under the null hypothesis. This reveals their asymptotic independence, allowing us to combine them efficiently through Fisher's method \citep{littell1971asymptotic}.
	\item \textit{Adaptive estimator}: We propose an adaptive p-value-based estimator that adjusts to the type of data shifts---whether in the mean, covariance, or both.
\end{itemize}
We establish rigorous theoretical properties, demonstrating the consistency of both our joint testing and adaptive estimation procedures. Extensive simulation studies illustrate that our approach significantly improves detection power and localization accuracy compared to existing methods. In particular, Figure~\ref{fig-intro} illustrates the superior test power of our method across different scenarios.

The remainder of this paper is organized as follows. Section~\ref{section_test} introduces the proposed procedure for the simultaneous detection of changes in the mean and covariance, and establishes its asymptotic properties. Section~\ref{section_estimate} presents the adaptive change-point estimator and establishes its consistency. Numerical experiments including a real data example are conducted in Section~\ref{section_simulation}. Further discussions are given in Section~\ref{sec_conclusion}. All theoretical proofs are relegated to the Supplementary Material.

\textit{Notations}.
For a scalar $a$, we denote its integer part by $\lfloor a\rfloor$.
For a vector $\bm g$, $\|\bm g\|$ is its $L_2$ norm.
For a matrix $\bm A$, $\tr(\bm A)$ denotes its trace.
For random variables $z_1,\ldots, z_n$, $\cum(z_1,\ldots,z_n)$ represents their joint cumulant.
For two sequences $\{a_n\}$ and $\{b_n\}$, we write $a_n\sim b_n$ to indicate that they are asymptotically of the same order.
For nonnegative integers $k\leq n$, define $P_n^k=n!/(n-k)!$ as the number of $k$-permutations of $n$. 
The operator $\sum_{1\leq i_1,\ldots,i_m\leq n}^{\ast}$ denotes summation over mutually distinct indices $i_1,\ldots,i_m$.
Finally, $\Phi(\cdot)$ denotes the cumulative distribution function of a standard normal random variable.

\section{Simultaneous Changepoint Detection}\label{section_test}

\subsection{Proposed Test via Asymptotic Independence}

To jointly detect changes in both the mean and covariance, we develop two tailored test statistics---one for each parameter---and then fuse the information using a p-value combination approach justified by asymptotic theory.

If a potential changepoint were known, the problem would reduce to a two-sample comparison. For any candidate changepoint $\tau\in\{4,\ldots,n-4\}$, we divide the observations into two groups: $\{\X_1,\ldots,\X_{\tau}\}$ and $\{\X_{\tau+1},\ldots,\X_{n}\}$. To test whether the means of these two groups are equal, a natural statistic is the squared $L_2$-distance between group means $\bar{\X}_{\tau-}$ and $\bar{\X}_{\tau+}$, i.e., $\|\bar{\X}_{\tau-}-\bar{\X}_{\tau+}\|^2$. We employ a U-statistic-based variant \citep{chen2010two}:
\begin{equation*}
	M_n(\tau)=\frac{1}{\tau(\tau-1)(n-\tau)(n-\tau-1)}\sum_{1\leq i_1,i_2\leq \tau}^{\ast}\sum_{\tau+1\leq j_1,j_2\leq n}^{\ast}(\X_{i_1}-\X_{j_1})^{\top}(\X_{i_2}-\X_{j_2}).
\end{equation*}
Note that $M_n(\tau)$ can be expressed as
\begin{equation*}
	M_n(\tau)=\|\bar{\X}_{\tau-}-\bar{\X}_{\tau+}\|^2-\tau^{-1}\tr(\bm S_{\tau-})-(n-\tau)^{-1}\tr(\bm S_{\tau+}),
\end{equation*}
where $\bm S_{\tau-}$ and $\bm S_{\tau+}$ are the group covariance matrices,
\begin{equation*}
	\bm S_{\tau-}=\frac{\sum_{1\leq i_1,i_2\leq\tau}(\X_{i_1}-\X_{i_2})(\X_{i_1}-\X_{i_2})^{\top}}{2\tau(\tau-1)},  \quad
	\bm S_{\tau+}=\frac{\sum_{\tau+1\leq i_1,i_2\le n}(\X_{i_1}-\X_{i_2})(\X_{i_1}-\X_{i_2})^{\top}}{2(n-\tau)(n-\tau-1)},
\end{equation*}
respectively. Therefore, $M_n(\tau)$ isolates the mean shift by removing within‑group variance terms. Under the null hypothesis of no change, $\E\{M_n(\tau)\}=0$ for all $\tau$, while under the alternative, $\E\{M_n(\tau^{\ast})\}=\|\bm\mu_1-\bm\mu_n\|^2$. Since the actual changepoint $\tau^{\ast}$ is unknown, we aggregate information across all candidate changepoints. There are two common aggregation strategies in the literature: one based on the maximum of $M_n(\tau)$ or its variants \citep[e.g.,][]{wang2022inference}, and another based on the summation \citep{wang2018change,wang2019multiple,wang2023computationally}. We adopt the latter approach due to its favorable asymptotic normality. Specifically, we define the aggregated mean-shift statistic as:
\begin{equation*}
	M_n=\sum_{\tau=2}^{n-2}\left\{\frac{\tau(n-\tau)}{n}M_n(\tau)\right\}.
\end{equation*}

For covariance test, a direct measure of covariance change at $\tau$ is the squared Frobenius norm of the difference of group covariances $\bm S_{\tau-}$ and $\bm S_{\tau+}$, i.e., $\|\bm S_{\tau-}-\bm S_{\tau+}\|^2_F$.
We adopt a U-statistic version \citep{li2012two}:
\begin{align*}
	V_n(\tau)=\frac{1}{P_{\tau}^4}&\sum_{1\leq i,j,k,l\leq \tau}^{\ast}H(\X_i,\X_j,\X_k,\X_l)+\frac{1}{P_{n-\tau}^4}\sum_{\tau+1\leq i,j,k,l\leq n}^{\ast}H(\X_i,\X_j,\X_k,\X_l)\\
	&-\frac{2}{P_{\tau}^2P_{n-\tau}^2}\sum_{1\leq i,j\leq \tau}^{\ast}\sum_{\tau+1\leq k,l\leq n}^{\ast}H(\X_i,\X_j,\X_k,\X_l),
\end{align*}
where $H(\X_i,\X_j,\X_k,\X_l)=\{(\X_i-\X_j)^{\top}(\X_k-\X_l)\}^2/4$. Similarly, $V_n(\tau)$ can also be regarded as the squared Frobenius norm difference between the group covariance matrices, adjusted by subtracting terms that are irrelevant to covariance test. Under the null hypothesis of no change,
$\E\{H(\X_i,\X_j,\X_k,\X_l)\}=\tr(\bms_1)$ for all distinct $i,j,k,l$. Consequently, $\E\{V_n(\tau)\} = 0$ for all $\tau$. Under the alternative, we have
\begin{align*}
	\frac{1}{P_{\tau^{\ast}}^4}\sum_{1\leq i,j,k,l\leq \tau^{\ast}}^{\ast}\E\{H(\X_i,\X_j,\X_k,\X_l)\} &= \tr(\bms_1^2);\\
	\frac{1}{P_{n-\tau^{\ast}}^4}\sum_{\tau^{
			\ast
		}+1\leq i,j,k,l\leq n}^{\ast}\E\{H(\X_i,\X_j,\X_k,\X_l)\} &= \tr(\bms_n^2);\\
	\frac{1}{P_{\tau^{\ast}}^2P_{n-\tau^{\ast}}^2}\sum_{1\leq i,j\leq \tau^{\ast}}^{\ast}\sum_{\tau^{\ast}+1\leq k, l\leq n}^{\ast}\E\{H(\X_i,\X_j,\X_k,\X_l)\} &= \tr(\bms_1\bms_n).
\end{align*}
Hence, $\E\{V_n(\tau^{\ast})\}=\tr\{(\bms_1-\bms_n)^2\}$. Similar to the mean statistic, we define the aggregated covariance-shift statistic as:
\begin{equation*}
	V_n=\sum_{\tau=4}^{n-4}\left\{\frac{\tau(n-\tau)}{n}V_n(\tau)\right\}.
\end{equation*}

\begin{remark}
	Although direct computations of $M_n(\tau)$ and $V_n(\tau)$ may be computationally intensive, efficient reformulations can significantly reduce computational costs (see the Supplementary Material). Specifically, calculating $M_n$ and $V_n$ take $O(np)$ and $O(n^2p)$ operations, respectively.
\end{remark}

To jointly assess changes in the mean and covariance, we examine the joint asymptotic behavior of $M_n$ and $V_n$. Define $\Y_i=\X_i-\bm\mu_i$, $i=1,\ldots,n$, and consider the following moment and cumulant assumptions.

\begin{assumption}\label{cond_sigma}
	$\tr(\bms_i^4)=o\{\tr^2(\bms_i^2)\}$.
\end{assumption}

\begin{assumption}\label{cond_cum_mean}
	There exists a constant $C$ independent of $n, p$ such that
	\begin{equation*}
		\sum_{l_1,\ldots,\l_h}^{p}\cum^2(y_{il_1},\ldots,y_{il_h})\leq C\tr^{h/2}(\bms_i^2),
	\end{equation*}
	for $h=2,3,4$.
\end{assumption}

\begin{assumption}\label{cond_cum_var}
	There exists a constant $C$ independent of $n, p$ such that
	\begin{equation*}
		\sum_{l_1,\ldots,l_h}^{p}\cum^2(y_{il_1},\ldots,y_{il_h}) \leq C\tr^{h/2}(\bms_i^2),
	\end{equation*}
	for $h=5,\ldots,8$. Further,
	\begin{equation*}
		\sum_{l_1,\ldots,l_4}^{p}\cum^2(y_{il_1},\ldots,y_{il_4})=o\{\tr^{2}(\bms_i^2)\}.
	\end{equation*}
\end{assumption}

\begin{assumption}\label{cond_cum_independent}
	\begin{equation*}
		\sum_{l_1,\ldots,l_3}^{p}\cum^2(y_{il_1},\ldots,y_{il_3})=o\{\tr^{3/2}(\bms_i^2)\}.
	\end{equation*}
\end{assumption}

Assumption~\ref{cond_sigma} is a standard moment condition commonly used in high-dimensional two-sample and changepoint testing problems for both mean and covariance structures. Assumption~\ref{cond_cum_mean} controls the joint cumulants and is critical for establishing the asymptotic normality of the mean statistic $M_n$. Similar conditions have been utilized in previous works, such as \cite{wang2022inference}, for high-dimensional mean changepoint detection; our adoption of a sum-type aggregation allows for slightly weaker conditions than those required by max-type aggregation approaches. Assumption~\ref{cond_cum_var} plays a parallel role in ensuring the asymptotic normality of the covariance statistic $V_n$ by controlling higher-order cumulants. Finally, Assumption~\ref{cond_cum_independent} guarantees the asymptotic independence between $M_n$ and $V_n$, which underpins the validity of our joint testing procedure.

\begin{remark}
	If $\Y_i$ follows a $p$-dimensional multivariate normal distribution, all cumulants of order $h\geq 3$ vanish, and $\sum_{l_1,l_2=1}^{p}\cum^2(y_{il_1},y_{il_2})=\tr^2(\bms_i^2)$.  Thus, Assumptions \ref{cond_cum_mean}--\ref{cond_cum_independent} are automatically satisfied as $p\rightarrow \infty$. These assumptions also hold for pseudo-factor models with finite moments up to order $8$ \citep{chen2010two, li2012two}.
\end{remark}

\begin{theorem}\label{thm_independent}
	Under Assumptions~\ref{cond_sigma}--\ref{cond_cum_independent} and the null hypothesis, as $n,p\rightarrow \infty$,
	\begin{equation*}
		\Pr\left(\frac{M_n}{\sqrt{\var_0(M_n)}}\leq x_1, \frac{V_n}{\sqrt{\var_0(V_n)}}\leq x_2\right)\longrightarrow \Phi(x_1)\Phi(x_2),
	\end{equation*}
	for any $x_1,x_2\in \bR$, where 
	\begin{equation*}
		\var_0(M_n)=\frac{2\pi^2-18}{3} n^2\tr(\bms_1^2)\{1+o(1)\}, \quad  \var_0(V_n)=\frac{4\pi^2-36}{3} n^2\tr^2(\bms_1^2)\{1+o(1)\}.
	\end{equation*}
\end{theorem}

\begin{remark}
	In high-dimensional two-sample testing, \cite{yu2023power} showed that the mean statistic proposed by \cite{chen2010two} and the covariance statistic by \cite{li2012two} jointly converge in distribution to a standard bivariate normal under the null hypothesis of equal means and covariances. We extend this result to the changepoint testing framework, which is nontrivial due to the dependence introduced by aggregating correlated two-sample test statistics across multiple candidate changepoints.
\end{remark}

Both $\var_0(M_n)$ and $\var_0(V_n)$ depend on the unknown quantity $\tr(\bms_1^2)$. To address this, we propose the following difference-based estimator:
\begin{equation*}
	\widehat{\tr(\bms_1^2)}=\frac{1}{4(n-3)}\sum_{i=1}^{n-3} \left\{(\X_{i}-\X_{i+1})^{\top}(\X_{i+2}-\X_{i+3})\right\}^2,
\end{equation*}
which is unbiased for $\tr(\bms_1^2)$ under the null hypothesis $H_{0}$. By replacing $\tr(\bms_1^2)$ in $\var_0(M_n)$ and $\var_0(V_n)$ with $\widehat{\tr(\bms_1^2)}$, we obtain the plug-in estimators, denoted as $\widehat{\sigma}_{1}^2$ and $\widehat{\sigma}_{2}^2$, respectively.

\begin{proposition}\label{pro_sigma_estimator}
	Under Assumptions~\ref{cond_sigma}--\ref{cond_cum_mean} and the null hypothesis, as $n,p \rightarrow \infty$,
	\begin{equation*}
		\{\tr(\bms_1^2)\}^{-1}\widehat{\tr(\bms_1^2)}\rightarrow 1
	\end{equation*}
	 in probability. Moreover, Theorem~\ref{thm_independent} remains valid when $\var_0(M_n)$ and $\var_0(V_n)$ are replaced by their corresponding plug-in estimators $\widehat{\sigma}_{1}^2$ and $\widehat{\sigma}_{2}^2$.
\end{proposition}

According to Proposition~\ref{pro_sigma_estimator}, the scaled statistics $\widehat{\sigma}_{1}^{-1}M_n$ and $\widehat{\sigma}_{2}^{-1}V_n$ are asymptotically independent. To jointly test for a changepoint in both the mean and covariance, we combine their respective p-values via Fisher's method. Specifically, define $p_{M_n}=1-\Phi(\widehat{\sigma}_{1}^{-1}M_n)$ and $p_{V_n}=1-\Phi(\widehat{\sigma}_{2}^{-1}V_n)$. The combined test statistic is then given by
\begin{equation*}
	T_{n}=-2\log(p_{M_n})-2\log(p_{V_n}).
\end{equation*}
Under the null hypothesis, $T_n$ converges in distribution to a chi-squared distribution with four degrees of freedom. Therefore, for a significance level $\alpha\in(0,1)$, we reject the null hypothesis if $T_n$ exceeds the upper-$\alpha$ quantile of this limiting distribution.

\subsection{Consistency of the Proposed Test}

To establish the consistency of our proposed test, we introduce the following standard assumption on the changepoint location, commonly used in high-dimensional changepoint literature \citep[e.g.,][]{liu2020unified,wang2022inference,zhang2022adaptive}.

\begin{assumption}\label{cond_length}
	$\tau^{\ast}=\lfloor \rho n \rfloor$ for some constant $0<\rho <1$.
\end{assumption}

\begin{theorem}\label{thm_power}
	Under Assumptions~\ref{cond_sigma}--\ref{cond_cum_var} and \ref{cond_length}, as $n,p \rightarrow \infty$, the following results hold:
	\begin{itemize}
		\item[](i) If $n\|\bm\mu_1-\bm\mu_n\|^2/\sqrt{\tr\{(\bms_1+\bms_n)^2\}}\rightarrow \infty$, then $p_{M_n}\rightarrow 0$ and $T_n\rightarrow \infty$ in probability.
		\item[](ii) If $n\|\bm\mu_1-\bm\mu_n\|^2/\sqrt{\tr\{(\bms_1+\bms_n)^2\}}=O(1)$, $n\tr\{(\bms_1-\bms_n)^2\}/\tr\{(\bms_1+\bms_n)^2\}\rightarrow \infty$, then $p_{V_n}\rightarrow 0$ and $T_n\rightarrow \infty$ in probability. 
	\end{itemize}
\end{theorem}

\begin{remark}
	The conditions in $(i)$ and $(ii)$ are analogous to those in \cite{yu2023power} for simultaneous testing of high-dimensional mean and covariance in a two-sample setting. Furthermore, when only the mean changes (i.e., $\bms_1=\bms_n$), the condition in (i) matches the optimal detection rate (up to a logarithmic factor) for dense mean changes in high dimensional data \citep{liu2021minimax, wang2022inference, zhang2022adaptive}. When only the covariance changes (i.e., $\bm\mu_1=\bm\mu_n$), the condition in (ii) aligns with the change magnitude required in \cite{dornemann2024detecting} for consistent test of high-dimensional covariance changes.
\end{remark}

\section{Changepoint Localization}\label{section_estimate}

\subsection{Limitations in Localization Accuracy of Existing Methods}\label{subsetion_accuracyloss}

Changepoint estimation, like testing, can suffer from reduced accuracy when methods designed for identifying changes in a single parameter (mean or covariance) are applied to cases where both parameters shift simultaneously. For example, applying a mean changepoint localization method in the presence of a covariance shift can result in a notable degradation in localization accuracy, and vice versa.

\begin{figure}[htbp]
	\centering{
		\includegraphics[]{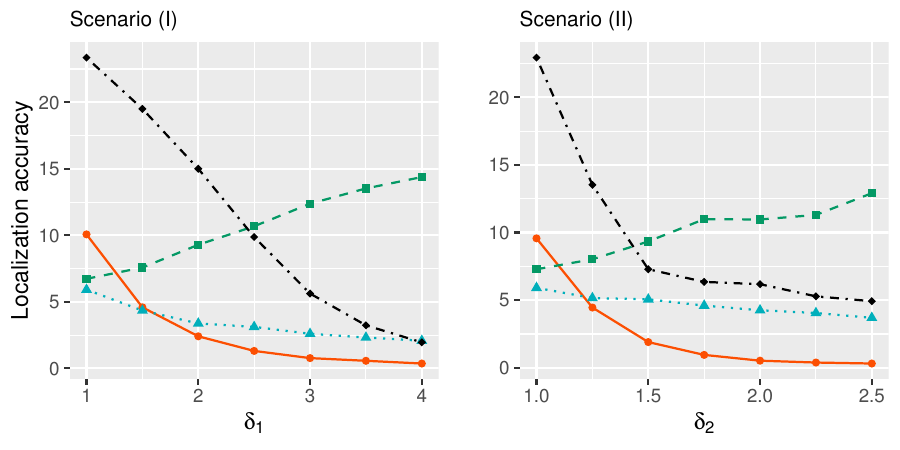}
	}
	\caption{\label{fig-estimation}Empirical localization accuracy comparison under scenarios (I) and (II). Shown are the individual mean or covariance changepoint localization method from \cite{wang2022inference} (dashed line with square markers), their minimal p-value combination procedure (dotted line with triangle markers), the nonparametric estimator from \cite{chu2019asymptotic} (dash-dot line with diamond markers), and the proposed estimator (solid line with circle markers).}
\end{figure}

To illustrate these limitations, we revisit the simulated example from Section~\ref{subsection_powerloss} with the same setup. We compare individual localization procedures for mean and covariance changes, as proposed by \citep{wang2022inference}. A natural strategy to accommodate different change patterns---whether in the mean or covariance---is to select the estimator corresponding to the test with the smaller p-value. This minimal p-value strategy has been previously employed for sparsity-adaptive mean changepoint estimation in high-dimensional settings \citep{liu2020unified,zhang2022adaptive,wang2023computationally}. We apply this strategy to the two individual tests. In addition, we include the graph-based nonparametric changepoint estimator from \cite{chu2019asymptotic} for comparison. For fairness, all methods restrict the grid search interval to $[\lambda_n,n-\lambda_n]$ with $\lambda_n=\lfloor 0.2n \rfloor$.

Figure~\ref{fig-estimation} presents a comparison of empirical localization accuracy, measured as the absolute distance between the estimated and true changepoint. The results show that the nonparametric approach loses accuracy for moderate changes. Additionally, the individual mean changepoint method is heavily influenced by covariance shifts (and vice versa). Although the minimal p-value strategy alleviates this issue, it still relies on applying a single estimator designed for either mean or covariance changes. This leaves substantial room for improving localization accuracy when both parameters shift simultaneously.

\subsection{Proposed Adaptive Changepoint Estimator}

Rather than selecting or combining individual changepoint estimators, our method integrates both the mean and covariance changes into a unified monitoring statistic.

For each candidate changepoint $\tau$, define
\begin{equation*}
	\tilde{M}_{n}(\tau)=\left\{2\widehat{\tr(\bms_1^2)}\right\}^{-1/2}\frac{\tau(n-\tau)}{n}M_n(\tau), \quad
	\tilde{V}_{n}(\tau)=\left\{2\widehat{\tr(\bms_1^2)}\right\}^{-1}\frac{\tau(n-\tau)}{n}V_n(\tau).
\end{equation*}
Under the null hypothesis, and given Assumptions~\ref{cond_sigma}--\ref{cond_cum_independent}, it can be shown that for each $\tau\sim n$, $\tilde{M}_{n}(\tau)$ and $\tilde{V}_{n}(\tau)$ jointly converge in distribution to a bivariate standard normal distribution, as $n,p \rightarrow \infty$.

To estimate the changepoint, we then compute the combined test statistic for each candidate location $\tau$:
\begin{equation*}
	T_n(\tau)=-2\log(p_{\tilde{M}_n(\tau)})-2\log(p_{\tilde{V}_n(\tau)}),
\end{equation*}
where $p_{\tilde{M}_n(\tau)}=1-\Phi(\tilde{M}_{n}(\tau))$ and $p_{\tilde{V}_n(\tau)}=1-\Phi(\tilde{V}_{n}(\tau))$. The adaptive changepoint estimator is then defined as:
\begin{equation*}
	\widehat{\tau}=\underset{\lambda_n\leq \tau\leq n-\lambda_n}{\arg\max}T_n(\tau),
\end{equation*}
where $\lambda_n = \lfloor\lambda n\rfloor$ for a constant $\lambda\in (0,0.5)$.

\begin{theorem}\label{thm_estimation_consistency}
	Suppose that Assumptions~\ref{cond_sigma}--\ref{cond_cum_var} and Assumption~\ref{cond_length} hold. Further assume that one of the following conditions is satisfied:
	\begin{itemize}
		\item[](i) $n\|\bm\mu_1-\bm\mu_n\|^2/\sqrt{\tr\{(\bms_1+\bms_n)^2\}}\rightarrow \infty$, $\|\bm\mu_1-\bm\mu_n\|^2/\sqrt{\tr\{(\bms_1+\bms_n)^2\}}=o(1)$, and
		$n\tr\{(\bms_1-\bms_n)^2\}/\tr\{(\bms_1+\bms_n)^2\}=O(1)$;
		\item[](ii) $n\|\bm\mu_1-\bm\mu_n\|^2/\sqrt{\tr\{(\bms_1+\bms_n)^2\}}=O(1)$  and $n\tr\{(\bms_1-\bms_n)^2\}/\tr\{(\bms_1+\bms_n)^2\}\rightarrow \infty$;
		\item[](iii) $n\|\bm\mu_1-\bm\mu_n\|^2/\sqrt{\tr\{(\bms_1+\bms_n)^2\}}\rightarrow \infty$,  $\|\bm\mu_1-\bm\mu_n\|^2/\sqrt{\tr\{(\bms_1+\bms_n)^2\}}=o(1)$, and $n\tr\{(\bms_1-\bms_n)^2\}/\tr\{(\bms_1+\bms_n)^2\}\rightarrow \infty$.
	\end{itemize}
	Then, $\widehat{\tau}$ is a consistent estimator of $\tau^{\ast}$, that is, $(\widehat{\tau}-\tau^{\ast})/n=o_p(1)$. 
\end{theorem}

\begin{remark}
	If only the mean changes, the condition $n\|\bm\mu_1-\bm\mu_n\|^2/\sqrt{\tr\{(\bms_1+\bms_n)^2\}}\rightarrow \infty$ aligns with the change magnitude required for high-dimensional mean changepoint estimation as proposed in \cite{zhang2022adaptive}. In addition, the condition $n\tr\{(\bms_1-\bms_n)^2\}/\tr\{(\bms_1+\bms_n)^2\}\rightarrow \infty$ matches the change magnitude outlined in \cite{dornemann2024detecting} for high-dimensional covariance changepoint estimation.
\end{remark}

The empirical results, such as those presented in Figure~\ref{fig-estimation}, further demonstrate that our proposed estimator generally outperforms existing methods in terms of localization accuracy.

\section{Numerical Studies}\label{section_simulation}

\subsection{Simulation}

In this section, we further examine the finite-sample performance of the proposed method and various competitors across various simulation settings.

For notational convenience, when there is no changepoint, we define $\tau^{\ast}=n$. We generate the data according to the model $\X_i=\bm\mu_i+\bm\varepsilon_i$, $i=1,\ldots,n$. For the mean vector $\bm\mu_i$, we set $\bm\mu_i=(0,\ldots,0)^{\top}$ for $i\leq \tau^{\ast}$; and $\bm\mu_{i}=\delta_1(p^{-1/2},\ldots,p^{-1/2})^{\top}$ for $i> \tau^{\ast}$. For the random error vector $\bm\varepsilon_i$, given the covariance matrix $\bms_i$, we generate $\bm\varepsilon_i$ as $\bm\varepsilon_i= \bms_{i}^{1/2}\bm e_i$, where $\bm e_i=(e_{i1},\ldots,e_{ip})^{\top}$, and $e_{i1},\ldots,e_{ip}$ are independently identically distributed. We consider the standard normal distribution $N(0,1)$ and the $t(9)$ distribution ($t$-distribution with nine degrees of freedom) standardized with unit variance.

For the covariance structures, we consider two scenarios:
\begin{itemize}
	\item[] (I) AR(1) structure: $\bms_1=(0.3^{\vert i-j\vert})$ and $\bms_{\tau^{\ast}+1}=\delta_2\cdot(0.5^{\vert i-j\vert})$.
	\item[] (II) Blocked diagonal structure: $\bms_1=\text{blk}(5,0.3)$ and $\bms_{\tau^{\ast}+1}=\delta_2\cdot\text{blk}(5,0.5)$, where 
	$\text{blk}(5,r)$ denotes a matrix $\bms$ such that $\Sigma_{jj}=1$ for $j=1,\ldots,p$; $\Sigma_{ij}=r$ for $k = 1,\ldots, \lfloor p/5\rfloor$ and $5(k-1)< i \neq j \leq 5k$, and $\Sigma_{ij}=0$ otherwise. 
\end{itemize}

To investigate the impact of location of the changepoint, we consider two cases $\tau^{\ast}=0.3n$ and $\tau^{\ast}=0.5n$. The search interval is set to $[\lambda_n,n-\lambda_n]$, where $\lambda_n=\lfloor 0.2n \rfloor$, for all methods to ensure a fair comparison. In all simulation experiments, the significance level is set to $\alpha=5\%$. Each empirical size, power, and localization accuracy is calculated based on $1,000$ replications.

\subsubsection{Test}\label{subsub_test}

To evaluate the performance of the proposed method in changepoint detection, we consider several competing methods: the individual mean and covariance changepoint tests introduced in \cite{wang2022inference}, their Bonferroni-type combination, and the graph-based nonparametric test proposed by \cite{chu2019asymptotic}. For brevity, we refer to these methods as Mean, Covariance, Bonferroni, and Graph, respectively. The proposed method is denoted as Ours.

We first examine whether the proposed method can effectively control the Type-I error under the specified settings. We fix $n=200$ and consider $p \in \{100, 200, 300\}$. Table \ref{table-S-H0} shows the empirical sizes of various changepoint tests under different sample dimensions and covariance structures, when $e_{ij}$'s follow the $N(0,1)$ distribution and the standardized $t(9)$ distribution, respectively. We observe that all methods effectively control the Type-I error. The Bonferroni method tends to be slightly conservative, and the other methods roughly maintain the nominal significance level, which is consistent with the corresponding theoretical results.

\begin{table}[htbp]
	\centering
	\caption{Empirical size (in \%) comparison under scenarios (I) and (II), with $p \in \{100, 200, 300\}$ and $n = 200$.}
	{	\begin{tabular}{cccccccc}
			\toprule
			& \multicolumn{3}{c}{Scenario (I)} & &\multicolumn{3}{c}{Scenario (II)} \\
			\cmidrule{2-4} \cmidrule{6-8} 
			Method & $p=100$ & $p=200$ & $p=300$ & &$p=100$ & $p=200$ & $p=300$ \\ 
			\hline
			&\multicolumn{6}{c}{Normal}\\
			Mean & 5.3 & 3.7 & 4.9 && 6.2 & 6.4 & 4.2 \\ 
			Covariance & 5.2 & 5.2 & 5.4 && 5.9 & 4.1 & 5.2 \\ 
			Bonferroni & 4.3 & 3.0 & 4.2 && 5.3 & 4.1 & 3.8 \\ 
			Graph & 4.3 & 4.9 & 5.4 && 4.7 & 5.3 & 4.9 \\ 
			Ours & 6.2 & 5.8 & 6.6 && 6.8 & 6.2 & 5.7 \\ 
			&\multicolumn{6}{c}{Standardized $t(9)$}\\
			Mean & 5.3 & 5.0 & 4.7 && 5.8 & 5.1 & 5.6 \\ 
			Covariance &4.5 & 6.5 & 5.6 && 5.5 & 4.6 & 5.3 \\ 
			Bonferroni &3.1 & 4.2 & 3.7 && 5.3 & 3.8 & 3.0 \\ 
			Graph & 4.3 & 4.6 & 6.0 && 5.2 & 4.9 & 4.1 \\ 
			Ours &  5.9 & 5.5 & 7.4 && 6.5 & 5.2 & 4.2 \\ 
			\bottomrule
	\end{tabular}}
	\label{table-S-H0}
\end{table}

\begin{figure}[htbp]
	\centering{
		\includegraphics[]{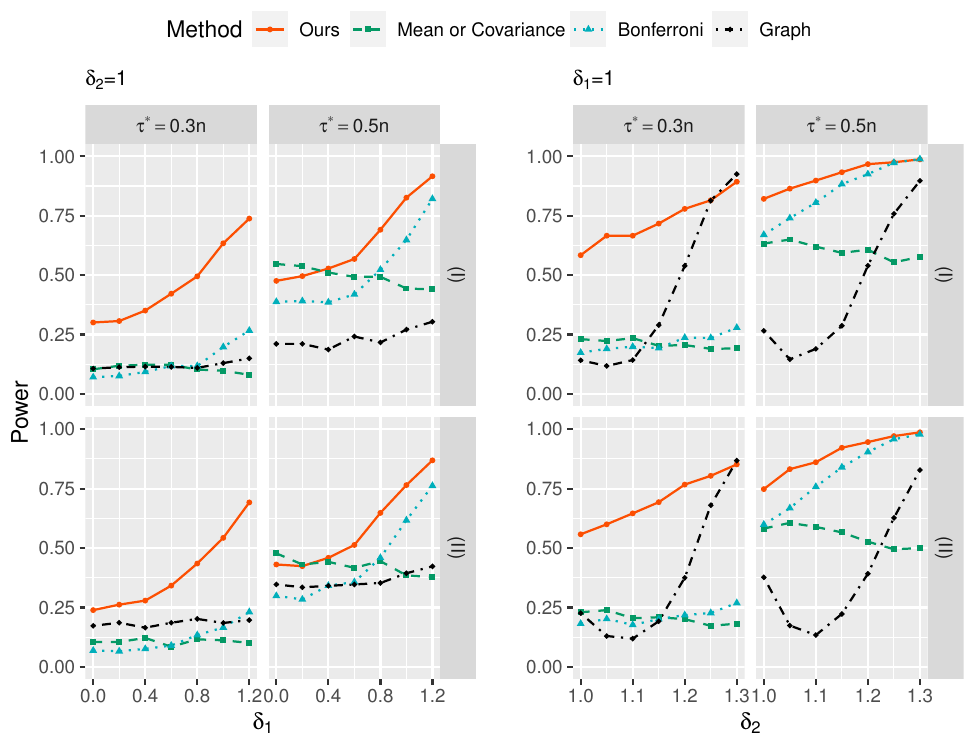}
	}
	\caption{\label{fig-S-power}Empirical power comparison under scenarios (I) and (II), with $\tau^{\ast} \in \{0.3n, 0.5n\}$, where $n = 200$, $p = 100$, and the $e_{ij}$'s follow the standardized $t(9)$ distribution.
		The method represented by the dashed line with square markers corresponds to the Covariance method in the left panel, and to the Mean method in the right panel.}
\end{figure}

Next, we examine the power performance of the proposed method under scenarios (I) and (II). We fix $n = 200$, $p = 100$, and vary $\tau^{\ast}\in \{ 0.3n, 0.5n\}$. To investigate the impact of mean and covariance change magnitudes on detection power, we consider two settings: (i) fixing the covariance change size ($\delta_2 = 1$) while varying the mean change size by adjusting $\delta_1$; and (ii) fixing the mean change size ($\delta_1 = 1$) while varying the covariance change size by adjusting $\delta_2$. Since the results under scenario (I) with $\tau^{\ast}=0.5n$ have already been presented in Figure~\ref{fig-intro} in Section~\ref{subsection_powerloss} for the case where the $e_{ij}$'s follow the $N(0,1)$ distribution, we report here only the empirical power when the $e_{ij}$'s follow the standardized $t(9)$ distribution. The corresponding empirical powers are presented in the left and right panels of Figure~\ref{fig-S-power}, respectively. When $\tau^{\ast} = 0.5n$, the power of the Mean method decreases as the magnitude of the covariance change increases, and vice versa. The Bonferroni method improves to some extent. The Graph method underperforms the Bonferroni method. In contrast, our method outperforms the Bonferroni method.
When the changepoint is close to the boundary, i.e., $\tau^{\ast} = 0.3n$, the effective sample size decreases, leading to reduced power across all methods compared to the case where $\tau^{\ast} = 0.5n$. Despite this, our method demonstrates significantly higher power than all competing methods.

\subsubsection{Estimation}

Finally, we evaluate the localization accuracy of the proposed method, measured as the absolute distance between the estimated and true changepoint. For comparison, we consider several benchmark methods: the individual mean and covariance changepoint localization procedures from \cite{wang2022inference}, their minimal p-value combination procedure, and the graph-based nonparametric estimator proposed by \cite{chu2019asymptotic}. For brevity, we refer to these methods as Mean, Covariance, Min-p, and Graph, respectively, and denote the proposed estimator as Ours.

\begin{figure}[htbp]
	\centering{
		\includegraphics[]{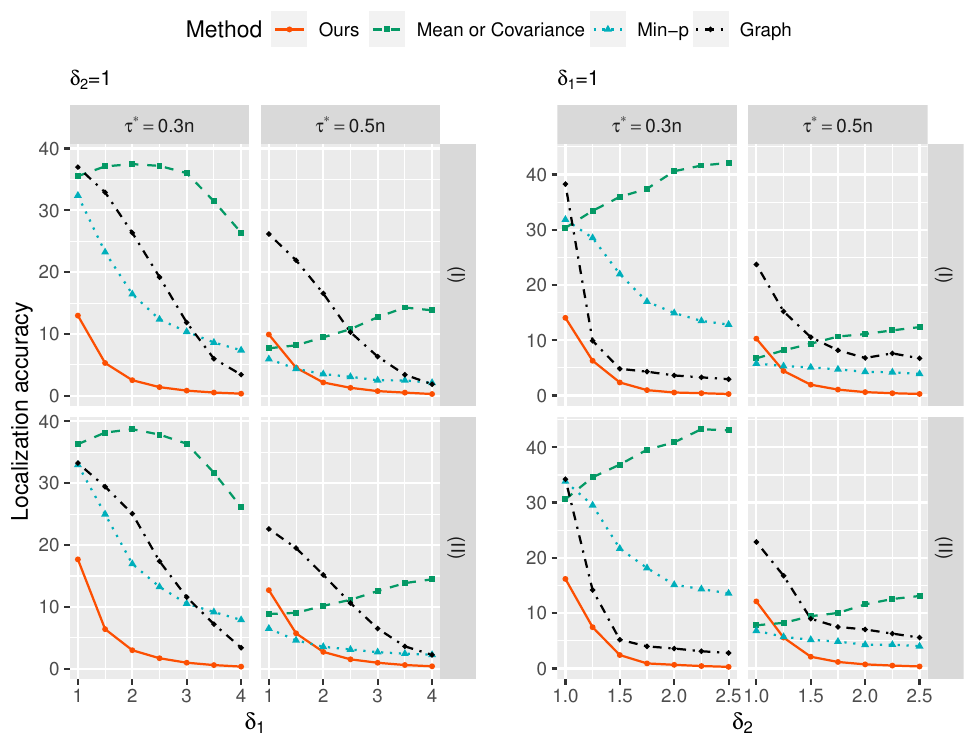}
	}
	\caption{\label{fig-S-estimation}Empirical localization accuracy comparison under scenarios (I) and (II), with changepoint locations $\tau^{\ast} \in \{0.3n, 0.5n\}$, where $n = 200$, $p = 100$, and the $e_{ij}$'s follow the standardized $t(9)$ distribution. The method represented by the dashed line with square markers corresponds to the Covariance method in the left panel, and to the Mean method in the right panel.}
\end{figure}

The simulation settings are similar to those used in the power comparison in Section~\ref{subsub_test}. Since the results for the case where the $e_{ij}$'s follow the $N(0,1)$ distribution have been partially presented in Figure~\ref{fig-estimation} in Section~\ref{subsetion_accuracyloss}, we only report the localization accuracy when the $e_{ij}$'s follow the standardized $t(9)$ distribution. The corresponding results are illustrated in Figure~\ref{fig-S-estimation}. It is evident that the localization accuracy of the Mean method is significantly affected by changes in covariance, and vice versa. The Min-p strategy alleviates this issue and improves overall performance. The Graph method achieves performance comparable to, or even better than, the Min-p method when the changepoint is close to the boundary. Among all competing methods, our method demonstrates the best estimation accuracy, particularly when the changes in both mean and covariance are not too weak.

\subsection{Real Data}

In this section, we evaluate the performance of the proposed method using a diagnostic breast cancer dataset containing $n = 569$ observations. Each instance includes $p=30$ real-valued features related to the cell nucleus, such as radius, perimeter, and smoothness, along with a categorical label indicating whether the tumor is benign or malignant. The label ``B'' denotes a benign tumor, while ``M'' denotes a malignant one. Among the 569 observations, 357 are labeled as ``B'' and the remaining 212 as ``M''. This data set is available through the ftp server in the Computer Sciences Department at UW-Madison: \texttt{http://ftp.cs.wisc.edu/math-prog/cpo-dataset/machine-learn/cancer/WDBC/}.
Based on the labels, we divide the dataset into two subsets. The mean vectors and covariance matrices of the features may differ between the two subsets. To investigate this, we apply the two-sample mean test proposed by \cite{chen2010two} and the covariance test introduced by \cite{li2012two} to compare them. The resulting p-values for both tests are extremely close to zero, indicating that the mean vectors and covariance matrices of the features differ significantly between benign and malignant tumors.

\begin{table}[ht]
	\centering	
	\caption{Estimated changepoint locations ($\tau^{\ast} = 357$).}
	{ 
		\begin{tabular}{cccccc}
			\toprule
			Mean & Covariance & Min-p & Graph & Ours \\  
			\hline    
			352 & 355 & 352 & 357 & 357 \\ 
			\bottomrule 
		\end{tabular}
	}\label{table-realdata}
\end{table}

To evaluate the effectiveness of the proposed changepoint detection method, we reorder the data by placing the 357 benign observations first, followed by the 212 malignant observations, implying that the true changepoint is at $\tau^{\ast} = 357$. Since the features are measured on different scales, we standardize each feature using its sample mean and sample standard deviation computed from the entire dataset. We then apply the Mean, Covariance, Bonferroni, Graph, and the proposed method to test for the existence of a changepoint. The p-values of all the tests are extremely close to zero, indicating strong evidence for the presence of a changepoint. Subsequently, we use the corresponding localization methods Mean, Covariance, Min-p, Graph, and our method to estimate the location of the changepoint. The results are reported in Table~\ref{table-realdata}. All the estimators are close to the true changepoint location $\tau^{\ast}$; notably, both the proposed method and the nonparametric Graph method successfully identify $\tau^{\ast}$. This provides additional empirical evidence for the effectiveness of our method.

\begin{figure}[htbp]
	\centering{
		\includegraphics[scale=0.8]{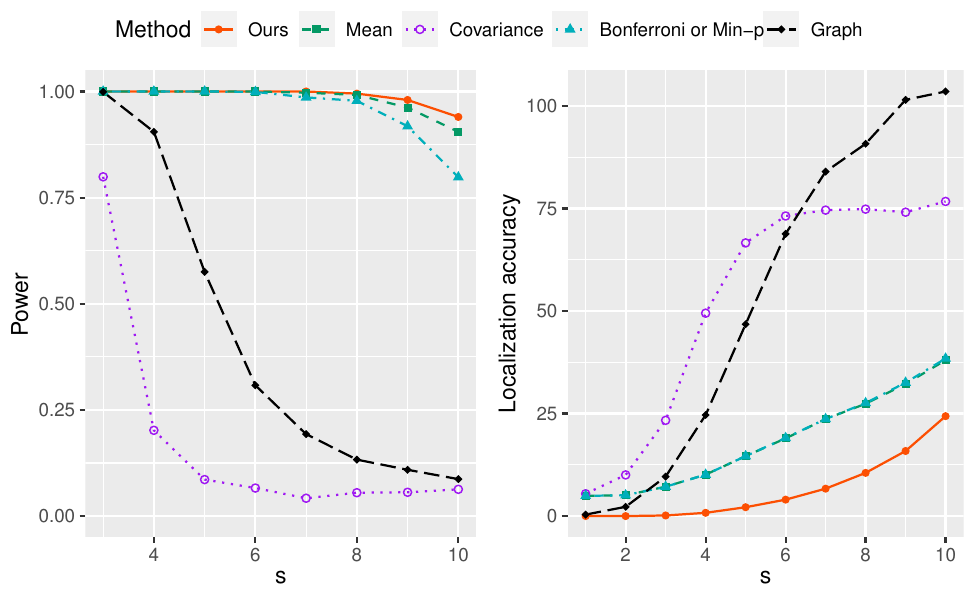}
	}
	\caption{\label{fig-realdata} Empirical power and localization accuracy comparison with additional noise added. The method represented by the dash-dot line with triangle markers corresponds to the Bonferroni method in the left panel, and to the Min-p method in the right panel.}
\end{figure}	
		
To more thoroughly evaluate the performance of the proposed method relative to competing methods, we consider a noise-perturbed setting in which independent and identically distributed noise from the $N(0,s^2)$ distribution is added to the original data. We repeated this process 1,000 times and computed the corresponding empirical power and localization accuracy. The results are presented in the left and right panels of Figure~\ref{fig-realdata}, respectively. We observe that as the noise level increases, the performance of all methods deteriorates to some extent. Among them, the Covariance and Graph methods are most severely affected by the noise. In contrast, our proposed method exhibits strong robustness to noise and achieves superior performance in both changepoint detection and estimation.

\section{Conclusion}\label{sec_conclusion}
	
In this paper, we propose a simple yet effective procedure for the simultaneous detection and localization of high-dimensional changes in both mean and covariance---a problem not explicitly studied to date. We introduce novel joint testing and adaptive estimation techniques, and rigorously establish the consistency of both components. Extensive simulations demonstrate that our method substantially outperforms existing approaches in terms of detection power and localization accuracy.
	
We conclude the article with three remarks. First, we investigate the asymptotic independence between the mean-change statistic and the covariance-change statistic under the assumption of temporally independent observations. While this form of temporal independence is standard in the mean- or covariance-changepoint literature (see, for example, \cite{enikeeva2019high}, \cite{liu2020unified}, \cite{yu2021finite}, \cite{zhang2022adaptive}, and \cite{dornemann2024detecting}), proving asymptotic independence in this setting is far from trivial. Extending our framework to more intricate settings, such as high-dimensional time series, remains an important direction for future work. Second, our test statistic and estimator rely on the $L_2$-norm aggregation, which excels at detecting dense, weak change signals characterized by many small yet nonzero entries in $\bm\mu_1 - \bm\mu_n$ and/or $\bms_1 - \bms_n$. In recent years, considerable effort has been devoted to procedures that adapt to various sparsity regimes in mean- or covariance- changepoint problems \citep[e.g.,][]{wang2018high,liu2020unified,liu2021minimax,zhang2022adaptive,hu2023likelihood,wang2023computationally}. Developing sparse-adaptive detection and localization of simultaneous mean and covariance changes in our setting is a promising line of future research. Third, an online or sequential extension is natural for streaming data: update the mean- and covariance-change components in real time with low memory, and combine them under a single alarm rule with data-driven thresholds to control false alarms and minimize delay. Establishing these guarantees and demonstrating scalability on high-throughput streams would broaden the framework's practical impact \cite{mei2010efficient, xie2013sequential,xie2020sequential}.

\section*{Supplementary Material}
	
Details on the efficient computation of the mean- and covariance-change statistics and all theoretical proofs are provided in the Supplementary Material.

\vspace{0.5cm}
{\small \baselineskip 10pt
\bibliographystyle{asa}
\bibliography{paper-ref}
}
\newpage
\setcounter{equation}{0}
\def\thelemma{S.\arabic{lemma}}
\def\thepro{S.\arabic{pro}}
\def\theequation{S.\arabic{equation}}
\def\thetable{S.\arabic{table}}
\def\thefigure{S.\arabic{figure}}
\renewcommand{\thealgocf}{S.\arabic{algocf}}
\setcounter{lemma}{0}
\setcounter{figure}{0}
\setcounter{table}{0}
\setcounter{algocf}{0}

\def\thesection{S.\arabic{section}}
\setcounter{section}{0}

\noindent{\bf\Large Supplementary Material for ``Simultaneous Detection and Localization of Mean and Covariance Changes in High Dimensions"}

\bigskip
The supplementary material contains all the theoretical proofs and details on the efficient computation of the mean- and covariance-change statistics.

\section{Proof of Main Results}

\subsection{Proof of  Theorem~1}

Recall that $M_n$ can be written as
\begin{equation*}
	M_{n}=\sum_{k=2}^{n}\sum_{i=1}^{k-1}a_{i,k}\X_{i}^{\top}\X_{k},
\end{equation*}
where $a_{i,k}=\sum_{\tau=2}^{i-1}2(n-\tau-1)^{-1}(1-n^{-1})+\sum_{\tau=k}^{n-2}2(\tau-1)^{-1}(1-n^{-1})+6n^{-1}-2$.
In the proof of Lemma~\ref{lemma_normal}, we introduce the variable $V_{n,1}$, where 
\begin{align*}
	V_{n,1}=\sum_{k=2}^{n}\sum_{i=1}^{k-1}a_{i,k}(\X_{i}^{\top}\X_{k})^2,
\end{align*}
and show that $V_n-V_{n,1}=o_p\{\var^{1/2}(V_{n,1})\}$.
To prove Theorem~1, by Slutsky’s Theorem, we only need to show that for any $x_1, x_2\in \bR$, as $n,p\rightarrow \infty$,
\begin{equation*}
	\Pr\left(\frac{M_{n}}{\sigma_{01}}\leq x_1,\frac{V_{n,1}}{\sigma_{02}}\leq x_2\right)\rightarrow \Phi(x_1)\Phi(x_2),
\end{equation*}
where $\sigma_{01}=\var^{1/2}(M_{n})$, $\sigma_{02}=\var^{1/2}(V_{n,1})$.
Under $\bH_0$, without loss of generality, we assume $\bm\mu_i= \bm 0$.
Note that
\begin{equation*}
	\cov\left(\frac{T_{n,1}}{\sigma_{01}},\frac{T_{n,2}}{\sigma_{02}}\right)=\sum_{k=2}^{n}\sum_{i=1}^{k-1}a_{i,k}^2\frac{\cov\left\{\X_{i}^{\top}\X_{k},(\X_{i}^{\top}\X_{k})^2\right\}}{\sigma_{01}\sigma_{02}}.
\end{equation*}
By Lemma~\ref{lemma_basic_results} (3),
\begin{equation*}
	\cov\left\{\X_{i}^{\top}\X_{k},(\X_{i}^{\top}\X_{k})^2\right\}=o\{\tr^{3/2}(\bms_1^2)\}.
\end{equation*}
Since $\sigma_{01}\sim n\tr^{1/2}(\bms_1^2)$, $\sigma_{02} \sim n\tr(\bms_1^2)$, and $\sum_{k=2}^{n}\sum_{i=1}^{k-1}a_{i,k}^2=O(n^2)$, we have 
\begin{equation*}
	\cov\left(\frac{T_{n,1}}{\sigma_{01}},\frac{T_{n,2}}{\sigma_{02}}\right)=o(1).
\end{equation*}
Below, it suffices to show that for any $a,b\in \mathbb{R}$, $aM_{n}/\sigma_{01}+bV_{n,1}/\sigma_{02}$ converges to a normal distribution.
Define
\begin{equation*}
	G_{n,k}=\X_{k}^{\top}\sum_{i=1}^{k-1}a_{i,k}\X_l
\end{equation*}
and 
\begin{equation*}
	D_{n,k}=\tr\left\{(\X_{k}\X_{k}^{\top}-\bms_1)\sum_{i=1}^{k-1}a_{i,k}(\X_i\X_{i}^{\top}-\bms_1) \right\}.
\end{equation*}
Clearly,
\begin{equation*}
	aT_{n,1}/\sigma_{01}+bT_{n,2}/\sigma_{02}=\sum_{k=1}^{n}aG_{n,k}/\sigma_{01}+bD_{n,k}/\sigma_{02}. 
\end{equation*}
In addition, it is easy to verify that for any $n$, $\{aG_{n,k}/\sigma_{01}+bD_{n,k}/\sigma_{02}\}_{k=1}^{n}$ is a martingale difference sequence with respect to the $\sigma$-fields $\{\mF_{k}\}_{k=1}^{n}$, where $\mathcal{F}_0=\{\varnothing ,\Omega\}$, $\mathcal{F}_k=\sigma\{\X_1,\ldots \X_k\}$, for $k=1,\ldots,n$.
By the proof of Lemma~\ref{lemma_normal} and the proof of Theorem~4 in \cite{yu2023power}, we only need to show that as $n,p\rightarrow \infty$,
\begin{equation*}
	\frac{1}{\sigma_{01}\sigma_{02}}\sum_{k=1}^{n}\E(G_{n,k}D_{n,k}\mid \mF_{k-1})\rightarrow 0
\end{equation*}
in probability.
For $1\leq i<k\leq n$, let $\psi_{i,k}=\X_{k}^{\top}\X_{i}$ and $\xi_{i,k}=\tr\left\{
(\X_k\X_{k}^{\top}-\bms_1)(\X_i\X_{i}^{\top}-\bms_1) \right\}$. 
Denote $\zeta_i=\E(\psi_{i,k}\xi_{i,k}\mid \X_i)$, $\phi_{i_1,i_2}=\E(\psi_{i_1,k}\xi_{i_2,k}\mid \X_{i_1},\X_{i_2})$.
Then,
\begin{equation*}
	\sum_{k=1}^{n}\E( G_{n,k}D_{n,k} \mid \mathcal{F}_{k-1})=\sum_{k=1}^{n}\sum_{i=1}^{k-1}a_{i,k}^{2}\zeta_i+\sum_{k=1}^{n}\sum_{i_1\neq i_2}^{k-1}a_{i_1,k}a_{i_2,k}\phi_{i_1,i_2},
\end{equation*}
and
\begin{align*}
	&\var\left\{  \sum_{k=1}^{n}\E( G_{n,k}D_{n,k} \mid \mathcal{F}_{k-1})  \right\}\\
	=&\sum_{i=1}^{n-1}\left(\sum_{k=i+1}^{n}a_{i,k}^2\right)^2\var(\zeta_i)+4\sum_{i_1<i_2}\left(\sum_{k=i_2+1}^{n}a_{i_1,k}a_{i_2,k}\right)^2\var(\phi_{i_1,i_2})\\
	=&n^3O\{\tr^{3}(\bms_1^2)\}+n^4o\{\tr^{3}(\bms_1^2)\},
\end{align*}
where we use the facts that $\sum_{i=1}^{n-1}\left(\sum_{k=i+1}^{n}a_{i,k}^2\right)^2=O(n^3)$,
$\sum_{i_1<i_2}\left(\sum_{k=i_2+1}^{n}a_{i_1,k}a_{i_2,k}\right)^2=O(n^4)$, and the results of Lemma~\ref{lemma_basic_results} (3) that $\E\zeta_i^2=O\{\tr^3(\bms_1^2)\}$, $\E(\phi_{i_1,i_2}^2)=o\{\tr^3(\bms_1^2)\}$ for $i_1\neq i_2$ and $\E(\phi_{i_1,i_2}\phi_{i_3,i_4})=0$ if $i_1 \notin \{i_3,i_4\}$ or $i_2 \notin \{i_3,i_4\}$.
Note that $\sigma_{01}^2\sigma_{02}^2 \sim n^4\tr^3(\bms_1^2)$, thus
\begin{equation*}
	\frac{\var\{\sum_{k=1}^{n}\E( G_{n,k}D_{n,k} \mid \mathcal{F}_{k-1})\}}{\sigma_{01}^2\sigma_{02}^2}\rightarrow 0,
\end{equation*}
which implies
\begin{equation*}
	\frac{1}{\sigma_{01}\sigma_{02}}\sum_{k=1}^{n}\E(G_{n,k}D_{n,k}\mid \mF_{k-1})\rightarrow 0
\end{equation*}
in probability.
The theorem follows. \done

\subsection{Proof of Proposition~1}

Recall that
\begin{equation*}
	\widehat{\tr(\bms_1^2)}=\frac{1}{4(n-3)}\sum_{i=1}^{n-3} \left\{(\X_{i}-\X_{i+1})^{\top}(\X_{i+2}-\X_{i+3})\right\}^2.
\end{equation*}
Under $\bH_0$, we have $\E\{\widehat{\tr(\bms_1^2)}\}=\tr(\bms_1^2)$ and $\var\{\widehat{\tr(\bms_1^2)}\}=O\{n^{-1}\tr^2(\bms_1^2)\}$.
Since 
\begin{equation*}
	\var\left\{\widehat{\tr(\bms_1^2)}\right\}\big/\tr^2(\bms_1^2)\rightarrow 0, 
\end{equation*}
it holds that $\widehat{\tr(\bms_1^2)}/\tr(\bms_1^2)\rightarrow 1$ in probability as $n,p\rightarrow \infty$.
By Slutsky’s Theorem, the results in Theorem~1 remain valid when the plug-in estimators $\widehat{\sigma}_1^2$ and $\widehat{\sigma}_2^2$ are used in place of $\var_0(M_n)$ and $\var_0(V_n)$, respectively. \done

\subsection{Proof of Theorem~2}

Denote $\bm\delta=\bm\mu_n-\bm\mu_1$.
Note that the values of $M_n(\tau)$ and $V_n(\tau)$ do not depend on $\bm\mu_1$.
Without loss of generality, we assume $\bm\mu_1=\bm 0$ and $\bm\mu_n=\bm\delta$.

We first prove $(i)$.
For $\tau\leq \tau^{\ast}$, we can write $M_n(\tau)$ as 
\begin{align*}
	M_n(\tau)=&\frac{1}{P_{\tau}^{2}}\sum_{i\neq j}^{\tau}\Y_{i}^{\top}\Y_{j}+\frac{1}{P_{n-\tau}^{2}}\sum_{i,j=\tau+1,i\neq j}^{n}\Y_{i}^{\top}\Y_{j}-\frac{2}{\tau(n-\tau)}\sum_{i=1}^{\tau}\sum_{j=\tau+1}^{n}\Y_{i}^{\top}\Y_{j}\\
	&+\frac{2(n-\tau^{\ast}-1)\sum_{i=\tau^{\ast}+1}^{n}\Y_{i}^{\top}\bm\delta +2(n-\tau^{\ast})\sum_{i=\tau+1}^{\tau^{\ast}}\Y_{i}^{\top}\bm\delta}{(n-\tau)(n-\tau-1)}-\frac{2(n-\tau^{\ast})\sum_{i=1}^{\tau}\Y_{i}^{\top}\bm\delta}{\tau(n-\tau)}\\
	&+\frac{(n-\tau^{\ast})(n-\tau^{\ast}-1)\|\bm\delta\|^2}{(n-\tau)(n-\tau-1)}.
\end{align*}
The decomposition of $M_n(\tau)$ for $\tau>\tau^{\ast}$ is similar. 
Based on the decomposition above, by straightforward calculation, we have $\E(M_n)\sim n^2\|\bm\delta\|^2$ and 
\begin{equation*}
	\var(M_n)=O[n^2\tr\{(\bms_1+\bms_n)^2\}]+O\{n^3\bm\delta^{\top}(\bms_1+\bms_n)\bm\delta\},
\end{equation*}
by noting $\var(\Y_{i}^{\top}\Y_j)\leq \tr\{(\bms_1+\bms_n)^2\}$ for $i\neq j$.
Since $\tr(\bms_1^4)=o\{\tr^2(\bms_1^2)\}$ and $\tr(\bms_n^4)=o\{\tr^2(\bms_n^2)\}$, we have $\lambda_{\max}(\bms_1)=o\{\tr^{1/2}(\bms_1^2)\}$ and $\lambda_{\max}(\bms_n)=o\{\tr^{1/2}(\bms_n^2)\}$, where $\lambda_{\max}(\bms_1)$ denotes the largest eigenvalue of the matrix $\bms_1$.
It is easy to show $\lambda_{\max}(\bms_1+\bms_n)=o[\tr^{1/2}\{(\bms_1+\bms_n)^2\}]$.
Thus, $\bm\delta^{\top}(\bms_1+\bms_n)\bm\delta=o[\|\bm\delta\|^2\tr^{1/2}\{(\bms_1+\bms_n)^2\}]$.
Under the condition
\begin{equation*}
	n\|\bm\delta\|^2/\sqrt{\tr\{(\bms_1+\bms_n)^2\}}\rightarrow \infty,
\end{equation*}
we have $\E(M_n)/\var_{0}^{1/2}(M_n)\rightarrow \infty$ and $M_n=\E(M_n)\{1+o_p(1)\}$.
Recall that
\begin{align*}
	\widehat{\tr(\bms_1^2)}=&\frac{1}{4(n-3)}\sum_{i=1}^{n-3} \left\{(\X_{i}-\X_{i+1})^{\top}(\X_{i+2}-\X_{i+3})\right\}^2\\
	=&\frac{\sum_{i=1}^{\tau^{\ast}-3} \left\{(\X_{i}-\X_{i+1})^{\top}(\X_{i+2}-\X_{i+3})\right\}^2+\sum_{i=\tau^{\ast}+1}^{n-3} \left\{(\X_{i}-\X_{i+1})^{\top}(\X_{i+2}-\X_{i+3})\right\}^2}{4(n-3)}\\
	&+\frac{1}{4(n-3)}\sum_{i=\tau^{\ast}-2}^{\tau^{\ast}} \left\{(\X_{i}-\X_{i+1})^{\top}(\X_{i+2}-\X_{i+3})\right\}^2\\
	:=& \Delta_1+\Delta_2.
\end{align*}
Since $\Delta_1>0$ and $\E(\Delta_1) \sim \tr\{(\bms_1+\bms_n)^2\}$, we have $\Delta_1=\E(\Delta_1)\{1+o_p(1)\}$.
Next, we consider $\Delta_2$.
For $i=\tau^{\ast}$,
\begin{align*}
	&\E\left\{(\X_{i}-\X_{i+1})^{\top}(\X_{i+2}-\X_{i+3})\right\}^2\\
	=&\E\left\{(\Y_{i}-\Y_{i+1}-\bm\delta)^{\top}(\Y_{i+2}-\Y_{i+3})\right\}^2\\
	=&\E\left\{(\Y_{i}-\Y_{i+1})^{\top}(\Y_{i+2}-\Y_{i+3})-\bm\delta^{\top}(\Y_{i+2}-\Y_{i+3})\right\}^2\\
	=&O_p[\tr\{(\bms_1+\bms_n)^2\}]+O_p\{n^{-1}\bm\delta^{\top}(\bms_1+\bms_n)\bm\delta\}.
\end{align*}
Similarly, for $i=\tau^{\ast}-2,\tau^{\ast}-1$, we also have
\begin{equation*}
	\E[\{(\X_{i}-\X_{i+1})^{\top}(\X_{i+2}-\X_{i+3})\}^2]=O[\tr\{(\bms_1+\bms_n)^2\}]+O\{\bm\delta^{\top}(\bms_1+\bms_n)\bm\delta\}.
\end{equation*}
Thus $\Delta_2=O_p[n^{-1}\tr\{(\bms_1+\bms_n)^2\}]+O_p\{n^{-1}\bm\delta^{\top}(\bms_1+\bms_n)\bm\delta\}$ and
\begin{equation*}
	\widehat{\tr(\bms_1^2)}=O_p[\tr\{(\bms_1+\bms_n)^2\}]+O_p\{n^{-1}\bm\delta^{\top}(\bms_1+\bms_n)\bm\delta\}.
\end{equation*}
Then,
\begin{equation*}
	\widehat{\sigma}_1=O(n)\left\{\widehat{\tr(\bms_1^2)}\right\}^{1/2}=o_p\{\E(M_n)\}.
\end{equation*}
Together with $M_n=\E(M_n)\{1+o_p(1)\}$, we have $\widehat{\sigma}_{1}^{-1}M_n\rightarrow \infty$, $p_{M_n}\rightarrow 0$, and $T_n\rightarrow \infty$ in probability as $n,p\rightarrow \infty$.

~\\
Next, we prove $(ii)$. 
Write
\begin{align*}
	&\left\{(\X_{i_1}-\X_{i_3})^{\top}(\X_{i_2}-\X_{i_4})\right\}^2\\
	=&\left\{(\Y_{i_1}-\Y_{i_3}+\E\X_{i_1}-\E\X_{i_3})^{\top}(\Y_{i_2}-\Y_{i_4}+\E\X_{i_2}-\E\X_{i_4})\right\}^2\\
	= &\left\{(\Y_{i_1}-\Y_{i_3})^{\top}(\Y_{i_2}-\Y_{i_4})\right\}^2+\Delta_{i_1,i_2,i_3,i_4,1}+\Delta_{i_1,i_2,i_3,i_4,2}
\end{align*}
where
\begin{align*}
	&\Delta_{1}(i_1,i_2,i_3,i_4)\\
	=&2(\Y_{i_1}-\Y_{i_3})^{\top}(\Y_{i_2}-\Y_{i_4})\E^{\top}(\X_{i_1}-\X_{i_3})(\Y_{i_2}-\Y_{i_4})\\
	&+2(\Y_{i_1}-\Y_{i_3})^{\top}(\Y_{i_2}-\Y_{i_4})\E^{\top}(\X_{i_2}-\X_{i_4})(\Y_{i_1}-\Y_{i_3})\\
	&+2(\Y_{i_1}-\Y_{i_3})^{\top}(\Y_{i_2}-\Y_{i_4})\E^{\top}(\X_{i_1}-\Y_{i_3})\E(\X_{i_2}-\X_{i_4}),
\end{align*}
and 
\begin{align*}
	&\Delta_{2}(i_1,i_2,i_3,i_4)\\
	=&\left\{\E^{\top}(\X_{i_1}-\X_{i_3})(\Y_{i_2}-\Y_{i_4})+\E^{\top}(\X_{i_2}-\X_{i_4})(\Y_{i_1}-\Y_{i_3})+\E^{\top}(\X_{i_1}-\X_{i_3})\E(\X_{i_2}-\X_{i_4})\right\}^2\\
	\leq& 3\{\E^{\top}(\X_{i_1}-\X_{i_3})(\Y_{i_2}-\Y_{i_4})\}^2+3\{\E^{\top}(\X_{i_2}-\X_{i_4})(\Y_{i_1}-\Y_{i_3})\}^2\\
	&+3\{\E^{\top}(\X_{i_1}-\X_{i_3})\E(\X_{i_2}-\X_{i_4})\}^2.
\end{align*}
Then, $V_n(\tau)$ can be decomposed as 
\begin{equation}\label{equation_V_tau_decomposition}
	V_n(\tau):=V_{n,1}(\tau)+V_{n,2}(\tau)+V_{n,3}(\tau),
\end{equation}
where
\begin{align*}
	V_{n,1}(\tau)=&\frac{\sum_{1\leq i_1,i_2,i_3,i_4\leq \tau}^{\ast}\left\{(\Y_{i_1}-\Y_{i_3})^{\top}(\Y_{i_2}-\Y_{i_4})\right\}^2}{4P_{\tau}^{4}} \\
	&+ \frac{\sum_{\tau+1\leq i_1,i_2,i_3,i_4\leq n}^{\ast}\left\{(\Y_{i_1}-\Y_{i_3})^{\top}(\Y_{i_2}-\Y_{i_4})\right\}^2}{4P_{n-\tau}^{4}}\\
	&-\frac{\sum_{i_1\neq i_3}^{\tau}\sum_{\tau+1\leq i_2,i_4\leq n}^{\ast} \left\{(\Y_{i_1}-\Y_{i_3})^{\top}(\Y_{i_2}-\Y_{i_4})\right\}^2 }{2\tau(\tau-1)(n-\tau)(n-\tau-1)},\\
	V_{n,2}(\tau)=&\frac{\sum_{1\leq i_1,i_2,i_3,i_4\leq \tau}^{\ast}\Delta_{1}(i_1,i_2,i_3,i_4)}{4P_{\tau}^{4}}+ \frac{\sum_{\tau+1\leq i_1,i_2,i_3,i_4\leq n}^{\ast}\Delta_{1}(i_1,i_2,i_3,i_4)}{4P_{n-\tau}^{4}}\\
	&-\frac{\sum_{i_1\neq i_3}^{\tau}\sum_{\tau+1\leq i_2,i_4\leq n}^{\ast} \Delta_{1}(i_1,i_2,i_3,i_4) }{2\tau(\tau-1)(n-\tau)(n-\tau-1)},\\
	V_{n,3}(\tau)=&\frac{\sum_{1\leq i_1,i_2,i_3,i_4\leq \tau}^{\ast}\Delta_{2}(i_1,i_2,i_3,i_4)}{4P_{\tau}^{4}}+ \frac{\sum_{\tau+1\leq i_1,i_2,i_3,i_4\leq n}^{\ast}\Delta_{2}(i_1,i_2,i_3,i_4)}{4P_{n-\tau}^{4}}\\
	&-\frac{\sum_{i_1\neq i_3}^{\tau}\sum_{\tau+1\leq i_2,i_4\leq n}^{\ast} \Delta_{2}(i_1,i_2,i_3,i_4) }{2\tau(\tau-1)(n-\tau)(n-\tau-1)}.
\end{align*}
Define 
\begin{equation*}
	V_{n,i}=\sum_{\tau=4}^{n-4}\left\{\frac{\tau(n-\tau)}{n}V_{n,i}(\tau)\right\}, \quad i=1,2,3. 
\end{equation*}
Then $V_n=V_{n,1}+V_{n,2}+V_{n,3}$.
Below we bound each term respectively.

For $V_{n,2}$, by the facts
\begin{align*}
	&\var\left\{(\Y_{i_1}-\Y_{i_3})^{\top}(\Y_{i_2}-\Y_{i_4})\E^{\top}(\X_{i_1}-\X_{i_3})(\Y_{i_2}-\Y_{i_4})\right\}\\
	= &O\left[\tr\{(\bms_1+\bms_n)^2\}\bm\delta^{\top}(\bms_1+\bms_n)\bm\delta\right]+O\{\bm\delta^{\top}(\bms_1+\bms_n)^3\bm\delta\}+o\left[\tr^{3/2}\{(\bms_{1}+\bms_n)^2\}\|\bm\delta\|^2\right]\\
	=&o\left[\tr^{3/2}\{(\bms_{1}+\bms_n)^2\}\|\bm\delta\|^2\right],\\
	&\var\left\{(\Y_{i_1}-\Y_{i_3})^{\top}(\Y_{i_2}-\Y_{i_4})\E^{\top}(\X_{i_1}-\Y_{i_3})\E(\X_{i_2}-\Y_{i_4})\right\}\leq 4\tr\{(\bms_1+\bms_n)^2\}\|\bm\delta\|^4,
\end{align*}
and the condition
\begin{equation*}
	n\|\bm\delta\|^2/\sqrt{\tr\{(\bms_1+\bms_n)^2\}}=O(1),
\end{equation*}
we conclude that
\begin{align*}
	\var(V_{n,2})=&o\left[n^3\tr^{3/2}\{(\bms_{1}+\bms_n)^2\}\|\bm\delta\|^2\right]+O\left[n^2\tr\{(\bms_1+\bms_n)^2\}\|\bm\delta\|^4\right]\\
	=&o\left[n^2\tr^{2}\{(\bms_{1}+\bms_n)^2\}\right]
\end{align*}
and 
\begin{equation*}
	V_{n,2}=o_p\left[n\tr\{(\bms_{1}+\bms_n)^2\}\right].
\end{equation*}

For $V_{n,3}$, based on the facts
\begin{align*}
	&\E[\{\E^{\top}(\X_{i_1}-\X_{i_3})(\Y_{i_2}-\Y_{i_4})\}^2] \leq 2\bm\delta^{\top}(\bms_1+\bms_n)\bm\delta,\\
	& \{\E^{\top}(\X_{i_1}-\X_{i_3})\E(\X_{i_2}-\X_{i_4})\}^2\leq \|\bm\delta\|^4,
\end{align*}
and the condition 
\begin{equation*}
	n\|\bm\delta\|^2/\sqrt{\tr\{(\bms_1+\bms_n)^2\}}=O(1),
\end{equation*}
we have
\begin{equation*}
	\E|V_{n,3}|=O(n^2\|\bm\delta\|^4)+O\{n^2\bm\delta^{\top}(\bms_1+\bms_n)\bm\delta\}=o\left[n\tr\{(\bms_{1}+\bms_n)^2\}\right].
\end{equation*}
Thus
\begin{equation*}
	V_{n,3}=o_p\left[n\tr\{(\bms_{1}+\bms_n)^2\}\right].
\end{equation*}

For $V_{n,1}$, since 
\begin{equation*}
	\E\{V_{n,1}(\tau)\}=\frac{(n-\tau^{\ast})(n-\tau^{\ast}-1)}{(n-\tau)(n-\tau-1)}\tr\{(\bms_1-\bms_n)^2\}
\end{equation*}
for $\tau\leq \tau^{\ast}$, and 
\begin{equation*}
	\E\{V_{n,1}(\tau)\}=\frac{\tau^{\ast}(\tau^{\ast}-1)}{\tau(\tau-1)}\tr\{(\bms_1-\bms_n)^2\}
\end{equation*}
for $\tau>\tau^{\ast}$, we have $\E(V_{n,1})\sim n^2\tr\{(\bms_1-\bms_n)^2\}$.
Decompose $V_{n,1}$ as
\begin{equation*}
	V_{n,1}=\sum_{\tau=4}^{n-4}\left\{\frac{\tau(n-\tau)}{n}V_{n,1,1}(\tau)\right\}+\sum_{\tau=4}^{n-4}\left\{\frac{\tau(n-\tau)}{n}V_{n,1,2}(\tau)\right\} := V_{n,1,1}+V_{n,1,2},
\end{equation*}
where 
\begin{equation*}
	V_{n,1,1}(\tau)=\frac{1}{P_{\tau}^{2}}\sum_{i\neq j}^{\tau}(\Y_{i}^{\top}\Y_{j})^2+\frac{1}{P_{n-\tau}^{2}}\sum_{i,j=\tau+1,i\neq j}^{n}(\Y_{i}^{\top}\Y_{j})^2-\frac{2}{\tau(n-\tau)}\sum_{i=1}^{\tau}\sum_{j=\tau+1}^{n}(\Y_{i}^{\top}\Y_{j})^2,
\end{equation*}
and $V_{n,1,2}(\tau)=V_{n,1}(\tau)-V_{n,1,1}(\tau)$.
Similar to the proof of Lemma~\ref{lemma_normal}, it is easy to check that $\var(V_{n,1,2})=o\left[{n^2\tr^2\{(\bms_1+\bms_n)^2\}}\right]$.
By Lemma~\ref{lemma_two_sample_cov}, it holds uniformly that
\begin{align*}
	\var\{V_{n,1,1}(\tau)\}\leq & C_1\left[\left(\frac{1}{\tau}+\frac{1}{n-\tau}\right)\tr\{(\bms_1+\bms_n)^2\}\right]^2\{1+o(1)\}\\
	&+C_2\left(\frac{1}{\tau}+\frac{1}{n-\tau}\right)\tr\{(\bms_1+\bms_n)^2\}\tr\{(\bms_1-\bms_n)^2\}\{1+o(1)\},
\end{align*}
for some constants $C_1, C_2>0$.
Then, by the Cauchy-Schwarz inequality,
\begin{align*}
	\var(V_{n,1,1})\leq& n\sum_{\tau=4}^{n-4}\left[\frac{\tau^2(n-\tau)^2}{n^2}\var\{V_{n,1,1}(\tau)\}\right]\\
	=&O\left[n^2\tr^2\{(\bms_1+\bms_n)^2\}\right] +O\left[n^3\tr\{(\bms_1+\bms_n)^2\}\tr\{(\bms_1-\bms_n)^2\}\right].
\end{align*}
By $n\tr\{(\bms_1-\bms_n)^2\}/\tr\{(\bms_1+\bms_n)^2\}\rightarrow \infty$, it holds that $\var^{1/2}(V_{n,1})=o_p\{E(V_{n,1})\}$.
Thus, $V_{n,1}=\E(V_{n,1})\{1+o_p(1)\}$.
Together with $V_{n,2}=o_p\left[n\tr\{(\bms_{1}+\bms_n)^2\}\right]$ and $V_{n,3}=o_p\left[n\tr\{(\bms_{1}+\bms_n)^2\}\right]$, we conclude that
\begin{equation*}
	V_n= \E(V_{n,1})\{1+o_p(1)\}.
\end{equation*}
By
\begin{equation*}
	\widehat{\tr(\bms_1^2)}=O_p(\tr\{(\bms_1+\bms_n)^2\})+O_p\{n^{-1}\bm\delta^{\top}(\bms_1+\bms_n)\bm\delta\},
\end{equation*}
we have
\begin{equation*}
	\widehat{\sigma}_2=O_p(n\tr\{(\bms_1+\bms_n)^2\})+O_p\{\bm\delta^{\top}(\bms_1+\bms_n)\bm\delta\}=o_p\{\E(V_{n,1})\}.
\end{equation*}
Thus we have $\widehat{\sigma}_{2}^{-1}V_n\rightarrow \infty$, $p_{V_n}\rightarrow 0$, and $T_n\rightarrow \infty$ in probability as $n,p\rightarrow \infty$. \done

\subsection{Proof of Theorem~3}

Without loss of generality, we assume $\widehat{\tau}\leq \tau^{\ast}$.
To establish the consistency of the estimator $\widehat{\tau}$, we need to show that, for any $\epsilon>0$, as $n,p\rightarrow \infty$,
\begin{equation*}
	\Pr\left\{\max_{\lambda_n\leq \tau\leq \tau^{\ast}-n\epsilon}T_n(\tau)\geq T_n(\tau^{\ast})\right\}\rightarrow 0.
\end{equation*}
Recall
\begin{equation*}
	T_n(\tau)=-2\log\{1-\Phi(\tilde{M}_n(\tau)\}-2\log\{1-\Phi(\tilde{V}_n(\tau)\}.
\end{equation*}
Note that the function $f(x)=-2\log\{1-\Phi(x)\}$ is an increasing function, here we assume that for all $\tau\in[\lambda_n, \tau^{\ast}-n\epsilon]$, $\tilde{M}_{n}(\tau)\geq c_0$ and $\tilde{V}_{n}(\tau)\geq c_0$, for some constant $c_0>0$. 
If $\tilde{M}_{n}(\tau)< c_0$ or $\tilde{V}_{n}(\tau)< c_0$, then the corresponding terms $-2\log\{1-\Phi(\tilde{M}_n(\tau)\}$ or $-2\log\{1-\Phi(\tilde{V}_n(\tau)\}$ will have an upper bound $-2\log\{1-\Phi(c_0)\}$, which is negligible in the comparison of the size of $T_{n}(\tau)$ and $T_n(\tau^{\ast})$ due to the fact that $T_n(\tau^{\ast})\rightarrow \infty$ in probability under the designed conditions (we will show it later). 
Then, by Lemma~\ref{lemma_normal_inequality}, we have that for all $\tau\in[\lambda_n, \tau^{\ast}-n\epsilon]$,
\begin{equation}\label{equation_T_n}
	T_n(\tau)\leq 2\log(2\pi)+4\log(1+c_{0}^{-2})+2\log\{\tilde{M}_n(\tau)\}+\tilde{M}_{n}^{2}(\tau)+2\log\{\tilde{V}_n(\tau)\}+\tilde{V}_{n}^{2}(\tau).
\end{equation}
We first investigate the size of $\tilde{M}_n(\tau)$ and $\tilde{V}_n(\tau)$, for $\lambda_n\leq \tau\leq \tau^{\ast}$.
Let $\bm\delta=\bm\mu_n-\bm\mu_1$. 

By Lemma~\ref{lemma_M_tau_V_tau}, we have that
\begin{equation}\label{equation_M_tau}
	M_n(\tau)= \frac{(n-\tau^{\ast})(n-\tau^{\ast}-1)}{(n-\tau)(n-\tau-1)}\|\bm\delta\|^2 +O_p[n^{-1}\tr^{1/2}\{(\bms_1+\bms_n)^2\}]+O_p[\{n^{-1}\bm\delta(\bms_1+\bms_n) \bm\delta\}^{1/2}],
\end{equation}
and
\begin{align}\label{equation_V_tau}
	V_n(\tau)=&\frac{(n-\tau^{\ast})(n-\tau^{\ast}-1)}{(n-\tau)(n-\tau-1)}\tr\{(\bms_1-\bms_n)^2\}\\
	&+O_p[n^{-1}\tr\{(\bms_1+\bms_n)^2\}]+O_p\left([ n^{-1}\tr\{(\bms_1+\bms_n)^2\}\tr\{ (\bms_1-\bms_n)^2 \}]^{1/2}\right) \nonumber \\ 
	&+o_p\left[n^{-1/2}\tr^{3/4}\{(\bms_{1}+\bms_n)^2\}\|\bm\delta\|\right]+o_p\left[\|\bm\delta\|^2\sqrt{\tr\{(\bms_1+\bms_n)^2\}} \right]+O_p(\|\bm\delta\|^4)\nonumber
\end{align}
hold uniformly for all $\tau\in [\lambda_n,\tau^{\ast}]$. 
Therefore,
\begin{small}
	\begin{align}\label{equation_M2_difference}
		&\left\{\frac{\tau(n-\tau)}{n}M_n(\tau)\right\}^2-\left\{\frac{\tau^{\ast}(n-\tau^{\ast})}{n}M(\tau^{\ast})\right\}^2\\  \nonumber
		=&\left(-\frac{(n-1)(n-\tau^{\ast})(\tau^{\ast}-\tau)}{n(n-\tau-1)}\|\bm\delta\|^2+O_p[\tr^{1/2}\{(\bms_1+\bms_n)^2\}]+O_p[\{n\bm\delta(\bms_1+\bms_n) \bm\delta\}^{1/2}]\right)\\ \nonumber
		&\times\left( \left\{\frac{\tau(n-\tau^{\ast})(n-\tau^{\ast}-1)}{n(n-\tau-1)}+\frac{\tau^{\ast}(n-\tau^{\ast})}{n}\right\}\|\bm\delta\|^2+O_p[\tr^{1/2}\{(\bms_1+\bms_n)^2\}]+O_p[\{n\bm\delta(\bms_1+\bms_n) \bm\delta\}^{1/2}]  \right),
	\end{align}
\end{small}
and
\begin{align}\label{equation_V2_difference}
	&\left\{\frac{\tau(n-\tau)}{n}V_n(\tau)\right\}^2-\left\{\frac{\tau^{\ast}(n-\tau^{\ast})}{n}V(\tau^{\ast})\right\}^2 \\
	=&\left[-\frac{(n-1)(n-\tau^{\ast})(\tau^{\ast}-\tau)}{n(n-\tau-1)}\tr\{ (\bms_1-\bms_n)^2\} +R_n \right] \nonumber \\ 
	&\times \left[\frac{(n-\tau^{\ast})(n-\tau^{\ast}-1)}{n}\left( 
	\frac{\tau}{n-\tau-1}+\frac{\tau^{\ast}}{n-\tau^{\ast}-1}\right)\tr\{ (\bms_1-\bms_n)^2\} +R_n \right],\nonumber
\end{align}
where 
\begin{align*}
	R_n=&O_p[\tr\{(\bms_1+\bms_n)^2\}]+O_p\left([n\tr\{(\bms_1+\bms_n)^2\}\tr\{ (\bms_1-\bms_n)^2 \}]^{1/2}\right)\\
	&+o_p\left[n^{1/2}\tr^{3/4}\{(\bms_{1}+\bms_n)^2\}\|\bm\delta\|\right]+o_p\left[n\|\bm\delta\|^2\sqrt{\tr\{(\bms_1+\bms_n)^2\}} \right]+O_p(n\|\bm\delta\|^4).
\end{align*}

Under condition $(i)$, $n\|\bm\delta\|^2/\sqrt{\tr\{(\bms_1+\bms_n)^2\}}\rightarrow \infty$.
By Equation (\ref{equation_M_tau}), we have
\begin{equation*}
	\frac{\tau^{\ast}(n-\tau^{\ast})}{n}M_n(\tau^{\ast})=c n\|\bm\delta\|^2\{1+o_p(1)\},
\end{equation*}
for some constant $c>0$.
Since
\begin{equation*}
	\widehat{\tr(\bms_1^2)}=O_p(\tr\{(\bms_1+\bms_n)^2\})+O_p\{n^{-1}\bm\delta^{\top}(\bms_1+\bms_n)\bm\delta\}=o_p(n^2\|\bm\delta\|^4),
\end{equation*}
it holds that $\tilde{M}_n(\tau^{\ast})\rightarrow \infty$ in probability. 
Thus, by Lemma~\ref{lemma_normal_inequality},
\begin{equation*}
	T_n(\tau^{\ast})> \log(2\pi)+2\log\{\tilde{M}_n(\tau^{\ast})\}+\tilde{M}_{n}^{2}(\tau^{\ast})
\end{equation*}
holds in probability.
Then, by Inequality (\ref{equation_T_n}),
\begin{align*}
	T_n(\tau)- T_n(\tau^{\ast}) <&  \tilde{M}_{n}^{2}(\tau) 
	- \tilde{M}_{n}^{2}(\tau^{\ast})+ 2\log\{\tilde{M}_n(\tau)\}-2\log\{\tilde{M}_n(\tau^{\ast})\}\\
	& + 2\log\{\tilde{V}_n(\tau)\}+\tilde{V}_{n}^{2}(\tau)+C,
\end{align*}
where $C>0$ is a constant.
By the condition $n\|\bm\delta\|^2/\sqrt{\tr\{(\bms_1+\bms_n)^2\}}\rightarrow \infty$ and Equation (\ref{equation_M2_difference}),
we have 
\begin{equation*}
	\left\{\frac{\tau(n-\tau)}{n}M_n(\tau)\right\}^2-\left\{\frac{\tau^{\ast}(n-\tau^{\ast})}{n}M(\tau^{\ast})\right\}^2 \leq -c_1n^2\epsilon\|\bm\delta\|^4+o_p(n^2\|\bm\delta\|^4)\rightarrow -\infty
\end{equation*}
in probability, where $c_1>0$ is a constant.
Thus, we have $\tilde{M}_n(\tau)-\tilde{M}_n(\tau^{\ast})<0$ and $\log\{\tilde{M}_n(\tau)\}-\log\{\tilde{M}_n(\tau^{\ast})\}<0$ hold in probability.
Under the condition $n\tr\{(\bms_1-\bms_n)^2\}/\tr\{(\bms_1+\bms_n)^2\}=O(1)$, together with the fact $\{\widehat{\tr(\bms_1^2)}\}^{-1}\tr\{(\bms_1+\bms_n)^2\}=O_p(1)$, we have
\begin{align*}
	\tilde{V}_{n}^{2}(\tau)=
	\frac{o_p(n\|\bm\delta\|^2)}{\{\widehat{\tr(\bms_1^2)}\}^{1/2}}+
	\frac{o_p(n^2\|\bm\delta\|^4)}{\widehat{\tr(\bms_1^2)}}
	+\frac{O_p(n^2\|\bm\delta\|^8)}{\{\widehat{\tr(\bms_1^2)}\}^2}+O_p(1).
\end{align*}
Further, since $\widehat{\tr(\bms_1^2)}/(n^2\|\bm\delta\|^4)=o_p(1)$ and $\|\bm\delta\|^4=o_p[\tr\{(\bms_1+\bms_n)^2\}]$,
we have 
\begin{equation*}
	\frac{n\|\bm\delta\|^2}{\{\widehat{\tr(\bms_1^2)}\}^{1/2}}=\frac{o_p(n^2\|\bm\delta\|^4)}{\widehat{\tr(\bms_1^2)}}, \quad \frac{O_p(\|\bm\delta\|^4)}{\widehat{\tr(\bms_1^2)}}=o_p(1).
\end{equation*}
Thus,
\begin{equation*}
	\tilde{V}_{n}^{2}(\tau)= \frac{o_p(n^2\|\bm\delta\|^4)}{\widehat{\tr(\bms_1^2)}}.
\end{equation*}
Then,
\begin{equation*}
	\max_{\lambda_n\leq \tau \leq \tau^{\ast}-n\epsilon } T_n(\tau)-T_n(\tau^{\ast})<  \frac{-c_1n^2\epsilon\|\bm\delta\|^4+o_p(n^2\|\bm\delta\|^4)}{\widehat{\tr(\bms_1^2)}}.  
\end{equation*}
We conclude that, for any $\epsilon>0$, as $n,p\rightarrow \infty$,
\begin{equation*}
	\Pr \left\{\max_{\lambda_n\leq \tau \leq \tau^{\ast}-n\epsilon } T_{n}(\tau)-T_{n}(\tau^{\ast})<0 \right\}\rightarrow 1.
\end{equation*}

Under condition $(ii)$, $n\tr\{(\bms_1-\bms_n)^2\}/\tr\{(\bms_1+\bms_n)^2\}\rightarrow \infty$. 
By Equation (\ref{equation_V_tau}), we have
\begin{equation*}
	\frac{\tau^{\ast}(n-\tau^{\ast})}{n}V_n(\tau^{\ast})= c n\tr\{(\bms_1-\bms_n)^2\}\{1+o_p(1)\},
\end{equation*}
for some constant $c>0$.
Since $n\|\bm\delta\|^2/\sqrt{\tr\{(\bms_1+\bms_n)^2\}}=O(1)$, we have
\begin{equation*}
	\widehat{\tr(\bms_1^2)}=O_p[\tr\{(\bms_1+\bms_n)^2\}].
\end{equation*}
It holds that $\tilde{V}_n(\tau^{\ast})\rightarrow \infty$ in probability. 
By Lemma~\ref{lemma_normal_inequality}, we have
\begin{equation*}
	T_n(\tau^{\ast})> \log(2\pi)+2\log\{\tilde{V}_n(\tau^{\ast})\}+\tilde{V}_{n}^{2}(\tau^{\ast})
\end{equation*}
holds in probability.
Then, by Inequality (\ref{equation_T_n}),
\begin{align*}
	T_n(\tau)- T_n(\tau^{\ast}) <&  \tilde{V}_{n}^{2}(\tau) 
	- \tilde{V}_{n}^{2}(\tau^{\ast})+ 2\log\{\tilde{V}_n(\tau)\}-2\log\{\tilde{V}_n(\tau^{\ast})\}\\
	& + 2\log\{\tilde{M}_n(\tau)\}+\tilde{M}_{n}^{2}(\tau)+C,
\end{align*}
where $C>0$ is a constant.
By Equation (\ref{equation_M_tau}), it holds that $2\log\{\tilde{M}_n(\tau)\}+\tilde{M}_{n}^{2}(\tau)=O_p(1)$.
By Equation (\ref{equation_V2_difference}), 
\begin{align*}
	&\left\{\frac{\tau(n-\tau)}{n}V_n(\tau)\right\}^2-\left\{\frac{\tau^{\ast}(n-\tau^{\ast})}{n}V(\tau^{\ast})\right\}^2 \\
	\leq& -c_2n^2\epsilon \tr^2\{(\bms_1-\bms_n)^2\}+o_p[n^2\ \tr^2\{(\bms_1-\bms_n)^2\}],
\end{align*}
where $c_2>0$ is a constant.
Thus,
\begin{equation*}
	\tilde{V}_{n}^{2}(\tau) 
	- \tilde{V}_{n}^{2}(\tau^{\ast})\leq \frac{-c_2n^2\epsilon \tr^2\{(\bms_1-\bms_n)^2\}+o_p[n^2\tr^2\{(\bms_1-\bms_n)^2\}]}{\{\widehat{\tr(\bms_1^2)}\}^2}\rightarrow -\infty, 
\end{equation*}
and $2\log\{\tilde{V}_n(\tau)\}<2\log\{\tilde{V}_n(\tau^{\ast})\}$ hold in probability.
Then,
\begin{equation*}
	\max_{\lambda_n\leq \tau \leq \tau^{\ast}-n\epsilon } T_n(\tau)-T_n(\tau^{\ast})<  
	\frac{-c_2n^2\epsilon \tr^2\{(\bms_1-\bms_n)^2\}+o_p[n^2\ \tr^2\{(\bms_1-\bms_n)^2\}]}{\{\widehat{\tr(\bms_1^2)}\}^2}+O_p(1).
\end{equation*}
We conclude that, for any $\epsilon>0$, as $n,p\rightarrow \infty$,
\begin{equation*}
	\Pr \left\{\max_{\lambda_n\leq \tau \leq \tau^{\ast}-n\epsilon } T_{n}(\tau)-T_{n}(\tau^{\ast})<0 \right\}\rightarrow 1.
\end{equation*}

Under condition $(iii)$, we have $\tilde{M}_n(\tau^{\ast})\rightarrow \infty$, $\tilde{V}_n(\tau^{\ast})\rightarrow \infty$ in probability, by noting that
\begin{equation*}
	\widehat{\tr(\bms_1^2)} = c\tr\{(\bms_1+\bms_n)^2\}\{1+o_p(1)\},
\end{equation*}
for some constant $c>0$.
By Lemma~\ref{lemma_normal_inequality}, we have
\begin{equation*}
	T_n(\tau^{\ast})\geq 2\log(2\pi)+2\log\{\tilde{M}_n(\tau^{\ast})\}+\tilde{M}_{n}^{2}(\tau^{\ast})+2\log\{\tilde{V}_n(\tau^{\ast})\}+\tilde{V}_{n}^{2}(\tau^{\ast})
\end{equation*}
holds in probability.

Under the condition that $n\|\bm\delta\|^2/\sqrt{\tr\{(\bms_1+\bms_n)^2\}}\rightarrow \infty$,
we have shown that
\begin{equation*}
	\left\{\frac{\tau(n-\tau)}{n}M_n(\tau)\right\}^2-\left\{\frac{\tau^{\ast}(n-\tau^{\ast})}{n}M(\tau^{\ast})\right\}^2 \leq -c_1n^2\epsilon\|\bm\delta\|^4+o_p(n^2\|\bm\delta\|^4)\rightarrow -\infty,
\end{equation*}
for some constant $c_1>0$.
Thus, $\tilde{M}_n(\tau)-\tilde{M}_n(\tau^{\ast})<0$ and $\log\{\tilde{M}_n(\tau)\}-\log\{\tilde{M}_n(\tau^{\ast})\}<0$ hold in probability.

Under the condition $n\tr\{(\bms_1-\bms_n)^2\}/\tr\{(\bms_1+\bms_n)^2\}\rightarrow \infty$, by Equation (\ref{equation_V2_difference}), we have
\begin{align*}
	&\left\{\frac{\tau(n-\tau)}{n}V_n(\tau)\right\}^2-\left\{\frac{\tau^{\ast}(n-\tau^{\ast})}{n}V(\tau^{\ast})\right\}^2 \\
	\leq &-c_2n^2\epsilon \tr^2\{(\bms_1-\bms_n)^2\}+o_p[n^2\ \tr^2\{(\bms_1-\bms_n)^2\}]\\
	&+o_p[n^2\|\bm\delta\|^4\tr\{(\bms_1+\bms_n)^2\}]+o_p[n^2\|\bm\delta\|^2\sqrt{\tr\{(\bms_1+\bms_n)^2\}}\tr\{(\bms_1-\bms_n)^2\}]\\
	&+O_p(n^2\|\bm\delta\|^8)+O_p[n^2\|\bm\delta\|^4\tr\{(\bms_1-\bms_n)^2\}],
\end{align*}
for some constant $c_2>0$.
Note that $o_p[n^2\|\bm\delta\|^2\sqrt{\tr\{(\bms_1+\bms_n)^2\}}\tr\{(\bms_1-\bms_n)^2\}]$ is also $o_p[n^2\|\bm\delta\|^4\tr\{(\bms_1+\bms_n)^2\}+n^2\tr^2\{(\bms_1-\bms_n)^2\}]$, and $n^2\|\bm\delta\|^8=o[n^2\|\bm\delta\|^4\tr\{(\bms_1+\bms_n)^2\}]$.
By $\|\bm\delta\|^4=o(\tr\{(\bms_1+\bms_n)^2\})$, we have
\begin{equation*}
	\|\bm\delta\|^4\tr\{(\bms_1-\bms_n)^2\}=o(\tr^2\{(\bms_1-\bms_n)^2\}+\|\bm\delta\|^4\tr\{(\bms_1+\bms_n)^2\}).
\end{equation*}
Based on the facts, it is easy to verify that
\begin{align}{\label{equation_main_term}}
	&\tilde{M}_{n}^{2}(\tau)+\tilde{V}_{n}^{2}(\tau) - \tilde{M}_{n}^{2}(\tau^{\ast})-\tilde{V}_{n}^{2}(\tau^{\ast}) \nonumber \\
	\leq &\frac{-c_3 n^2 \epsilon\left[\|\bm\delta\|^4\tr\{(\bms_1+\bms_n)^2\}+\tr^2\{(\bms_1-\bms_n)^2\}\right]\{1+o_p(1)\}}{\{\widehat{\tr(\bms_1^2)}\}^2},
\end{align}
for some constant $c_3>0$.
By Equation (\ref{equation_V_tau}), we have
\begin{equation*}
	\tilde{V}_n(\tau)=O_p\left[ \frac{ n\tr\{(\bms_1-\bms_n)^2\}}{\tr\{(\bms_1+\bms_n)^2\}}\right]+o_p\left[\frac{n\|\bm\delta\|^2}{\sqrt{\tr\{(\bms_1+\bms_n)^2\}}}\right].
\end{equation*}
Then, 
\begin{equation}{\label{equation_minor_term}}
	\log\{\tilde{V}_n(\tau)\}=o_p\left\{\tilde{M}_{n}^{2}(\tau)+\tilde{V}_{n}^{2}(\tau) - \tilde{M}_{n}^{2}(\tau^{\ast})-\tilde{V}_{n}^{2}(\tau^{\ast})\right\}.
\end{equation}
Since
\begin{align*}
	T_n(\tau)-T_n(\tau^{\ast}) <&  \tilde{M}_{n}^{2}(\tau) 
	- \tilde{M}_{n}^{2}(\tau^{\ast})+
	\tilde{V}_{n}^{2}(\tau) 
	- \tilde{V}_{n}^{2}(\tau^{\ast})
	+ 2\log\{\tilde{M}_n(\tau)\}-2\log\{\tilde{M}_n(\tau^{\ast})\}\\
	& +2\log\{\tilde{V}_n(\tau)\}-2\log\{\tilde{V}_n(\tau^{\ast})\}+C,
\end{align*}
we have shown that
$\tilde{V}_n(\tau^{\ast})\rightarrow \infty$ and $\log\{\tilde{M}_n(\tau)\}-2\log\{\tilde{M}_n(\tau^{\ast})\}<0$ in probability,
together with Inequality (\ref{equation_main_term}) and Equation (\ref{equation_minor_term}), it holds that
\begin{equation*}
	\max_{\lambda_n\leq \tau \leq \tau^{\ast}-n\epsilon } T_n(\tau)-T_n(\tau^{\ast})< \frac{-c_3 n^2 \epsilon\left[\|\bm\delta\|^4\tr\{(\bms_1+\bms_n)^2\}+\tr^2\{(\bms_1-\bms_n)^2\}\right]\{1+o_p(1)\}}{\{\widehat{\tr(\bms_1^2)}\}^2},
\end{equation*}
for some constant $c_3>0$.
We conclude that, for any $\epsilon>0$,
\begin{equation*}
	\Pr \left\{\max_{\lambda_n\leq \tau \leq \tau^{\ast}-n\epsilon } T_{n}(\tau)-T_{n}(\tau^{\ast})<0 \right\}\rightarrow 1.
\end{equation*}

The proof for $\widehat{\tau}>\tau^{\ast}$ is similar, we omit it here. 
Hence, under the designed conditions, we have $(\widehat{\tau}-\tau^{\ast})/n=o_p(1)$. \done

\section{Auxiliary lemmas}

The first lemma contains some basic results for the theoretical analysis.
\begin{lemma}\label{lemma_basic_results}
	Let $\Y_0$, $\Y_0^{\prime}$ and $\Y_{0}^{\prime\prime}$ be independent copies of $\Y_1=\X_1-\E(\X_1)$. Suppose $\bH_0$ holds.
	\begin{itemize}
		\item[(1)] Suppose Assumptions~1--2 hold. Then, \\
		$\E\{(\Y_{0}^{\top}\Y_{0}^{\prime})^4\}=O\{\tr^2(\bms_1^2)\}$,\\
		$\E\{(\Y_0^{\top}\bms_1\Y_0)^2\}=O\{\tr^2(\bms_1^2)\}$,\\
		$\var(\Y_0^{\top}\bms_1\Y_0)=o\{\tr^2(\bms_1^2)\}$;
		
		\item[(2)] Suppose Assumptions~1--3 hold. Then,\\
		$\E\{(\Y_{0}^{\top}\Y_{0}^{\prime})^6\}=O\{\tr^3(\bms_1^2)\}$,\\
		$\E\{(\Y_{0}^{\top}\Y_{0}^{\prime})^8\}=O\{\tr^4(\bms_1^2)\}$,\\
		$\E\{(\Y_0^{\top}\bms_1\Y_0)^3\}=O\{\tr^3(\bms_1^2)\}$,\\
		$\E\{(\Y_0^{\top}\bms_1\Y_0)^4\}=O\{\tr^4(\bms_1^2)\}$,\\
		$\var\{(\Y_{0}^{\top}\Y_{0}^{\prime})^2\}=2\tr^2(\bms_1^2)+6\tr(\bms_1^4)+o\{\tr^2(\bms_1^2)\}$,\\
		$\cov\{(\Y_{0}^{\top}\Y_{0}^{\prime})^2,(\Y_{0}^{\top}\Y_{0}^{\prime\prime})^2\}
		=2\tr(\bms_1^4)+o\{\tr^2(\bms_1^2)\}$,\\
		$\E\left[\tr^4\left\{(\Y_0\Y_0^{\top}-\bms_1)(\Y_{0}^{\prime}\Y_0^{\prime\top}-\bms_1) \right\}\right]= O\{\tr^4(\bms_1^2)\}$,\\
		$\var\left(\E\left[\tr^2\left\{(\Y_{0}^{\prime}\Y_0^{\prime\top}-\bms_1)(\Y_0\Y_0^{\top}-\bms_1) \right\}\mid \Y_0\right]\right)=O\{\tr^4(\bms_1^2)\}$,\\
		$\var\left(\E\left[\tr\left\{(\Y_0^{\prime\prime}\Y_0^{\prime\prime\top}-\bms_1)(\Y_{0}\Y_0^{\top}-\bms_1) \right\} \tr\left\{(\Y_0^{\prime\prime}\Y_0^{\prime\prime\top}-\bms_1)(\Y_{0}^{\prime}\Y_0^{\prime\top}-\bms_1) \right\}\mid \Y_0,\Y_{0}^{\prime}\right]\right)=o\{\tr^4(\bms_1^2)\}$;
		\item[(3)] Suppose Assumptions~1--4 hold. Then,\\
		$\var\left(\E\left[\Y_{0}^{\prime\top}\Y_0\tr\left\{(\Y_{0}^{\prime}\Y_0^{\prime\top}-\bms_1)(\Y_0\Y_0^{\top}-\bms_1) \right\}\mid \Y_0\right]\right)= O\{\tr^3(\bms_1^2)\}$,\\
		$\var\left(\E\left[\Y_0^{\prime\prime\top}\Y_0 \tr\left\{(\Y_0^{\prime\prime}\Y_0^{\prime\prime\top}-\bms_1)(\Y_{0}^{\prime}\Y_0^{\prime\top}-\bms_1) \right\}\mid \Y_0,\Y_{0}^{\prime}\right]\right)=o\{\tr^3(\bms_1^2)\}$.
	\end{itemize}
\end{lemma}

For $2\leq \tau \leq n-2$, let 
\begin{small}
	\begin{equation*}
		H(\tau)=\frac{1}{\tau(\tau-1)}\sum_{i\neq j}^{\tau}(\Y_{i}^{\top}\Y_{j})^2+\frac{1}{(n-\tau)(n-\tau-1)}\sum_{i,j=\tau+1,i\neq j}^{n}(\Y_{i}^{\top}\Y_{j})^2-\frac{2}{\tau(n-\tau)}\sum_{i=1}^{\tau}\sum_{j=\tau+1}^{n}(\Y_{i}^{\top}\Y_{j})^2.
	\end{equation*}
\end{small}
The next lemma bounds the variance of $H(\tau)$. 
\begin{lemma}\label{lemma_two_sample_cov}
	Suppose Assumptions~1--3 and Assumption~5 hold. 
	Then,
	\begin{align*}
		\var\{H(\tau)\}\leq & C_1\left[\left(\frac{1}{\tau}+\frac{1}{n-\tau}\right)\tr\{(\bms_1+\bms_n)^2\}\right]^2\{1+o(1)\}\\
		&+C_2\left(\frac{1}{\tau}+\frac{1}{n-\tau}\right)\tr\{(\bms_1+\bms_n)^2\}\tr\{(\bms_1-\bms_n)^2\}\{1+o(1)\}
	\end{align*}
	holds uniformly for $2\leq \tau \leq n-2$, where $C_1,C_2>0$ are two constants.
\end{lemma}

The lemma below bounds $M_n(\tau)$ and $V_n(\tau)$. 
\begin{lemma}\label{lemma_M_tau_V_tau}
	Suppose Assumptions~1--3 and Assumption~5 hold. 
	Then,
	\begin{equation*}
		M_n(\tau)= \frac{(n-\tau^{\ast})(n-\tau^{\ast}-1)}{(n-\tau)(n-\tau-1)}\|\bm\delta\|^2 +O_p[n^{-1}\tr^{1/2}\{(\bms_1+\bms_n)^2\}]+O_p[\{n^{-1}\bm\delta(\bms_1+\bms_n) \bm\delta\}^{1/2}],
	\end{equation*}
	and
	\begin{align*}
		V_n(\tau)=&\frac{(n-\tau^{\ast})(n-\tau^{\ast}-1)}{(n-\tau)(n-\tau-1)}\tr\{(\bms_1-\bms_n)^2\}\\
		&+O_p[n^{-1}\tr\{(\bms_1+\bms_n)^2\}]+O_p\left([ n^{-1}\tr\{(\bms_1+\bms_n)^2\}\tr\{ (\bms_1-\bms_n)^2 \}]^{1/2}\right) \nonumber \\ 
		&+o_p\left[n^{-1/2}\tr^{3/4}\{(\bms_{1}+\bms_n)^2\}\|\bm\delta\|\right]+o_p\left[\|\bm\delta\|^2\sqrt{\tr\{(\bms_1+\bms_n)^2\}} \right]+O_p(\|\bm\delta\|^4)
	\end{align*}
	hold uniformly for all $\tau\in [\lambda_n,\tau^{\ast}]$, where $\bm\delta=\bm\mu_n-\bm\mu_1$.
\end{lemma}

The following lemma establishes the asymptotic normality of $M_n$ and $V_n$ individually.
\begin{lemma}\label{lemma_normal}
	Suppose Assumptions~1--2 hold. Under $\bH_0$, $\{\var_0(M_n)\}^{-1/2}M_n \rightarrow N(0,1)$ in distribution, as $n, p\rightarrow \infty$; if Assumption~3 also holds, then $\{\var_0(V_n)\}^{-1/2}V_n \rightarrow N(0,1)$ in distribution.
\end{lemma}

The following are two necessary lemmas.
\begin{lemma}[$L_p$ maximum inequality]\label{lemma_maximum_inequality}
	Let $z_i$ be a martingale. Let $z_n^{\ast}=\max_{1\leq i\leq n} |z_i|$. 
	Then, for $1\leq p<\infty$, 
	\begin{equation*}
		\E(|z_n^{\ast}|^p)\leq \left( \frac{p}{p-1} \right)^p \E(|z_n|^p).
	\end{equation*}
	Particularly, if $\E(z_i)=0$ and $\var(z_i)<\infty$, it holds that
	\begin{equation*}
		\Pr(z_n^{\ast}>x) \leq \frac{4\var(z_n)}{x^2}.
	\end{equation*}
\end{lemma}

\begin{lemma}[Mill’s ratio inequality]\label{lemma_normal_inequality}
	Let $\Phi(x)$ and $\phi(x)$ be the distribution function and the density function of the standard normal distribution, respectively.
	For any $x>0$,
	\begin{equation*}
		\frac{x}{\phi(x)}\leq \frac{1}{1-\Phi(x)}\leq  \frac{x}{\phi(x)}\frac{1+x^2}{x^2}.
	\end{equation*}
\end{lemma}

\section{Proof of Lemmas}

In this section, we give the proofs of Lemmas~\ref{lemma_basic_results}--\ref{lemma_normal}.

\subsection{Proof of Lemma~\ref{lemma_basic_results}}

For (1), the proof can be found in the proof of S1.11 in \cite{wang2022inference}.
For (2), the proof of $\E\{(\Y_{0}^{\prime\top}\Y_0)^6\}=O\{\tr^3(\bms_1^2)\}$ can be found in Lemma~S1.3 in \cite{wang2022inference}.
Using a similar proof technique, below we show $\E\{(\Y_{0}^{\top}\Y_{0}^{\prime})^8\}=O\{\tr^4(\bms_1^2)\}$.
Let $\mA$ be any disjoint partition over the set $\{l_1,\ldots,l_8\}$ such that for any $\mB\in \mA$, the cardinality $|\mB|\neq 1$.
We have
\begin{align*}
	&\E(\Y_{0}^{\prime\top}\Y_0)^8\\
	=&\sum_{l_1,\ldots,l_8=1}^{p}\{\E (y_{0,l_1}y_{0,l_2}y_{0,l_3}y_{0,l_4}y_{0,l_5}y_{0,l_6}y_{0,l_7}y_{0,l_8})\}^2\\
	=&\sum_{l_1,\ldots,l_8=1}^{p}\left\{\sum_{\mA}\prod_{\mB\in \mA}\cum(y_{0,l_k},l_k\in \mB)\right\}^2\\
	\leq &  C\sum_{l_1,\ldots,l_8=1}^{p} \sum_{\mA}\prod_{\mB\in \mA} \cum^2(y_{0,l_k},l_k\in \mB)\\
	\leq & C\sum_{\mA}\prod_{\mB\in \mA}\left\{ 
	\sum_{l_k\in \mB, l_k=1}^{p}\cum^2(y_{0,l_k},l_k\in \mB)\right\}\\
	\leq & C\sum_{\mA}\prod_{\mB\in \mA}\|\bms_1\|_{F}^{|\mB|}\\
	\leq & C\tr^4(\bms_1^2),
\end{align*}
where $C>0$ is a constant that varies line by line. Note that 
\begin{align*}
	&\E\{(\Y_0^{\top}\bms\Y_0)^4\}\\
	=&\sum_{l_1,\ldots,l_8=1}^{p}\E (y_{0,l_1}y_{0,l_2}y_{0,l_3}y_{0,l_4}y_{0,l_5}y_{0,l_6}y_{0,l_7}y_{0,l_8}\bmss_{l_1,l_2}\bmss_{l_3,l_4}\bmss_{l_5,l_6}\bmss_{l_7,l_8})\\
	\leq&\left[\sum_{l_1,\ldots,l_8=1}^{p}\{\E(y_{0,l_1}y_{0,l_2}y_{0,l_3}y_{0,l_4}y_{0,l_5}y_{0,l_6}y_{0,l_7}y_{0,l_8})\}^2\sum_{l_1,\ldots,l_8=1}^{p}(\bmss_{1,l_1,l_2}^2\bmss_{1,l_3,l_4}^2\bmss_{1,l_5,l_6}^2\bmss_{1,l_7,l_8}^2) \right]^{1/2}\\
	\leq & C\tr^4(\bms_1^2),
\end{align*}
where $C>0$ is a constant. Similarly, $\E\{(\Y_0^{\top}\bms\Y_0)^3\}=O\{\tr^3(\bms_1^2)\}$. 

Using the facts 
\begin{equation*}
	\sum_{l_1,l_2,l_3,l_4=1}^{p}\Sigma_{1,l_1,l_2}^2\Sigma_{1,l_3,l_4}^2=\tr^2(\bms_1^2),
\end{equation*}
\begin{equation*}
	\sum_{l_1,l_2,l_3,l_4=1}^{p}\Sigma_{1,l_1,l_2}\Sigma_{1,l_3,l_4}\Sigma_{1,l_1,l_4}\Sigma_{1,l_2,l_3}=\sum_{l_1,l_3=1}^{p}\left(\sum_{l_2=1}^{p}{\Sigma_{1,l_1,l_2}\Sigma_{1,l_2,l_3}}\right)^2=\tr(\bms_1^4),
\end{equation*}
\begin{equation*}
	\sum_{l_1,l_2,l_3,l_4=1}^{p}\cum^2(y_{0,l_1},y_{0,l_2},y_{0,l_3},y_{0,l_4})=o\{\tr^2(\bms_1^2)\},
\end{equation*}
and the Cauchy--Schwarz inequality, we have
\begin{align*}
	&\var\{(\Y_{0}^{\top}\Y_{0}^{\prime})^2\}=\E\{(\Y_{0}^{\top}\Y_{0}^{\prime})^4\}-\tr^2(\bms_1^2)\\
	= & \E\left(\sum_{l_1,l_2,l_3,l_4=1}^{p}y_{0,l_1}y_{0,l_1}^{\prime}y_{0,l_2}y_{0,l_2}^{\prime}y_{0,l_3}y_{0,l_3}^{\prime}y_{0,l_4}y_{0,l_4}^{\prime}\right)-\tr^2(\bms_1^2)\\
	=&\sum_{l_1,l_2,l_3,l_4=1}^{p}\E^2(y_{0,l_1}y_{0,l_2}y_{0,l_3}y_{0,l_4})-\tr^2(\bms_1^2)\\
	=&\sum_{l_1,l_2,l_3,l_4=1}^{p}\left\{\Sigma_{1,l_1,l_2}\Sigma_{1,l_3,l_4}+\Sigma_{1,l_1,l_3}\Sigma_{1,l_2,l_4}+\Sigma_{1,l_1,l_4}\Sigma_{1,l_2,l_3}+\cum(y_{0,l_1},y_{0,l_2},y_{0,l_3},y_{0,l_4})\right\}^2-\tr^2(\bms_1^2)\\
	= & 2\tr^2(\bms_1^2)+6\tr(\bms_1^4)+o\{\tr^2(\bms_1^2)\}.
\end{align*}

Observe that
\begin{align*}
	&\cov\{(\Y_{0}^{\top}\Y_{0}^{\prime})^2,(\Y_{0}^{\top}\Y_{0}^{\prime\prime})^2\}\\
	=&\E\{(\Y_{0}^{\top}\Y_{0}^{\prime})^2(\Y_{0}^{\top}\Y_{0}^{\prime\prime})^2\}-\tr^2(\bms_1^2)\\
	=&\E(\Y_{0}^{\top}\bms_1\Y_{0}\Y_{0}^{\top}\bms_1\Y_{0})-\tr^2(\bms_1^2)\\
	=&\sum_{l_1,l_2,l_3,l_4=1}^{p}\E(y_{0,l_1}y_{0,l_2}y_{0,l_3}y_{0,l_4}\Sigma_{1,l_1,l_2}\Sigma_{1,l_3,l_4})-\tr^2(\bms_1^2)\\
	=&\sum_{l_1,l_2,l_3,l_4=1}^{p}\left[\left\{\Sigma_{1,l_1,l_2}\Sigma_{1,l_3,l_4}+\Sigma_{1,l_1,l_3}\Sigma_{1,l_2,l_4}+\Sigma_{1,l_1,l_4}\Sigma_{1,l_2,l_3}+\cum(y_{0,l_1},y_{0,l_2},y_{0,l_3},y_{0,l_4})\right\}\Sigma_{1,l_1,l_2}\Sigma_{1,l_3,l_4}\right]\\
	&-\tr^2(\bms_1^2)\\
	=&2\tr(\bms_1^4)+o\{\tr^2(\bms_1^2)\}.
\end{align*}
From the $C_r$ inequality,
\begin{align*}
	&\E\left[\tr^4\left\{ 
	(\Y_0\Y_0^{\top}-\bms_1)(\Y_{0}^{\prime}\Y_0^{\prime\top}-\bms_1) \right\}\right]\\
	=&\E\left\{(\Y_0^{\top}\Y_{0}^{\prime})^2-\Y_{0}^{\top}\bms_1\Y_{0}-\Y_{0}^{\prime\top}\bms_1\Y_{0}^{\prime}+\tr(\bms_1^2)\right\}^4\\
	\leq &C\{\E\{(\Y_0^{\top}\Y_{0}^{\prime})^8\}+\E\{(\Y_{0}^{\top}\bms_1\Y_{0})^4\}+\E\{(\Y_{0}^{\prime\top}\bms_1\Y_{0}^{\prime})^4\}+\tr^4(\bms_1^2)\}\\
	\leq & C\tr^4(\bms_1^2),
\end{align*}
where $C>0$ is a constant that varies line by line. 
By Jensen's inequality,
\begin{align*}
	&\var\left(\E\left[\tr^2\left\{ 
	(\Y_{0}^{\prime}\Y_0^{\prime\top}-\bms_1)(\Y_0\Y_0^{\top}-\bms_1) \right\}\mid \Y_0\right]\right)\\
	\leq& \E\left[\tr^4\left\{ 
	(\Y_0\Y_0^{\top}-\bms_1)(\Y_{0}^{\prime}\Y_0^{\prime\top}-\bms_1) \right\}\right]\leq C\tr^4(\bms_1^2).
\end{align*}
Next, we will show 
\begin{small}
	\begin{equation*}
		\var\left(\E\left[\tr\left\{ 
		(\Y_0^{\prime\prime}\Y_0^{\prime\prime\top}-\bms_1)(\Y_{0}\Y_0^{\top}-\bms_1) \right\} \tr\left\{ 
		(\Y_0^{\prime\prime}\Y_0^{\prime\prime\top}-\bms_1)(\Y_{0}^{\prime}\Y_0^{\prime\top}-\bms_1) \right\}\mid \Y_0,\Y_{0}^{\prime}\right]\right\}=o\{\tr^4(\bms_1^2)\}.
	\end{equation*}
\end{small}
Note that
\begin{align*}
	&\E\left[\tr\left\{ 
	(\Y_0^{\prime\prime}\Y_0^{\prime\prime\top}-\bms_1)(\Y_{0}\Y_0^{\top}-\bms_1) \right\} \tr\left\{ 
	(\Y_0^{\prime\prime}\Y_0^{\prime\prime\top}-\bms_1)(\Y_{0}^{\prime}\Y_0^{\prime\top}-\bms_1) \right\}\mid \Y_0,\Y_{0}^{\prime}\right]\\
	=&\E\{ (\Y_0^{\prime\prime\top}\Y_0)^2(\Y_0^{\prime\prime\top}\Y_{0}^{\prime})^2\mid \Y_0,\Y_{0}^{\prime}\}-(\Y_0^{\prime\top}\bms_1\Y_{0}^{\prime})(\Y_0^{\top}\bms_1\Y_0)\\
	&+\Y_0^{\top}\bms_1\Y_0\tr(\bms_1^2)+\Y_0^{\prime\top}\bms_1\Y_{0}^{\prime}\tr(\bms_1^2)-\E\{(\Y_0^{\prime\prime\top}\Y_0)^2\Y_0^{\prime\prime\top}\bms_1\Y_0^{\prime\prime}\mid \Y_0 \}\\
	&-\E\{(\Y_0^{\prime\prime\top}\Y_{0}^{\prime})^2\Y_0^{\prime\prime\top}\bms_1\Y_0^{\prime\prime}\mid \Y_{0}^{\prime} \}
	+\E(\Y_{0}^{\prime\prime\top}\bms_1\Y_{0}^{\prime\prime})^2-\tr^2(\bms_1^2).
\end{align*}
By $\var(\Y_0^{\top}\bms_1\Y_0)=o\{\tr^2(\bms_1^2)\}$, obviously, we have
\begin{equation*}
	\var\{\Y_0^{\top}\bms_1\Y_0\tr(\bms_1^2)\}=o\{\tr^4(\bms_1^2)\}
\end{equation*}
and 
\begin{equation*}
	\var\{(\Y_0^{\top}\bms_1\Y_0)(\Y_0^{\prime\top}\bms_1\Y_{0}^{\prime})\}=o\{\tr^4(\bms_1^2)\}.
\end{equation*}    
Note that
\begin{small}
	\begin{align*}
		&\E\{ (\Y_0^{\prime\prime\top}\Y_0)^2(\Y_0^{\prime\prime\top}\Y_{0}^{\prime})^2\mid \Y_0,\Y_{0}^{\prime}\}\\
		=&\sum_{l_1,l_2,l_3,l_4=1}^{p}\left\{y_{0,l_1}y_{0,l_2}y_{0,l_3}^{\prime}y_{0,l_4}^{\prime} \E(y_{0,l_1}^{\prime\prime}y_{0,l_2}^{\prime\prime}y_{0,l_3}^{\prime\prime}y_{0,l_4}^{\prime\prime}) \right\}\\
		=&\sum_{l_1,l_2,l_3,l_4=1}^{p}\left[y_{0,l_1}y_{0,l_2}y_{0,l_3}^{\prime}y_{0,l_4}^{\prime}\{\bmss_{1,l_1,l_2}\bmss_{1,l_3,l_4}+\bmss_{1,l_1,l_3}\bmss_{1,l_2,l_4}+\bmss_{1,l_1,l_4}\bmss_{1,l_2,l_3}+\cum(y_{0,l_1},y_{0,l_2},y_{0,l_3},y_{0,l_4})\}\right].
	\end{align*}
\end{small}
Then,
\begin{align*}
	&\E\left[\E\{ (\Y_0^{\prime\prime\top}\Y_0)^2(\Y_0^{\prime\prime\top}\Y_{0}^{\prime})^2\mid \Y_0,\Y_{0}^{\prime}\}\right]^2\\
	=&\sum_{l_1,\ldots,l_8}^{p}\bigg[\E(y_{0,l_1}y_{0,l_2}y_{0,l_5}y_{0,l_6})\E(y_{0,l_3}^{\prime}y_{0,l_4}^{\prime}y_{0,l_7}^{\prime}y_{0,l_8}^{\prime}) \\
	& \times\{\bmss_{1,l_1,l_2}\bmss_{1,l_3,l_4}+\bmss_{1,l_1,l_3}\bmss_{1,l_2,l_4}+\bmss_{1,l_1,l_4}\bmss_{1,l_2,l_3}+\cum(y_{0,l_1},y_{0,l_2},y_{0,l_3},y_{0,l_4})\}\\
	&\times\{\bmss_{1,l_5,l_6}\bmss_{1,l_7,l_8}+\bmss_{1,l_5,l_7}\bmss_{1,l_6,l_8}+\bmss_{1,l_5,l_8}\bmss_{1,l_6,l_7}+\cum(y_{0,l_5},y_{0,l_6},y_{0,l_7},y_{0,l_8})\}\bigg]\\
	=&\sum_{l_1,\ldots,l_8}^{p}\bigg[ \{\bmss_{1,l_1,l_2}\bmss_{1,l_5,l_6}+\bmss_{1,l_1,l_5}\bmss_{1,l_2,l_6}+\bmss_{1,l_1,l_6}\bmss_{1,l_2,l_5}+\cum(y_{0,l_1},y_{0,l_2},y_{0,l_5},y_{0,l_6})\} \\
	& \times  \{\bmss_{1,l_3,l_4}\bmss_{1,l_7,l_8}+\bmss_{1,l_3,l_7}\bmss_{1,l_4,l_8}+\bmss_{1,l_3,l_8}\bmss_{1,l_4,l_7}+\cum(y_{0,l_3},y_{0,l_4},y_{0,l_7},y_{0,l_8})\} \\
	& \times\{\bmss_{1,l_1,l_2}\bmss_{1,l_3,l_4}+\bmss_{1,l_1,l_3}\bmss_{1,l_2,l_4}+\bmss_{1,l_1,l_4}\bmss_{1,l_2,l_3}+\cum(y_{0,l_1},y_{0,l_2},y_{0,l_3},y_{0,l_4})\}\\
	&\times\{\bmss_{1,l_5,l_6}\bmss_{1,l_7,l_8}+\bmss_{1,l_5,l_7}\bmss_{1,l_6,l_8}+\bmss_{1,l_5,l_8}\bmss_{1,l_6,l_7}+\cum(y_{0,l_5},y_{0,l_6},y_{0,l_7},y_{0,l_8})\}\bigg]\\
	=&\tr^4(\bms_1^2)+o\{\tr^4(\bms_1^2)\},
\end{align*}
by using the facts
\begin{align*}
	\sum_{l_1,\ldots,l_8}^{p}\bmss_{1,l_1,l_2}\bmss_{1,l_5,l_6}\bmss_{1,l_3,l_4}\bmss_{1,l_7,l_8}\bmss_{1,l_1,l_2}\bmss_{1,l_3,l_4}\bmss_{1,l_5,l_6}\bmss_{1,l_7,l_8}&=\tr^4(\bms_1^2),\\
	\sum_{l_1,\ldots,l_8}^{p}\bmss_{1,l_1,l_2}\bmss_{1,l_5,l_6}\bmss_{1,l_3,l_7}\bmss_{1,l_4,l_8}\bmss_{1,l_1,l_2}\bmss_{1,l_3,l_4}\bmss_{1,l_5,l_6}\bmss_{1,l_7,l_8}&=\tr^2(\bms_1^2)\tr(\bms_1^4),\\
	\sum_{l_1,\ldots,l_8}^{p}\bmss_{1,l_1,l_2}\bmss_{1,l_5,l_6}\bmss_{1,l_3,l_4}\bmss_{1,l_7,l_8}\bmss_{1,l_1,l_5}\bmss_{1,l_2,l_6}\bmss_{1,l_3,l_7}\bmss_{1,l_4,l_8}&=\tr^2(\bms_1^4),\\
	\sum_{l_1,\ldots,l_8}^{p}\bmss_{1,l_1,l_2}\bmss_{1,l_5,l_7}\bmss_{1,l_3,l_4}\bmss_{1,l_6,l_8}\bmss_{1,l_1,l_2}\bmss_{1,l_5,l_6}\bmss_{1,l_3,l_7}\bmss_{1,l_4,l_8}&=\tr(\bms_1^2)\tr(\bms_1^6)\leq \tr^2(\bms_1^2)\tr(\bms_1^4) ,\\
	\sum_{l_1,\ldots,l_8}^{p}\bmss_{1,l_1,l_2}\bmss_{1,l_5,l_7}\bmss_{1,l_3,l_4}\bmss_{1,l_6,l_8}\bmss_{1,l_1,l_5}\bmss_{1,l_2,l_6}\bmss_{1,l_3,l_7}\bmss_{1,l_4,l_8}&=\tr(\bms_1^8)\leq \tr^2(\bms_1^4),\\
	\sum_{l_1,l_2,l_3,l_4=1}^{p}\cum^2(y_{0,l_1},y_{0,l_2},y_{0,l_3},y_{0,l_4})&=o\{\tr^2(\bms_1^2)\},
\end{align*}
and the Cauchy--Schwarz inequality. Thus $\var[\E\{ (\Y_0^{\prime\prime\top}\Y_0)^2(\Y_0^{\prime\prime\top}\Y_{0}^{\prime})^2\mid \Y_0,\Y_{0}^{\prime}\}]=o\{\tr^4(\bms_1^2)\}$.
Since 
\begin{equation*}
	\E[\E\{ (\Y_0^{\prime\prime\top}\Y_0)^2(\Y_0^{\prime\prime\top}\Y_{0}^{\prime})^2\mid \Y_0,\Y_{0}^{\prime}\}\mid \Y_0]=\E\{(\Y_0^{\prime\prime\top}\Y_0)^2\Y_0^{\prime\prime\top}\bms_1\Y_0^{\prime\prime}\mid \Y_0 \},
\end{equation*}
we have
\begin{equation*}
	\var[\E\{(\Y_0^{\prime\prime\top}\Y_0)^2\Y_0^{\prime\prime\top}\bms_1\Y_0^{\prime\prime}\mid \Y_0 \}] \leq\var[ \E\{ (\Y_0^{\prime\prime\top}\Y_0)^2(\Y_0^{\prime\prime\top}\Y_{0}^{\prime})^2\mid \Y_0,\Y_{0}^{\prime}\} ]=o\{\tr^4(\bms_1^2)\}.
\end{equation*}
Putting the results together, we conclude
\begin{equation*}
	\var(\E[\tr\{ 
	(\Y_0^{\prime\prime}\Y_0^{\prime\prime\top}-\bms_1)(\Y_{0}\Y_0^{\top}-\bms_1) \} \tr\{ 
	(\Y_0^{\prime\prime}\Y_0^{\prime\prime\top}-\bms_1)(\Y_{0}^{\prime}\Y_0^{\prime\top}-\bms_1) \mid \Y_0,\Y_{0}^{\prime}])=o\{\tr^4(\bms_1^2)\}.
\end{equation*}

For (3),
\begin{align*}
	&\var(\E[\Y_{0}^{\prime\top}\Y_0\tr\{ 
	(\Y_{0}^{\prime}\Y_0^{\prime\top}-\bms_1)(\Y_0\Y_0^{\top}-\bms_1)\}\mid \Y_0])\\
	\leq & \E[\Y_0^{\prime\top}\Y_0\tr\{ 
	(\Y_{0}^{\prime}\Y_0^{\prime\top}-\bms)(\Y_0\Y_0^{\top}-\bms)\}]^2 \\
	\leq & C[\E\{(\Y_0^{\prime\top}\Y_0)^6\}+\E\{(\Y_0^{\top}\bms_1\Y_0)^3\}+\E\{(\Y_{0}^{\prime\top}\bms_1\Y_{0}^{\prime})^3\}+\tr^2(\bms_1^2)\E\{(\Y_0^{\prime\top}\Y_0)^2\}]\\
	=&O\{\tr^3(\bms_1^2)\},
\end{align*}
where $C>0$ is a constant.
Note that
\begin{align*}
	&\E[
	\Y_0^{\prime\prime\top}\Y_0 \tr\{ 
	(\Y_0^{\prime\prime}\Y_0^{\prime\prime\top}-\bms_1)(\Y_{0}^{\prime}\Y_0^{\prime\top}-\bms_1)\}\mid \Y_0,\Y_{0}^{\prime}]\\
	=&\E\{ \Y_{0}^{\prime\prime\top}\Y_0(\Y_0^{\prime\prime\top}\Y_{0}^{\prime})^2\mid \Y_0,\Y_{0}^{\prime}\}-\E\{\Y_0^{\prime\prime\top}\Y_0\Y_0^{\prime\prime\top}\bms_1\Y_{0}^{\prime\prime}\mid \Y_0 \},
\end{align*}
and that
\begin{align*}
	&\E\left[\E\{ \Y_{0}^{\prime\prime\top}\Y_0(\Y_0^{\prime\prime\top}\Y_{0}^{\prime})^2\mid \Y_0,\Y_{0}^{\prime}\}\right]^2\\
	=&\sum_{l_1,\ldots,l_6=1}^{p}\E(y_{0,l_1}y_{0,l_4}) \E(y_{0,l_2}^{\prime}y_{0,l_3}^{\prime}y_{0,l_5}^{\prime}y_{0,l_6}^{\prime}) \cum(y_{0,l_1}^{\prime\prime},y_{0,l_2}^{\prime\prime},y_{0,l_3}^{\prime\prime}) \cum(y_{0,l_4}^{\prime\prime},y_{0,l_5}^{\prime\prime},y_{0,l_6}^{\prime\prime})\\
	=&\sum_{l_1,\ldots,l_6=1}^{p}\bmss_{1,l_1,l_4}\cum(y_{0,l_1}^{\prime\prime},y_{0,l_2}^{\prime\prime},y_{0,l_3}^{\prime\prime})\cum(y_{0,l_4}^{\prime\prime},y_{0,l_5}^{\prime\prime},y_{0,l_6}^{\prime\prime})\\
	&\times\{\bmss_{1,l_2,l_3}\bmss_{1,l_5,l_6}+\bmss_{1,l_2,l_5}\bmss_{1,l_3,l_6}+\bmss_{1,l_2,l_6}\bmss_{1,l_3,l_5}+\cum(y_{0,l_2}^{\prime},y_{0,l_3}^{\prime},y_{0,l_5}^{\prime},y_{0,l_6}^{\prime})\}\\
	\leq &\left\{\sum_{l_1,\ldots,l_6=1}^{p}\cum^2(y_{0,l_1}^{\prime\prime},y_{0,l_2}^{\prime\prime},y_{0,l_3}^{\prime\prime})\cum^2(y_{0,l_4}^{\prime\prime},y_{0,l_5}^{\prime\prime},y_{0,l_6}^{\prime\prime})\right\}^{1/2}\\
	&\times \left[\sum_{l_1,\ldots,l_6=1}^{p}\bmss_{1,l_1,l_4}^2 \{\bmss_{1,l_2,l_3}\bmss_{1,l_5,l_6}+\bmss_{1,l_2,l_5}\bmss_{1,l_3,l_6}+\bmss_{1,l_2,l_6}\bmss_{1,l_3,l_5}+\cum(y_{0,l_2}^{\prime},y_{0,l_3}^{\prime},y_{0,l_5}^{\prime},y_{0,l_6}^{\prime})\}^2 \right]^{1/2}\\
	=&o\{\tr^3(\bms_1^2)\}.
\end{align*}
Since 
\begin{equation*}
	\E[\E\{ \Y_{0}^{\prime\prime\top}\Y_0(\Y_0^{\prime\prime\top}\Y_{0}^{\prime})^2\mid \Y_0,\Y_{0}^{\prime}\}\mid \Y_0 ]=\E(\Y_0^{\prime\prime\top}\Y_0\Y_0^{\prime\prime\top}\bms_1\Y_{0}^{\prime\prime}\mid \Y_0 ),
\end{equation*}
\begin{equation*}
	\var\{\E(\Y_0^{\prime\prime\top}\Y_0\Y_0^{\prime\prime\top}\bms_1\Y_{0}^{\prime\prime}\mid \Y_0 )\} \leq \var[\E\{ \Y_{0}^{\prime\prime\top}\Y_0(\Y_0^{\prime\prime\top}\Y_{0}^{\prime})^2\mid \Y_0,\Y_{0}^{\prime}\}]=o\{\tr^3(\bms_1^2)\},
\end{equation*}
we obtain
\begin{equation*}
	\var(\E[
	\Y_0^{\prime\prime\top}\Y_0 \tr\{ 
	(\Y_0^{\prime\prime}\Y_0^{\prime\prime\top}-\bms_1)(\Y_{0}^{\prime}\Y_0^{\prime\top}-\bms_1) \}\mid \Y_0,\Y_{0}^{\prime}])=o\{\tr^3(\bms_1^2)\}. 
\end{equation*}
The lemma follows. \done

\subsection{Proof of Lemma~\ref{lemma_two_sample_cov}}

For $\tau\leq \tau^{\ast}$,  we can decompose $H(\tau)$ as
\begin{align*}
	H(\tau)=&\frac{\sum_{i_1\neq i_2}^{\tau} (\Y_{i_1}^{\top}\Y_{i_2})^2  }{\tau(\tau-1)}+\frac{\sum_{i_1\neq i_2,\tau+1}^{\tau^{\ast}} (\Y_{i_1}^{\top}\Y_{i_2})^2  }{(n-\tau)(n-\tau-1)}-\frac{2\sum_{i_1=1}^{\tau}\sum_{i_2=\tau+1}^{\tau^{\ast}} (\Y_{i_1}^{\top}\Y_{i_2})^2  }{\tau(n-\tau)}\\
	&+\frac{\sum_{i_1\neq i_2,\tau^{\ast}+1}^{n} (\Y_{i_1}^{\top}\Y_{i_2})^2  }{(n-\tau)(n-\tau-1)}+\frac{2\sum_{i_1=\tau+1}^{\tau^{\ast}} \sum_{i_2=\tau^{\ast}+1}^{n} (\Y_{i_1}^{\top}\Y_{i_2})^2  }{(n-\tau)(n-\tau-1)}-\frac{2\sum_{i_1=1}^{\tau} \sum_{i_2=\tau^{\ast}+1}^{n} (\Y_{i_1}^{\top}\Y_{i_2})^2  }{\tau(n-\tau)}.
\end{align*}
To more intuitively distinguish the data before and after the changepoint, we denote $\Y_{\tau^{\ast}+i}$ by $\Z_i$ for $i=1,\ldots,n-\tau^{\ast}$.
Let $\Z_0=(z_{0,1},\ldots,z_{0,p})^{\top}$ be an independent copy of $\Z_1,\ldots,\Z_{n-\tau^{\ast}}$. 
One may verify the following facts that
\begin{equation*}
	\var\left\{ \frac{\sum_{i_1\neq i_2}^{\tau} (\Y_{i_1}^{\top}\Y_{i_2})^2  }{\tau(\tau-1)}  \right\}
	=\frac{2\var\{(\Y_{1}^{\top}\Y_2)^2\}+4(\tau-2)
		\cov\{(\Y_{1}^{\top}\Y_2)^2,(\Y_{1}^{\top}\Y_3)^2\} }{\tau(\tau-1)},
\end{equation*}

\begin{align*}
	&\var\left\{ \frac{\sum_{i_1\neq i_2,\tau+1}^{\tau^{\ast}} (\Y_{i_1}^{\top}\Y_{i_2})^2  }{(n-\tau)(n-\tau-1)} \right\}\\
	=&(\tau^{\ast}-\tau)(\tau^{\ast}-\tau-1)\frac{2\var\{(\Y_{1}^{\top}\Y_2)^2\}+4(\tau^{\ast}-\tau-2)
		\cov\{(\Y_{1}^{\top}\Y_2)^2,(\Y_{1}^{\top}\Y_3)^2\} }{(n-\tau)^2(n-\tau-1)^2},
\end{align*}

\begin{align*}
	&\var\left\{ \frac{2\sum_{i_1=1}^{\tau}\sum_{i_2=\tau+1}^{\tau^{\ast}} (\Y_{i_1}^{\top}\Y_{i_2})^2  }{\tau(n-\tau)} \right\}\\
	=&4\tau(\tau^{\ast}-\tau)\frac{\var\{(\Y_{1}^{\top}\Y_2)^2\}+(\tau^{\ast}-2)
		\cov\{(\Y_{1}^{\top}\Y_2)^2,(\Y_{1}^{\top}\Y_3)^2\} }{\tau^2(n-\tau)^2},
\end{align*}    

\begin{align*}    
	&\var\left\{\frac{\sum_{i_1\neq i_2,\tau^{\ast}+1}^{n} (\Y_{i_1}^{\top}\Y_{i_2})^2  }{(n-\tau)(n-\tau-1)}\right\}\\
	=&(n-\tau^{\ast})(n-\tau^{\ast}-1)\frac{2\var\{(\Z_{1}^{\top}\Z_2)^2\}+4(n-\tau^{\ast}-2)
		\cov\{(\Z_{1}^{\top}\Z_2)^2,(\Z_{1}^{\top}\Z_3)^2\} }{(n-\tau)^2(n-\tau-1)^2},
\end{align*}

\begin{align*}
	&\var\left\{ \frac{2\sum_{i_1=\tau+1}^{\tau^{\ast}} \sum_{i_2=\tau^{\ast}+1}^{n} (\Y_{i_1}^{\top}\Y_{i_2})^2  }{(n-\tau)(n-\tau-1)} \right\}\\
	=&4\frac{(\tau^{\ast}-\tau)(n-\tau^{\ast})\var\{(\Y_{1}^{\top}\Z_1)^2\}}{(n-\tau)^2(n-\tau-1)^2}+4\frac{(\tau^{\ast}-\tau)(n-\tau^{\ast})(n-\tau^{\ast}-1)\cov\{(\Y_{1}^{\top}\Z_1)^2,(\Y_{1}^{\top}\Z_2)^2\}}{(n-\tau)^2(n-\tau-1)^2}\\
	&+4\frac{(\tau^{\ast}-\tau)(\tau^{\ast}-\tau-1)(n-\tau^{\ast})\cov\{(\Y_{1}^{\top}\Z_1)^2,(\Y_{2}^{\top}\Z_1)^2\}}{(n-\tau)^2(n-\tau-1)^2},\\
\end{align*}

\begin{align*}
	&\var\left\{ \frac{2\sum_{i_1=1}^{\tau} \sum_{i_2=\tau^{\ast}+1}^{n} (\Y_{i_1}^{\top}\Y_{i_2})^2  }{\tau(n-\tau)}\right\}\\
	=& 4\frac{\tau(n-\tau^{\ast})\var\{(\Y_{1}^{\top}\Z_1)^2\}}{\tau^2(n-\tau)^2}+4\frac{\tau(n-\tau^{\ast})(n-\tau^{\ast}-1)\cov\{(\Y_{1}^{\top}\Z_1)^2,(\Y_{1}^{\top}\Z_2)^2\}}{\tau^2(n-\tau)^2}\\
	&+4\frac{\tau(\tau-1)(n-\tau^{\ast})\cov\{(\Y_{1}^{\top}\Z_1)^2,(\Y_{2}^{\top}\Z_1)^2\}}{\tau^2(n-\tau)^2},
\end{align*}    

\begin{align*}
	&\cov\left\{ \frac{\sum_{i_1\neq i_2}^{\tau} (\Y_{i_1}^{\top}\Y_{i_2})^2  }{\tau(\tau-1)},-\frac{2\sum_{i_1=1}^{\tau}\sum_{i_2=\tau+1}^{\tau^{\ast}} (\Y_{i_1}^{\top}\Y_{i_2})^2  }{\tau(n-\tau)}\right\}\\
	=&-2\frac{\tau(\tau-1)2(\tau^{\ast}-\tau)\cov\{(\Y_{1}^{\top}\Y_2)^2,(\Y_{1}^{\top}\Y_3)^2\}}{\tau(\tau-1)\tau(n-\tau)},
\end{align*}

\begin{align*}
	&\cov\left\{ \frac{\sum_{i_1\neq i_2}^{\tau} (\Y_{i_1}^{\top}\Y_{i_2})^2  }{\tau(\tau-1)},-\frac{2\sum_{i_1=1}^{\tau} \sum_{i_2=\tau^{\ast}+1}^{n} (\Y_{i_1}^{\top}\Y_{i_2})^2  }{\tau(n-\tau)} \right\}\\
	=&-2\frac{\tau(\tau-1)2(n-\tau^{\ast})\cov\{(\Y_{1}^{\top}\Y_2)^2,(\Y_{1}^{\top}\Z_1)^2\}}{\tau(\tau-1)\tau(n-\tau)},
\end{align*}    

\begin{align*}
	&\cov\left\{  \frac{\sum_{i_1\neq i_2,\tau+1}^{\tau^{\ast}} (\Y_{i_1}^{\top}\Y_{i_2})^2  }{(n-\tau)(n-\tau-1)},-\frac{2\sum_{i_1=1}^{\tau}\sum_{i_2=\tau+1}^{\tau^{\ast}} (\Y_{i_1}^{\top}\Y_{i_2})^2  }{\tau(n-\tau)}  \right\}\\
	=&-2\frac{(\tau^{\ast}-\tau)(\tau^{\ast}-\tau-1)2\tau\cov\{(\Y_{1}^{\top}\Y_2)^2,(\Y_{1}^{\top}\Y_3)^2\}}{(n-\tau)(n-\tau-1)\tau(n-\tau)},
\end{align*}

\begin{align*}  
	&\cov\left\{  \frac{\sum_{i_1\neq i_2,\tau+1}^{\tau^{\ast}} (\Y_{i_1}^{\top}\Y_{i_2})^2  }{(n-\tau)(n-\tau-1)},\frac{2\sum_{i_1=\tau+1}^{\tau^{\ast}} \sum_{i_2=\tau^{\ast}+1}^{n} (\Y_{i_1}^{\top}\Y_{i_2})^2  }{(n-\tau)(n-\tau-1)} \right\}\\
	&=2\frac{(\tau^{\ast}-\tau)(\tau^{\ast}-\tau-1)2(n-\tau^{\ast})\cov\{(\Y_{1}^{\top}\Y_2)^2,(\Y_{1}^{\top}\Z_1)^2\}}{(n-\tau)^2(n-\tau-1)^2},\\
\end{align*}     

\begin{align*}
	&\cov\left\{ -\frac{2\sum_{i_1=1}^{\tau}\sum_{i_2=\tau+1}^{\tau^{\ast}} (\Y_{i_1}^{\top}\Y_{i_2})^2  }{\tau(n-\tau)}, \frac{2\sum_{i_1=\tau+1}^{\tau^{\ast}} \sum_{i_2=\tau^{\ast}+1}^{n} (\Y_{i_1}^{\top}\Y_{i_2})^2  }{(n-\tau)(n-\tau-1)} \right\}\\
	=&-4\frac{\tau(\tau^{\ast}-\tau)(n-\tau^{\ast})\cov\{(\Y_{1}^{\top}\Y_2)^2,(\Y_{1}^{\top}\Z_1)^2\}}{\tau(n-\tau)^2(n-\tau-1)},
\end{align*}

\begin{align*}
	&\cov\left\{ -\frac{2\sum_{i_1=1}^{\tau}\sum_{i_2=\tau+1}^{\tau^{\ast}} (\Y_{i_1}^{\top}\Y_{i_2})^2  }{\tau(n-\tau)}, -\frac{2\sum_{i_1=1}^{\tau} \sum_{i_2=\tau^{\ast}+1}^{n} (\Y_{i_1}^{\top}\Y_{i_2})^2  }{\tau(n-\tau)}  \right\}\\
	= & 4\frac{\tau(\tau^{\ast}-\tau)(n-\tau^{\ast})\cov\{(\Y_{1}^{\top}\Y_2)^2,(\Y_{1}^{\top}\Z_1)^2\}}{\tau^2(n-\tau)^2},
\end{align*} 

\begin{align*}
	&\cov\left\{ \frac{\sum_{i_1\neq i_2,\tau^{\ast}+1}^{n} (\Y_{i_1}^{\top}\Y_{i_2})^2  }{(n-\tau)(n-\tau-1)}, \frac{2\sum_{i_1=\tau+1}^{\tau^{\ast}} \sum_{i_2=\tau^{\ast}+1}^{n} (\Y_{i_1}^{\top}\Y_{i_2})^2  }{(n-\tau)(n-\tau-1)}\right\}\\
	=&2\frac{(n-\tau^{\ast})(n-\tau^{\ast}-1)2(\tau^{\ast}-\tau)\cov\{(\Z_{1}^{\top}\Z_2)^2,(\Y_{1}^{\top}\Z_1)^2\}}{(n-\tau)^2(n-\tau-1)^2},
\end{align*}    

\begin{align*}
	&\cov\left\{ \frac{\sum_{i_1\neq i_2,\tau^{\ast}+1}^{n} (\Y_{i_1}^{\top}\Y_{i_2})^2  }{(n-\tau)(n-\tau-1)}, -\frac{2\sum_{i_1=1}^{\tau} \sum_{i_2=\tau^{\ast}+1}^{n} (\Y_{i_1}^{\top}\Y_{i_2})^2  }{\tau(n-\tau)}\right\}\\
	= & -2\frac{(n-\tau^{\ast})(n-\tau^{\ast}-1)2\tau\cov\{(\Z_{1}^{\top}\Z_2)^2,(\Y_{1}^{\top}\Z_1)^2\}}{\tau(n-\tau)^2(n-\tau-1)}
\end{align*}    

\begin{align*}
	&\cov\left\{  \frac{2\sum_{i_1=\tau+1}^{\tau^{\ast}} \sum_{i_2=\tau^{\ast}+1}^{n} (\Y_{i_1}^{\top}\Y_{i_2})^2  }{(n-\tau)(n-\tau-1)},-\frac{2\sum_{i_1=1}^{\tau} \sum_{i_2=\tau^{\ast}+1}^{n} (\Y_{i_1}^{\top}\Y_{i_2})^2  }{\tau(n-\tau)}  \right\}\\
	=&-4\frac{(\tau^{\ast}-\tau)(n-\tau^{\ast})\tau \cov\{(\Y_{1}^{\top}\Z_1)^2,(\Y_{2}^{\top}\Z_1)^2\}}{\tau(n-\tau)^2(n-\tau-1)},
\end{align*}
and 
\begin{align*}
	&\var\{(\Y_{1}^{\top}\Y_2)^2\}=2\tr^2(\bms_1^2)\{1+o(1)\},\\
	&\var\{(\Z_{1}^{\top}\Z_2)^2\}=2\tr^2(\bms_n^2)\{1+o(1)\},\\
	&\var\{(\Y_{1}^{\top}\Z_1)^2\}=2\tr^2(\bms_1\bms_n)\{1+o(1)\},
\end{align*}
\begin{align*}
	&\cov\{(\Y_{1}^{\top}\Y_2)^2,(\Y_{1}^{\top}\Y_3)^2\}=2\tr(\bms_1^4)+\sum_{l_1,l_2,l_3,l_4}\left\{\cum(y_{0,l_1},y_{0,l_2},y_{0,l_3},y_{0,l_4})\bmss_{1,l_1,l_2}\bmss_{1,l_3,l_4}\right\},\\
	&\cov\{(\Z_{1}^{\top}\Z_2)^2,(\Z_{1}^{\top}\Z_3)^2\}=2\tr(\bms_n^4)+\sum_{l_1,l_2,l_3,l_4}\left\{\cum(z_{0,l_1},z_{0,l_2},z_{0,l_3},z_{0,l_4})\bmss_{n,l_1,l_2}\bmss_{n,l_3,l_4}\right\},\\
	&\cov\{(\Y_{1}^{\top}\Y_2)^2,(\Y_{1}^{\top}\Z_1)^2\}=2\tr(\bms_1^3\bms_n)+\sum_{l_1,l_2,l_3,l_4}\left\{\cum(y_{0,l_1},y_{0,l_2},y_{0,l_3},y_{0,l_4})\bmss_{1,l_1,l_2}\bmss_{n,l_3,l_4}\right\},\\
	&\cov\{(\Z_{1}^{\top}\Z_2)^2,(\Y_{1}^{\top}\Z_1)^2\}=2\tr(\bms_n^3\bms_1)+\sum_{l_1,l_2,l_3,l_4}\left\{\cum(z_{0,l_1},z_{0,l_2},z_{0,l_3},z_{0,l_4})\bmss_{1,l_1,l_2}\bmss_{n,l_3,l_4}\right\},\\
	&\cov\{(\Y_{1}^{\top}\Z_1)^2,(\Y_{1}^{\top}\Z_2)^2\}=2\tr(\bms_1\bms_n\bms_1\bms_n)+\sum_{l_1,l_2,l_3,l_4}\left\{\cum(y_{0,l_1},y_{0,l_2},y_{0,l_3},y_{0,l_4})\bmss_{n,l_1,l_2}\bmss_{n,l_3,l_4}\right\},\\
	&\cov\{(\Y_{1}^{\top}\Z_1)^2,(\Y_{2}^{\top}\Z_1)^2\}=2\tr(\bms_1\bms_n\bms_1\bms_n)+\sum_{l_1,l_2,l_3,l_4}\left\{\cum(z_{0,l_1},z_{0,l_2},z_{0,l_3},z_{0,l_4})\bmss_{1,l_1l_2}\bmss_{1,l_3l_4}\right\}.
\end{align*}
Moreover,
\begin{align*}
	&\sum_{l_1,l_2,l_3,l_4}\left\{\cum(y_{0,l_1},y_{0,l_2},y_{0,l_3},y_{0,l_4})\bmss_{1,l_1,l_2}\bmss_{1,l_3,l_4}\right\}\\
	&-2\sum_{l_1,l_2,l_3,l_4}\left\{\cum(y_{0,l_1},y_{0,l_2},y_{0,l_3},y_{0,l_4})\bmss_{1,l_1,l_2}\bmss_{n,l_3,l_4}\right\}\\
	&+\sum_{l_1,l_2,l_3,l_4}\left\{\cum(y_{0,l_1},y_{0,l_2},y_{0,l_3},y_{0,l_4})\bmss_{n,l_1,l_2}\bmss_{n,l_3,l_4}\right\}\\
	=&\sum_{l_1,l_2,l_3,l_4} \left\{\cum(y_{0,l_1},y_{0,l_2},y_{0,l_3},y_{0,l_4})(\bmss_{1,l_1,l_2}-\bmss_{n,l_1,l_2})(\bmss_{1,l_3,l_4}-\bmss_{n,l_3,l_4})\right\}\\
	\leq & \left[ \sum_{l_1,l_2,l_3,l_4}\cum^2(y_{0,l_1},y_{0,l_2},y_{0,l_3},y_{0,l_4})  \tr^2\{(\bms_1-\bms_n)^2\} \right]^{1/2}\\
	=&o\{\tr(\bms_1^2)\}\tr\{(\bms_1-\bms_n)^2\},
\end{align*}
\begin{align*}
	&\sum_{l_1,l_2,l_3,l_4}\left\{\cum(z_{0,l_1},z_{0,l_2},z_{0,l_3},z_{0,l_4})\bmss_{n,l_1,l_2}\bmss_{n,l_3,l_4}\right\}\\
	&-2\sum_{l_1,l_2,l_3,l_4}\left\{\cum(z_{0,l_1},z_{0,l_2},z_{0,l_3},z_{0,l_4})\bmss_{1,l_1,l_2}\bmss_{n,l_3,l_4}\right\}\\
	&+\sum_{l_1,l_2,l_3,l_4}\left\{\cum(z_{0,l_1},z_{0,l_2},z_{0,l_3},z_{0,l_4})\bmss_{1,l_1,l_2}\bmss_{1,l_3,l_4}\right\}\\
	=&o\{\tr(\bms_n^2)\}\tr\{(\bms_1-\bms_n)^2\}.
\end{align*}
The conclusion then follows. 
The proof for the case $\tau>\tau^{\ast}$ is similar and is therefore omitted. 
Hence, the lemma follows. \done

\subsection{Proof of Lemma~\ref{lemma_M_tau_V_tau}}

We first bound $M_n(\tau)$. 
Decompose $M_n(\tau)$ as 
\begin{align*}
	M_n(\tau)=&\frac{1}{P_{\tau}^{2}}\sum_{i\neq j}^{\tau}\Y_{i}^{\top}\Y_{j}+\frac{1}{P_{n-\tau}^{2}}\sum_{i,j=\tau+1,i\neq j}^{n}\Y_{i}^{\top}\Y_{j}-\frac{2}{\tau(n-\tau)}\sum_{i=1}^{\tau}\sum_{j=\tau+1}^{n}\Y_{i}^{\top}\Y_{j}\\
	&+\frac{2(n-\tau^{\ast}-1)\sum_{i=\tau^{\ast}+1}^{n}\Y_{i}^{\top}\bm\delta +2(n-\tau^{\ast})\sum_{i=\tau+1}^{\tau^{\ast}}\Y_{i}^{\top}\bm\delta}{(n-\tau)(n-\tau-1)}-\frac{2(n-\tau^{\ast})\sum_{i=1}^{\tau}\Y_{i}^{\top}\bm\delta}{\tau(n-\tau)}\\
	&+\frac{(n-\tau^{\ast})(n-\tau^{\ast}-1)}{(n-\tau)(n-\tau-1)}\|\bm\delta\|^2\\
	:=&R_1(\tau)+\E\{M_n(\tau)\},
\end{align*}
where 
\begin{equation*}
	\E\{M_n(\tau)\}=\frac{(n-\tau^{\ast})(n-\tau^{\ast}-1)}{(n-\tau)(n-\tau-1)}\|\bm\delta\|^2
\end{equation*}
and $R_1(\tau)=M_n(\tau)-\E\{M_n(\tau)\}$.

We will give a uniform bound for $R_1(\tau)$ by bounding each item in $R_1(\tau)$.
Let $R_{11}(\tau)=\sum_{i\neq j}^{\tau}\Y_{i}^{\top}\Y_j$, $R_{12}(\tau)=\sum_{i,j=\tau+1,i\neq j}^{n}\Y_{i}^{\top}\Y_{j}$,
$R_{13}(\tau)=\sum_{i=1}^{\tau}\sum_{j=\tau+1}^{n}\Y_{i}^{\top}\Y_{j}$,
$R_{14}(\tau)=\sum_{i=\tau^{\ast}+1}^{n}\Y_{i}^{\top}\bm\delta$,
$R_{15}(\tau)=\sum_{i=\tau+1}^{\tau^{\ast}}\Y_{i}^{\top}\bm\delta$, and
$R_{16}(\tau)=\sum_{i=1}^{\tau}\Y_{i}^{\top}\bm\delta$.

For $R_{11}(\tau)$, obviously, it is a martingale. By Lemma~\ref{lemma_maximum_inequality}, for  any $x>0$,
\begin{equation*}
	\Pr\left\{\max_{2\leq \tau \leq n-2} |R_{11}(\tau)|>x\right\}\leq \frac{4\var\{R_{11}(\tau)\}}{x^2}\leq \frac{8n^2\tr\{(\bms_1+\bms_n)^2\}}{x^2}.
\end{equation*}
Therefore, 
\begin{equation*}
	\max_{2\leq \tau \leq n-2} |R_{11}(\tau)|=O_p[n\tr^{1/2}\{(\bms_1+\bms_n)^2\}].
\end{equation*}

For $R_{12}(\tau)$, take $\Y_i=\Z_{n-i+1}$ and $H_{12}(\tau)= \sum_{i\neq j}^{\tau}\Z_{i}^{\top}\Z_{j} $.  
Clearly, $H_{12}(\tau)$ is a martingale and $R_{12}(\tau)=H_{12}(n-\tau)$.
Similarly to the proof for bounding $R_{11}(\tau)$, we have 
\begin{equation*}
	\max_{2\leq \tau \leq n-2} |R_{12}(\tau)|=O_p[n\tr^{1/2}\{(\bms_1+\bms_n)^2\}]. 
\end{equation*}

For $R_{13}$, observing that $R_{13}(\tau)=2^{-1}\{\sum_{i\neq j}^{n}\Y_{i}^{\top}\Y_{j}-R_{11}(\tau)-R_{12}(\tau)\}$, we obtain that 
\begin{equation*}
	\max_{2\leq \tau \leq n-2} |R_{13}(\tau)|=O_p[n\tr^{1/2}\{(\bms_1+\bms_n)^2\}]. 
\end{equation*}

For $R_{14}(\tau)$, $R_{14}(\tau)=\sum_{i=\tau^{\ast}+1}^{n}\Y_{i}^{\top}\bm\delta=O_p[\{n\bm\delta^{\top}(\bms_1+\bms_n) \bm\delta\}^{1/2}]$. 
Similar to the proofs for bounding $R_{11}(\tau)$ and $R_{12}(\tau)$, it's easy to show 
\begin{equation*}
	\max_{2\leq \tau \leq \tau^{\ast}-1} |R_{15}(\tau)|=O_p[\{n\bm\delta^{\top}(\bms_1+\bms_n) \bm\delta\}^{1/2}], \quad
	\max_{2\leq \tau \leq \tau^{\ast}} |R_{16}(\tau)|=O_p[\{n\bm\delta^{\top}(\bms_1+\bms_n) \bm\delta\}^{1/2}].
\end{equation*}
Since $\tau/n$ and $\tau^{\ast}/n$ are bounded away from 0 and 1, together with the results we prove above, we conclude that
\begin{equation*}
	\max_{\lambda_n\leq \tau \leq  \tau^{\ast}}\left|R_1(\tau)\right|=O_p[n^{-1}\tr^{1/2}\{(\bms_1+\bms_n)^2\}]+O_p[\{n^{-1}\bm\delta(\bms_1+\bms_n) \bm\delta\}^{1/2}].
\end{equation*}
Therefore, 
\begin{equation*}
	M_n(\tau)= \frac{(n-\tau^{\ast})(n-\tau^{\ast}-1)}{(n-\tau)(n-\tau-1)}\|\bm\delta\|^2 +O_p[n^{-1}\tr^{1/2}\{(\bms_1+\bms_n)^2\}]+O_p[\{n^{-1}\bm\delta(\bms_1+\bms_n) \bm\delta\}^{1/2}]
\end{equation*}
holds uniformly for $\lambda_n\leq \tau\leq \tau^{\ast}$.

Next, we bound $V_n(\tau)$. Under $H_1$, we can rewrite $V_n(\tau)$ as $V_n(\tau):=V_{n,1}(\tau)+V_{n,2}(\tau)+V_{n,3}(\tau)$, where
\begin{align*}
	V_{n,1}(\tau)=&\frac{\sum_{1\leq i_1,i_2,i_3,i_4\leq \tau}^{\ast}\left\{(\Y_{i_1}-\Y_{i_3})^{\top}(\Y_{i_2}-\Y_{i_4})\right\}^2}{4P_{\tau}^{4}} \\
	&+ \frac{\sum_{\tau+1\leq i_1,i_2,i_3,i_4\leq n}^{\ast}\left\{(\Y_{i_1}-\Y_{i_3})^{\top}(\Y_{i_2}-\Y_{i_4})\right\}^2}{4P_{n-\tau}^{4}}\\
	&-\frac{\sum_{i_1\neq i_3}^{\tau}\sum_{\tau+1\leq i_2,i_4\leq n}^{\ast} \left\{(\Y_{i_1}-\Y_{i_3})^{\top}(\Y_{i_2}-\Y_{i_4})\right\}^2 }{2\tau(\tau-1)(n-\tau)(n-\tau-1)},\\
	V_{n,2}(\tau)=&\frac{\sum_{1\leq i_1,i_2,i_3,i_4\leq \tau}^{\ast}\Delta_{1}(i_1,i_2,i_3,i_4)}{4P_{\tau}^{4}}+ \frac{\sum_{\tau+1\leq i_1,i_2,i_3,i_4\leq n}^{\ast}\Delta_{1}(i_1,i_2,i_3,i_4)}{4P_{n-\tau}^{4}}\\
	&-\frac{\sum_{i_1\neq i_3}^{\tau}\sum_{\tau+1\leq i_2,i_4\leq n}^{\ast} \Delta_{1}(i_1,i_2,i_3,i_4) }{2\tau(\tau-1)(n-\tau)(n-\tau-1)},\\
	V_{n,3}(\tau)=&\frac{\sum_{1\leq i_1,i_2,i_3,i_4\leq \tau}^{\ast}\Delta_{2}(i_1,i_2,i_3,i_4)}{4P_{\tau}^{4}}+ \frac{\sum_{\tau+1\leq i_1,i_2,i_3,i_4\leq n}^{\ast}\Delta_{2}(i_1,i_2,i_3,i_4)}{4P_{n-\tau}^{4}}\\
	&-\frac{\sum_{i_1\neq i_3}^{\tau}\sum_{\tau+1\leq i_2,i_4\leq n}^{\ast} \Delta_{2}(i_1,i_2,i_3,i_4) }{2\tau(\tau-1)(n-\tau)(n-\tau-1)},
\end{align*}
and
\begin{align*}
	\Delta_{1}(i_1,i_2,i_3,i_4)=&2(\Y_{i_1}-\Y_{i_3})^{\top}(\Y_{i_2}-\Y_{i_4})\E^{\top}(\X_{i_1}-\X_{i_3})(\Y_{i_2}-\Y_{i_4})\\
	&+2(\Y_{i_1}-\Y_{i_3})^{\top}(\Y_{i_2}-\Y_{i_4})\E^{\top}(\X_{i_2}-\X_{i_4})(\Y_{i_1}-\Y_{i_3})\\
	&+2(\Y_{i_1}-\Y_{i_3})^{\top}(\Y_{i_2}-\Y_{i_4})\E^{\top}(\X_{i_1}-\Y_{i_3})\E(\X_{i_2}-\X_{i_4}),
\end{align*} 
\begin{align*}
	\Delta_{2}(i_1,i_2,i_3,i_4)=&\big\{\E^{\top}(\X_{i_1}-\X_{i_3})(\Y_{i_2}-\Y_{i_4})+\E^{\top}(\X_{i_2}-\X_{i_4})(\Y_{i_1}-\Y_{i_3})\\
	&+\E^{\top}(\X_{i_1}-\X_{i_3})\E(\X_{i_2}-\X_{i_4})\big\}^2.
\end{align*}
We will bound $V_n(\tau)$ by bounding $V_{n,1}(\tau)$, $V_{n,2}(\tau)$, and $V_{n,3}(\tau)$, respectively.

To bound $V_{n,2}(\tau)$, we will bound its each term, respectively.
Note that for $\tau\leq \tau^{\ast}$,
\begin{equation*}
	\sum_{1\leq i_1,i_2,i_3,i_4\leq \tau}^{\ast}\Delta_{1}(i_1,i_2,i_3,i_4)=0,
\end{equation*}
and
\begin{align*}
	&\sum_{i_1\neq i_3}^{\tau}\sum_{\tau+1\leq i_2,i_4\leq n}^{\ast} \Delta_{1}(i_1,i_2,i_3,i_4)\\
	= &\sum_{i_1\neq i_3}^{\tau}\sum_{\tau+1\leq i_2,i_4\leq n}^{\ast} 2(\Y_{i_1}-\Y_{i_3})^{\top}(\Y_{i_2}-\Y_{i_4})\E^{\top}(\X_{i_2}-\X_{i_4})(\Y_{i_1}-\Y_{i_3}) \\
	=&-4\sum_{i_1\neq i_3}^{\tau}\sum_{i_2=\tau+1}^{\tau^{\ast}}\sum_{i_4=\tau^{\ast}+1}^{n}(\Y_{i_1}-\Y_{i_3})^{\top}(\Y_{i_2}-\Y_{i_4})\bm\delta^{\top}(\Y_{i_1}-\Y_{i_3})\\
	=&-4\sum_{i_1\neq i_3}^{\tau}\sum_{i_2=\tau+1}^{\tau^{\ast}}\sum_{i_4=\tau^{\ast}+1}^{n}\bigg( \Y_{i_1}^{\top}\Y_{i_2}\bm\delta^{\top}\Y_{i_1}
	-\Y_{i_1}^{\top}\Y_{i_2}\bm\delta^{\top}\Y_{i_3}
	-\Y_{i_1}^{\top}\Y_{i_4}\bm\delta^{\top}\Y_{i_1}
	+\Y_{i_1}^{\top}\Y_{i_4}\bm\delta^{\top}\Y_{i_3}\\
	&-\Y_{i_3}^{\top}\Y_{i_2}\bm\delta^{\top}\Y_{i_1}
	+\Y_{i_3}^{\top}\Y_{i_2}\bm\delta^{\top}\Y_{i_3}
	+\Y_{i_3}^{\top}\Y_{i_4}\bm\delta^{\top}\Y_{i_1}
	-\Y_{i_3}^{\top}\Y_{i_4}\bm\delta^{\top}\Y_{i_3}\bigg).
\end{align*}  
To bound $\sum_{i_1\neq i_3}^{\tau}\sum_{\tau+1\leq i_2,i_4\leq n}^{\ast} \Delta_{1}(i_1,i_2,i_3,i_4)$, we decompose the overall summation and consider each individual summation.
For the first term, decompose
\begin{align*}
	\sum_{i_1=1}^{\tau}\sum_{i_2=\tau+1}^{\tau^{\ast}} \Y_{i_1}^{\top}\Y_{i_2}\bm\delta^{\top}\Y_{i_1}=\sum_{i_1< i_2}^{\tau^{\ast}} \Y_{i_1}^{\top}\Y_{i_2}\bm\delta^{\top}\Y_{i_1}-\sum_{i_1< i_2}^{\tau} \Y_{i_1}^{\top}\Y_{i_2}\bm\delta^{\top}\Y_{i_1}-\sum_{i_1,i_2=\tau+1, i_1< i_2}^{\tau^{\ast}} \Y_{i_1}^{\top}\Y_{i_2}\bm\delta^{\top}\Y_{i_1}.
\end{align*}
Further, write
\begin{equation*}
	\sum_{i_1< i_2}^{\tau} \Y_{i_1}^{\top}\Y_{i_2}\bm\delta^{\top}\Y_{i_1}=\sum_{i_1< i_2}^{\tau} \Y_{i_2}^{\top}(\Y_{i_1}\Y_{i_1}^{\top}-\bms_1)\bm\delta+\tau\sum_{i_2=2}^{\tau}\Y_{i_2}^{\top}\bms_1\bm\delta.
\end{equation*}
By Lemma~\ref{lemma_maximum_inequality}, we have
\begin{align*}
	&\max_{2\leq \tau \leq \tau^{\ast}}\left|\sum_{i_1< i_2}^{\tau} \Y_{i_2}^{\top}(\Y_{i_1}\Y_{i_1}^{\top}-\bms_1)\bm\delta\right|\\
	=& n^{3/2} \bigg\{O_p\left( \left[ \tr\{(\bms_1+\bms_n)^2\}\bm\delta^{\top}(\bms_1+\bms_n)\bm\delta\right]^{1/2} \right)+O_p\left[ \{\bm\delta^{\top}(\bms_1+\bms_n)^3\bm\delta\}^{1/2}\right]\\
	&+o_p\left[\tr^{3/4}\{(\bms_{1}+\bms_n)^2\}\|\bm\delta\|\right]\bigg\}\\
	=&o_p\left[n^{3/2}\tr^{3/4}\{(\bms_{1}+\bms_n)^2\}\|\bm\delta\|\right],
\end{align*}
and 
\begin{equation*}
	\max_{2\leq \tau \leq \tau^{\ast}}\left| \sum_{i_2=2}^{\tau}\Y_{i_2}^{\top}\bms_1\bm\delta\right|=O_p\left\{n^{1/2}\bm\delta^{\top}(\bms_1+\bms_n)^3\bm\delta\right\}.
\end{equation*}
Hence, 
\begin{equation*}
	\max_{2\leq \tau \leq \tau^{\ast}}\left|\sum_{i_1< i_2}^{\tau} \Y_{i_2}^{\top}\Y_{i_1}\Y_{i_1}^{\top}\bm\delta\right|=o_p\left[n^{3/2}\tr^{3/4}\{(\bms_{1}+\bms_n)^2\}\|\bm\delta\|\right].
\end{equation*}
Similarly,
\begin{equation*}
	\max_{2\leq \tau\leq \tau^{\ast}-2}\left|\sum_{i_1,i_2=\tau+1, i_1< i_2}^{\tau^{\ast}} \Y_{i_1}^{\top}\Y_{i_2}\bm\delta^{\top}\Y_{i_1}\right|=o_p\left[n^{3/2}\tr^{3/4}\{(\bms_{1}+\bms_n)^2\}\|\bm\delta\|\right].
\end{equation*}
Thus
\begin{equation*}
	\max_{2\leq \tau\leq \tau^{\ast}-1} \left|\sum_{i_1=1}^{\tau}\sum_{i_2=\tau+1}^{\tau^{\ast}} \Y_{i_1}^{\top}\Y_{i_2}\bm\delta^{\top}\Y_{i_1}\right|= o_p\left[n^{3/2}\tr^{3/4}\{(\bms_{1}+\bms_n)^2\}\|\bm\delta\|\right].  
\end{equation*}
Similarly,
\begin{equation*}
	\max_{2\leq \tau\leq \tau^{\ast}} \left|\sum_{i_1=1}^{\tau}\sum_{i_4=\tau^{\ast}+1}^{n} \Y_{i_1}^{\top}\Y_{i_4}\bm\delta^{\top}\Y_{i_1}\right|= o_p\left[n^{3/2}\tr^{3/4}\{(\bms_{1}+\bms_n)^2\}\|\bm\delta\|\right],
\end{equation*}
\begin{equation*}
	\max_{2\leq \tau\leq \tau^{\ast}-1} \left|\sum_{i_3=1}^{\tau}\sum_{i_2=\tau+1}^{\tau^{\ast}} \Y_{i_3}^{\top}\Y_{i_2}\bm\delta^{\top}\Y_{i_3}\right|= o_p\left[n^{3/2}\tr^{3/4}\{(\bms_{1}+\bms_n)^2\}\|\bm\delta\|\right],
\end{equation*} 
and
\begin{equation*}
	\max_{2\leq \tau\leq \tau^{\ast}-1} \left|\sum_{i_3=1}^{\tau}\sum_{i_4=\tau^{\ast}+1}^{n} \Y_{i_3}^{\top}\Y_{i_4}\bm\delta^{\top}\Y_{i_3}\right|= o_p\left[n^{3/2}\tr^{3/4}\{(\bms_{1}+\bms_n)^2\}\|\bm\delta\|\right].  
\end{equation*}
For any $x>0$,
\begin{align*}
	&\Pr\left(\max_{2\leq \tau\leq \tau^{\ast}-1}  \left|\sum_{i_1\neq i_3}^{\ast}\sum_{i_2=\tau+1}^{\tau^{\ast}}\Y_{i_1}^{\top}\Y_{i_2}\bm\delta^{\top}\Y_{i_3}\right|\geq x \right)\\
	\leq &\sum_{\tau=2}^{\tau^{\ast}-1}\frac{\var\left( \sum_{i_1\neq i_3}^{\tau}\sum_{i_2=\tau+1}^{\tau^{\ast}}\Y_{i_1}^{\top}\Y_{i_2}\bm\delta^{\top}\Y_{i_3}  \right)}{x^2}\\
	\leq & \frac{n^4\tr\{(\bms_1+\bms_n)^2\}\bm\delta^{\top}(\bms_1+\bms_n)\bm\delta }{x^2}.
\end{align*}
Thus
\begin{equation*}
	\max_{2\leq \tau\leq \tau^{\ast}-1}  \left|\sum_{i_1\neq i_3}^{\tau}\sum_{i_2=\tau+1}^{\tau^{\ast}}\Y_{i_1}^{\top}\Y_{i_2}\bm\delta^{\top}\Y_{i_3}\right|=O_p\left( n^2\left[ \tr\{(\bms_1+\bms_n)^2\}\bm\delta^{\top}(\bms_1+\bms_n)\bm\delta \right]^{1/2}  \right).
\end{equation*}
Similarly,
\begin{equation*}
	\max_{2\leq \tau\leq \tau^{\ast}}  \left|\sum_{i_1\neq i_3}^{\tau}\sum_{i_4=\tau^{\ast}+1}^{n}\Y_{i_1}^{\top}\Y_{i_4}\bm\delta^{\top}\Y_{i_3}\right|=O_p\left( n^2\left[ \tr\{(\bms_1+\bms_n)^2\}\bm\delta^{\top}(\bms_1+\bms_n)\bm\delta \right]^{1/2}  \right),
\end{equation*}
\begin{equation*}
	\max_{2\leq \tau\leq \tau^{\ast}-1}  \left|\sum_{i_1\neq i_3}^{\tau}\sum_{i_2=\tau+1}^{\tau^{\ast}}\Y_{i_3}^{\top}\Y_{i_2}\bm\delta^{\top}\Y_{i_1}\right|=O_p\left( n^2\left[ \tr\{(\bms_1+\bms_n)^2\}\bm\delta^{\top}(\bms_1+\bms_n)\bm\delta \right]^{1/2}  \right),
\end{equation*}
and
\begin{equation*}
	\max_{2\leq \tau\leq \tau^{\ast}}  \left|\sum_{i_1\neq i_3}^{\tau}\sum_{i_4=\tau^{\ast}+1}^{n}\Y_{i_3}^{\top}\Y_{i_4}\bm\delta^{\top}\Y_{i_1}\right|=O_p\left( n^2\left[ \tr\{(\bms_1+\bms_n)^2\}\bm\delta^{\top}(\bms_1+\bms_n)\bm\delta \right]^{1/2}  \right).
\end{equation*}
Thus
\begin{equation*}
	\max_{2\leq \tau \leq \tau^{\ast}} \left|\sum_{i_1\neq i_3}^{\tau}\sum_{\tau+1\leq i_2,i_4\leq n}^{\ast} \Delta_{1}(i_1,i_2,i_3,i_4)\right|=o_p\left[n^{7/2}\tr^{3/4}\{(\bms_{1}+\bms_n)^2\}\|\bm\delta\|\right].
\end{equation*}
The proof for bounding $ \sum_{\tau+1\leq i_1,i_2,i_3,i_4\leq n}^{\ast}\Delta_{1}(i_1,i_2,i_3,i_4)$ is similar, here we omit it. 
We have that 
\begin{equation*}
	\max_{4\leq \tau \leq \tau^{\ast}} \left|\sum_{i_1\neq i_3}^{\tau}\sum_{\tau+1\leq i_2,i_4\leq n}^{\ast} \Delta_{1}(i_1,i_2,i_3,i_4)\right|=o_p\left[n^{7/2}\tr^{3/4}\{(\bms_{1}+\bms_n)^2\}\|\bm\delta\|\right].
\end{equation*}
Put the results for each term of $V_{n,2}(\tau)$ together, we have
\begin{equation*}
	\max_{\lambda_n\leq \tau \leq \tau^{\ast}}\left|V_{n,2}(\tau) \right|=o_p\left[n^{-1/2}\tr^{3/4}\{(\bms_{1}+\bms_n)^2\}\|\bm\delta\|\right].
\end{equation*}

For $V_{n,3}(\tau)$, observe that
\begin{align*}
	|V_{n,3}(\tau)|\leq & \frac{\sum_{1\leq i_1,i_2,i_3,i_4\leq \tau}^{\ast}\Delta_{2}(i_1,i_2,i_3,i_4)}{4P_{\tau}^{4}}+ \frac{\sum_{\tau+1\leq i_1,i_2,i_3,i_4\leq n}^{\ast}\Delta_{2}(i_1,i_2,i_3,i_4)}{4P_{n-\tau}^{4}}\\
	&+\frac{\sum_{i_1\neq i_3}^{\tau}\sum_{\tau+1\leq i_2,i_4\leq n}^{\ast} \Delta_{2}(i_1,i_2,i_3,i_4) }{2\tau(\tau-1)(n-\tau)(n-\tau-1)}.
\end{align*}
Further,
\begin{small}
	\begin{align*}
		&\sum_{1\leq i_1,i_2,i_3,i_4\leq \tau}^{\ast}\Delta_{2}(i_1,i_2,i_3,i_4)+\sum_{\tau+1\leq i_1,i_2,i_3,i_4\leq n}^{\ast}\Delta_{2}(i_1,i_2,i_3,i_4)+\sum_{i_1\neq i_3}^{\tau}\sum_{\tau+1\leq i_2,i_4\leq n}^{\ast} \Delta_{2}(i_1,i_2,i_3,i_4)\\
		\leq & \sum_{1\leq i_1,i_2,i_3,i_4\leq n}^{\ast}\Delta_{2}(i_1,i_2,i_3,i_4).
	\end{align*}
\end{small}
To bound $V_{n,3}(\tau)$,  it suffices to bound $\sum_{1\leq i_1,i_2,i_3,i_4\leq n}^{\ast}\Delta_{2}(i_1,i_2,i_3,i_4)$.
Since
\begin{align*}
	\Delta_2(i_1,i_2,i_3,i_4) \leq& 3\{\E^{\top}(\X_{i_1}-\X_{i_3})(\Y_{i_2}-\Y_{i_4})\}^2+3\{\E^{\top}(\X_{i_2}-\X_{i_4})(\Y_{i_1}-\Y_{i_3})\}^2\\
	&+3\{\E^{\top}(\X_{i_1}-\X_{i_3})\E(\X_{i_2}-\X_{i_4})\}^2,
\end{align*}
we have
\begin{equation*}
	\E\{\Delta_2(i_1,i_2,i_3,i_4)\} \leq 12\bm\delta^{\top}(\bms_1+\bms_n)\bm\delta+ 3\|\bm\delta\|^4.
\end{equation*}
Thus
\begin{align*}
	&\sum_{1\leq i_1,i_2,i_3,i_4\leq n}^{\ast}\Delta_{2}(i_1,i_2,i_3,i_4)\\
	=&O_p\{n^4\bm\delta^{\top}(\bms_1+\bms_n)\bm\delta\}+O_p(n^4\|\bm\delta\|^4)\\
	=& o_p\left[n^4\|\bm\delta\|^2\sqrt{\tr\{(\bms_1+\bms_n)^2\}} \right] +O_p(n^4\|\bm\delta\|^4) ,
\end{align*}
and 
\begin{equation*}
	\max_{\lambda_n\leq \tau \leq \tau^{\ast}}\left|V_{n,3}(\tau)\right|=o_p\left[\|\bm\delta\|^2\sqrt{\tr\{(\bms_1+\bms_n)^2\}} \right]+O_p(\|\bm\delta\|^4).
\end{equation*}

It remains to bound $V_{n,1}(\tau)$. Note that
\begin{equation*}
	V_{n,1}(\tau)= V_{n,1,1}(\tau)+ V_{n,1,2}(\tau)+ V_{n,1,3}(\tau),
\end{equation*}
where
\begin{align*}
	V_{n,1,1}(\tau)=&\frac{1}{P_{\tau}^{2}}\sum_{i\neq j}^{\tau}(\Y_{i}^{\top}\Y_{j})^2-\frac{2}{P_{\tau}^{3}}\sum_{1\leq i,j,k\leq \tau}^{\ast}\Y_{i}^{\top}\Y_{j}\Y_{j}^{\top}\Y_k+\frac{1}{P_{\tau}^4}\sum_{1\leq i,j,k,l\leq \tau}^{\ast}\Y_{i}^{\top}\Y_{j}\Y_{k}^{\top}\Y_{l};\\
	V_{n,1,2}(\tau)=&\frac{1}{P_{n-\tau}^{2}}\sum_{\tau+1\leq i,j\leq n}^{\ast}(\Y_{i}^{\top}\Y_{j})^2-\frac{2}{P_{n-\tau}^{3}}\sum_{\tau+1\leq i,j,k\leq n}^{\ast}\Y_{i}^{\top}\Y_{j}\Y_{j}^{\top}\Y_k\\
	&+\frac{1}{P_{n-\tau}^4}\sum_{\tau+1\leq i,j,k,l\leq n}^{\ast}\Y_{i}^{\top}\Y_{j}\Y_{k}^{\top}\Y_{l};\\
	V_{n,1,3}(\tau)=&\frac{1}{\tau(n-\tau)}\sum_{j=1}^{\tau}\sum_{j=\tau+1}^{n}(\Y_{i}^{\top}\Y_{j})^2-\frac{1}{\tau(n-\tau)(\tau-1)}\sum_{1\leq i,k\leq \tau}^{\ast}\sum_{j=\tau+1}^{n}\Y_{i}^{\top}\Y_{j}\Y_{j}^{\top}\Y_k\\
	&-\frac{1}{\tau(n-\tau)(n-\tau-1)}\sum_{\tau+1\leq i,k\leq n}^{\ast}\sum_{j=1}^{\tau}\Y_{i}^{\top}\Y_{j}\Y_{j}^{\top}\Y_k\\
	&+\frac{1}{\tau(\tau-1)(n-\tau)(n-\tau-1)}\sum_{1\leq i,k\leq \tau}^{\ast}\sum_{\tau+1\leq j,l\leq n}^{\ast}\Y_{i}^{\top}\Y_{j}\Y_{k}^{\top}\Y_l.
\end{align*}
Rearrange $V_{n,1}(\tau)$ as
\begin{align*}
	V_{n,1}(\tau)&=\frac{1}{P_{\tau}^{2}}\sum_{i\neq j}^{\tau}(\Y_{i}^{\top}\Y_{j})^2+\frac{1}{P_{n-\tau}^{2}}\sum_{\tau+1\leq i,j\leq n}^{\ast}(\Y_{i}^{\top}\Y_{j})^2-\frac{2}{\tau(n-\tau)}\sum_{j=1}^{\tau}\sum_{j=\tau+1}^{n}(\Y_{i}^{\top}\Y_{j})^2+V_{n,1,4}(\tau)\\
	&:=\tilde{V}_{n,1}(\tau)+V_{n,1,4}(\tau).
\end{align*}
We first bound $V_{n,1,4}(\tau)$ by bounding its each term.
Write
\begin{align*}
	\sum_{1\leq i,j,k\leq \tau}^{\ast}\Y_{i}^{\top}\Y_{j}\Y_{j}^{\top}\Y_k=\sum_{1\leq i,j,k\leq \tau}^{\ast} \Y_{i}^{\top}(\Y_{j}\Y_{j}^{\top}-\bms_1)\Y_k+(\tau-3)\sum_{1\leq i,k\leq \tau}^{\ast}\Y_{i}^{\top}\bms_1\Y_k.
\end{align*}
Since
\begin{align*}
	&\Pr\left\{ \max_{3\leq \tau \leq \tau^{\ast}} \left|\sum_{1\leq i,j,k\leq \tau}^{\ast} \Y_{i}^{\top}(\Y_{j}\Y_{j}^{\top}-\bms_1)\Y_k\right|>x \right\}\\
	\leq & \frac{\sum_{3\leq \tau \leq \tau^{\ast}} \var\left\{ \sum_{1\leq i,j,k\leq \tau}^{\ast} \Y_{i}^{\top}(\Y_{j}\Y_{j}^{\top}-\bms_1)\Y_k \right\}  }{x^2}= \frac{O[n^4\tr^2\{(\bms_1+\bms_n)^2\}]}{x^2},
\end{align*}
we have
\begin{equation*}
	\max_{3\leq \tau \leq \tau^{\ast}} \left|\sum_{1\leq i,j,k\leq \tau}^{\ast} \Y_{i}^{\top}(\Y_{j}\Y_{j}^{\top}-\bms_1)\Y_k\right|=O_p\left[n^{2}\tr\{(\bms_1+\bms_n)^2\}\right].
\end{equation*}
By Lemma~\ref{lemma_maximum_inequality},
\begin{equation*}
	\max_{2\leq \tau \leq \tau^{\ast}}\left|\sum_{1\leq i,k\leq \tau}^{\ast}\Y_{i}^{\top}\bms_1\Y_k\right|=O_p\left[n\sqrt{ \tr\{(\bms_1+\bms_n)^4\}}\right].
\end{equation*}
Hence
\begin{align*}
	&\max_{3\leq \tau \leq \tau^{\ast}} \left|\sum_{1\leq i,j,k\leq \tau}^{\ast} \Y_{i}^{\top}\Y_{j}\Y_{j}^{\top}\Y_k\right|\\
	=&O_p[n^{2}\tr\{(\bms_1+\bms_n)^2\}]+O_p\left[n^2\sqrt{ \tr\{(\bms_1+\bms_n)^4\}}\right]\\
	=&O_p[n^{2}\tr\{(\bms_1+\bms_n)^2\}].
\end{align*}
Similarly,
\begin{equation*}
	\max_{3\leq \tau \leq \tau^{\ast}} \left|\sum_{\tau+1\leq i,j,k\leq n}^{\ast}\Y_{i}^{\top}\Y_{j}\Y_{j}^{\top}\Y_k\right|=O_p[n^{2}\tr\{(\bms_1+\bms_n)^2\}],
\end{equation*}
\begin{equation*}
	\max_{3\leq \tau \leq \tau^{\ast}} \left|\sum_{\tau+1\leq i,k\leq n}^{\ast}\sum_{j=1}^{\tau}\Y_{i}^{\top}\Y_{j}\Y_{j}^{\top}\Y_k\right|=O_p[n^{2}\tr\{(\bms_1+\bms_n)^2\}],
\end{equation*}
and 
\begin{equation*}
	\max_{3\leq \tau \leq \tau^{\ast}} \left|\sum_{1\leq i,k\leq \tau}^{\ast}\sum_{j=\tau+1}^{n}\Y_{i}^{\top}\Y_{j}\Y_{j}^{\top}\Y_k\right|=O_p[n^{2}\tr\{(\bms_1+\bms_n)^2\}].
\end{equation*}
Besides, we can also show that
\begin{equation*}
	\max_{3\leq \tau \leq \tau^{\ast}} \left|\sum_{1\leq i,j,k,l\leq \tau}^{\ast}\Y_{i}^{\top}\Y_{j}\Y_{k}^{\top}\Y_{l}\right|=O_p[n^{5/2}\tr\{(\bms_1+\bms_n)^2\}],
\end{equation*}
\begin{equation*}
	\max_{3\leq \tau \leq \tau^{\ast}} \left|\sum_{\tau+1\leq i,j,k,l\leq n}^{\ast}\Y_{i}^{\top}\Y_{j}\Y_{k}^{\top}\Y_{l}\right|=O_p[n^{5/2}\tr\{(\bms_1+\bms_n)^2\}],
\end{equation*}
and
\begin{equation*}
	\max_{3\leq \tau \leq \tau^{\ast}} \left|\sum_{1\leq i,k\leq \tau}^{\ast}\sum_{\tau+1\leq j,l\leq n}^{\ast}\Y_{i}^{\top}\Y_{j}\Y_{k}^{\top}\Y_l\right|=O_p[n^{5/2}\tr\{(\bms_1+\bms_n)^2\}].
\end{equation*}
Based on the results,
\begin{equation*}
	\max_{ \lambda_n \leq \tau \leq \tau^{\ast}} \left|V_{n,1,4}(\tau)\right|=O_p[n^{-1}\tr\{(\bms_1+\bms_n)^2\}].
\end{equation*}
Next, we focus on $\tilde{V}_{n,1}(\tau)$.
Based on the decomposition
\begin{align*}
	&(\Y_{i}^{\top}\Y_{j})^2 -\E\{(\Y_{i}^{\top}\Y_{j})^2 \}\\
	=&\tr\{(\Y_i\Y_{i}^{\top}-\bms_i)(\Y_j\Y_{j}^{\top}-\bms_j)\}+\Y_{j}^{\top}\bms_i\Y_j+\Y_{i}^{\top}\bms_j\Y_i-2\tr(\bms_i\bms_j),
\end{align*}
we can write
\begin{align*}
	\tilde{V}_{n,1}(\tau)=&\frac{\sum_{i\neq j}^{\tau}\tr\{(\Y_i\Y_{i}^{\top}-\bms_i)(\Y_j\Y_{j}^{\top}-\bms_j)\}}{\tau(\tau-1) }+\frac{\sum_{\tau+1\leq i,j\leq n}^{\ast}\tr\{(\Y_i\Y_{i}^{\top}-\bms_i)(\Y_j\Y_{j}^{\top}-\bms_j)\}}{(n-\tau)(n-\tau-1)}\\
	&-\frac{2\sum_{i=1}^{\tau}\sum_{j=\tau+1}^{n}\tr\{(\Y_i\Y_{i}^{\top}-\bms_i)(\Y_j\Y_{j}^{\top}-\bms_j)\}}{\tau(n-\tau)}+\E\{\tilde{V}_{n,1}(\tau)\}\\
	=&\frac{2(n-\tau^{\ast})\sum_{i=1}^{\tau}[  \Y_{i}^{\top}(\bms_1-\bms_n)\Y_i-\tr\{\bms_1(\bms_1-\bms_n)\} ]}{\tau(n-\tau)}\\
	&-\frac{2(n-\tau^{\ast})\sum_{i=\tau+1}^{\tau^{\ast}}[\Y_{i}^{\top}(\bms_1-\bms_n)\Y_i-\tr\{\bms_1(\bms_1-\bms_n)\} ]}{(n-\tau)(n-\tau-1)}\\
	&-\frac{2(n-\tau^{\ast}-1)\sum_{i=\tau^{\ast}+1}^{n}[  \Y_{i}^{\top}(\bms_1-\bms_n)\Y_i-\tr\{\bms_n(\bms_1-\bms_n)\} ]}{(n-\tau)(n-\tau-1)}\\
	&+\E\{\tilde{V}_{n,1}(\tau)\}.
\end{align*}
By Lemma~\ref{lemma_maximum_inequality}, we have
\begin{equation*}
	\max_{2\leq \tau \leq \tau^{\ast}}\left|\sum_{i\neq j}^{\tau}\tr\{(\Y_i\Y_{i}^{\top}-\bms_i)(\Y_j\Y_{j}^{\top}-\bms_j)\}\right|=O_p[n\tr\{(\bms_1+\bms_n)^2\}],
\end{equation*}
\begin{equation*}
	\max_{1\leq \tau \leq \tau^{\ast}}\left|\sum_{\tau+1\leq i, j\leq n }^{\ast}\tr\{(\Y_i\Y_{i}^{\top}-\bms_i)(\Y_j\Y_{j}^{\top}-\bms_j)\}\right|=O_p[n\tr\{(\bms_1+\bms_n)^2\}],
\end{equation*}
Thus
\begin{equation*}
	\max_{1\leq \tau \leq \tau^{\ast}}\left|\sum_{i=1}^{\tau}\sum_{j=\tau+1}^{n}\tr\{(\Y_i\Y_{i}^{\top}-\bms_i)(\Y_j\Y_{j}^{\top}-\bms_j)\}\right|=O_p[n\tr\{(\bms_1+\bms_n)^2\}].
\end{equation*}
Besides,
\begin{align*}
	&\max_{1\leq \tau \leq \tau^{\ast}}\left|\sum_{i=1}^{\tau}[  \Y_{i}^{\top}(\bms_1-\bms_n)\Y_i-\tr\{\bms_1(\bms_1-\bms_n)\} ]\right|\\
	=&O_p\left([ n\tr\{(\bms_1+\bms_n)^2\}\tr\{ (\bms_1-\bms_n)^2 \}]^{1/2}\right),
\end{align*}
\begin{align*}
	&\max_{1\leq \tau \leq \tau^{\ast}-1}\left|\sum_{i=\tau+1}^{\tau^{\ast}}\{  \Y_{i}^{\top}(\bms_1-\bms_n)\Y_i-\tr\{\bms_1(\bms_1-\bms_n)\} \}\right|\\
	=&O_p\left([ n\tr\{(\bms_1+\bms_n)^2\}\tr\{ (\bms_1-\bms_n)^2 \}]^{1/2}\right).
\end{align*}
Together with 
\begin{align*}
	&\sum_{i=\tau^{\ast}+1}^{n}\{  \Y_{i}^{\top}(\bms_1-\bms_n)\Y_i-\tr\{\bms_n(\bms_1-\bms_n)\} \}\\
	=&O_p\left([ n\tr\{(\bms_1+\bms_n)^2\}\tr\{ (\bms_1-\bms_n)^2 \}]^{1/2}\right),
\end{align*}
we have
\begin{align*}
	&\max_{\lambda_n \leq \tau \leq \tau^{\ast}}\left|\tilde{V}_{n,1}(\tau)-\E\{\tilde{V}_{n,1}(\tau)\right|\\
	=&O_p[n^{-1}\tr\{(\bms_1+\bms_n)^2\}]+O_p\left([ n^{-1}\tr\{(\bms_1+\bms_n)^2\}\tr\{ (\bms_1-\bms_n)^2 \}]^{1/2}\right).
\end{align*}
Note that
\begin{equation*}
	\E\{\tilde{V}_{n,1}(\tau)\}=\frac{(n-\tau^{\ast})(n-\tau^{\ast}-1)}{(n-\tau)(n-\tau-1)}\tr\{(\bms_1-\bms_n)^2\},
\end{equation*}
together with the uniform bounds we obtain for $V_{n,2}(\tau)$, $V_{n,3}(\tau)$, $V_{n,1,4}(\tau)$, and $\tilde{V}_{n,1}(\tau)-\E\{\tilde{V}_{n,1}(\tau)\}$, we conclude that
\begin{align*}
	V_n(\tau)=&\frac{(n-\tau^{\ast})(n-\tau^{\ast}-1)}{(n-\tau)(n-\tau-1)}\tr\{(\bms_1-\bms_n)^2\}\\
	&+O_p[n^{-1}\tr\{(\bms_1+\bms_n)^2\}]+O_p\left([ n^{-1}\tr\{(\bms_1+\bms_n)^2\}\tr\{ (\bms_1-\bms_n)^2 \}]^{1/2}\right) \nonumber \\ 
	&+o_p\left[n^{-1/2}\tr^{3/4}\{(\bms_{1}+\bms_n)^2\}\|\bm\delta\|\right]+o_p\left[\|\bm\delta\|^2\sqrt{\tr\{(\bms_1+\bms_n)^2\}} \right]+O_p(\|\bm\delta\|^4)
\end{align*}
holds uniformly for all $\tau\in [\lambda_n,\tau^{\ast}]$. 
The lemma follows. \done

\subsection{Proof of Lemma~\ref{lemma_normal}}

We first establish the asymptotic normality of $M_n$.
Note that
\begin{equation*}
	M_n(\tau)=\frac{\sum_{i_1\neq i_2}^{\tau}\sum_{j_1,j_2=\tau+1,j_1\neq j_2}^{n}(\X_{i_1}-\X_{j_1})^{\top}(\X_{i_2}-\X_{j_2})}{\tau(\tau-1)(n-\tau)(n-\tau-1)},
\end{equation*}
and 
\begin{equation*}
	M_n=\sum_{\tau=2}^{n-2}\left\{ \frac{\tau(n-\tau)}{n} M_n(\tau)\right\}.
\end{equation*}
Under $\bH_0$, without loss of generality, we can assume $\bm\mu_i=\bm 0$. 
Rewrite $M_n$ as
\begin{equation*}
	M_n=\sum_{k=2}^{n}\sum_{i=1}^{k-1}a_{i,k}\X_{i}^{\top}\X_{k},
\end{equation*}
where $a_{i,k}=\sum_{\tau=2}^{i-1}2(n-\tau-1)^{-1}(1-n^{-1})+\sum_{\tau=k}^{n-2}2(\tau-1)^{-1}(1-n^{-1})+6n^{-1}-2$. (For convenience of notation, define $\sum_{i=a}^{b}c_i=0$ if $a>b$.)
We have 
\begin{equation*}
	\var(M_{n}) = \sum_{k=2}^{n}\sum_{i=1}^{k-1} a_{i,k}^{2}\var(\X_{i}^{\top}\X_{k})=\frac{2\pi^2-18}{3} n^2\tr(\bms_1^2)\{1+o(1)\},
\end{equation*}
by the facts $\var(\X_{i}^{\top}\X_{k})=\tr(\bms_1^2)$ and 
\begin{equation*}
	\sum_{k=2}^{n}\sum_{i=1}^{k-1} a_{i,k}^{2}= \frac{2\pi^2-18}{3} n^2\{1+o(1)\}.
\end{equation*}
Define
\begin{equation*}
	G_{n,k}=\X_{k}^{\top}\sum_{i=1}^{k-1}a_{i,k}\X_i.
\end{equation*}
We have $M_{n}=\sum_{k=1}^{n}G_{n,k}$.
Let $\mathcal{F}_0=\{\varnothing ,\Omega\}$, $\mathcal{F}_k=\sigma\{\X_1,\ldots \X_k\}$, for $k=1,\ldots,n$.
For each $n$, $\{G_{n,k}, 1\leq k\leq n\}$ is a martingale difference sequence with respect to the $\sigma$-fields $\{\mathcal{F}_k, 1\leq k\leq n\}$. To establish the asymptotic normality of $M_n$, we only need to show that as $n,p\rightarrow \infty$,
\begin{equation}\label{equation_proof_mean_CLT1}
	\frac{\sum_{k=1}^{n}\E( G_{n,k}^4)}{\var^2(M_{n})}\rightarrow 0,
\end{equation}
and
\begin{equation}\label{equation_proof_mean_CLT2}
	\frac{\sum_{k=1}^{n}\E( G_{n,k}^2 \mid \mathcal{F}_{k-1})}{\var(M_{n})}\rightarrow 1
\end{equation}
in probability.

For $1\leq i<k \leq n$, define $\xi_{i,k}=\X_{i}^{\top}\X_{k}$.
Then, $G_{n,k}=\sum_{i=1}^{k-1}a_{i,k}\xi_{i,k}$, and 
\begin{align*}
	\E(G_{n,k}^4)=\left(\sum_{i=1}^{k-1} a_{i,k}^4\E\xi_{i,k}^4 + \sum_{i_1\neq i_2}^{k-1}3a_{i_1,k}^2a_{i_2,k}^2\E\xi_{i_1,k}^2\xi_{i_2,k}^2\right).
\end{align*}
By Lemma~\ref{lemma_basic_results} (1),
\begin{equation*}
	\E(\xi_{i,k}^4)=O\{\tr^2(\bms_1^2)\}, \quad \E(\xi_{i_1,k}^2\xi_{i_2,k}^2)=O\{\tr^2(\bms_1^2)\},
\end{equation*}
we have
\begin{equation*}
	\sum_{k=1}^{n}\E( G_{n,k}^4)=O\{n^3\tr^2(\bms_1^2)\}=o\{\var^2(M_{n})\},
\end{equation*}
where we use the fact that $\sum_{k=2}^{n}\sum_{i_1,i_2}^{k-1}a_{i_1,k}^2a_{i_2,k}^2=O(n^3)$. Thus (\ref{equation_proof_mean_CLT1}) holds.

For (\ref{equation_proof_mean_CLT2}), it is sufficient to verify that as $n,p\rightarrow\infty$,
\begin{equation*}
	\var\left\{  \sum_{k=1}^{n}\E( G_{n,k}^2 \mid \mathcal{F}_{k-1})  \right\}/\var^2(M_{n})\rightarrow 0.
\end{equation*}
Note that
\begin{equation*}
	\sum_{k=1}^{n}\E( G_{n,k}^2 \mid \mathcal{F}_{k-1})=\sum_{k=1}^{n}\sum_{i=1}^{k-1}a_{i,k}^{2}\X_{i}^{\top}\bms_1\X_{i}+\sum_{k=1}^{n}\sum_{i_1\neq i_2}^{k-1}a_{i_1,k}a_{i_2,k}\X_{i_1}^{\top}\bms_1\X_{i_2}.
\end{equation*}
By Lemma~\ref{lemma_basic_results} (1), $\var(\X_{i}^{\top}\bms_1\X_{i})=O\{\tr^2(\bms_1^2)\}$ and  $\var(\X_{i_1}^{\top}\bms_1\X_{i_2})=O\{\tr(\bms_1^4)\}$, thus
\begin{align*}
	&\var\left\{  \sum_{k=1}^{n}\E( G_{n,k}^2 \mid \mathcal{F}_{k-1})  \right\}\\
	=&\sum_{i=1}^{n-1}\left(\sum_{k=i+1}^{n}a_{i,k}^2\right)^2\var(\X_{i}^{\top}\bms_1\X_{i})+4\sum_{i_1<i_2}\left(\sum_{k=i_2+1}^{n}a_{i_1,k}a_{i_2,k}\right)^2\var(\X_{i_1}^{\top}\bms_1\X_{i_2})\\
	=&O\{n^3\tr^2(\bms_1^2)\}+O\{n^4\tr(\bms_1^4)\}=o\{\var^2(M_{n})\},
\end{align*}
where we use the facts that $\sum_{i=1}^{n-1}\left(\sum_{k=i+1}^{n}a_{i,k}^2\right)^2=O(n^3)$ and 
$\sum_{i_1<i_2}\left(\sum_{k=i_2+1}^{n}a_{i_1,k}a_{i_2,k}\right)^2=O(n^4)$. Hence (\ref{equation_proof_mean_CLT2}) holds.

By the facts (\ref{equation_proof_mean_CLT1}) and (\ref{equation_proof_mean_CLT2}), we conclude that
$M_{n}/\var^{1/2}(M_{n})\rightarrow N(0,1)$ in distribution as $n,p\rightarrow \infty$. 
~\\

For $V_n$, under $\bH_0$, without loss of generality, we can assume $\bm\mu_i=\bm 0$. 
Recall that
\begin{equation*}
	V_n=\sum_{\tau=4}^{n-4}\left\{\frac{\tau(n-\tau)}{n}V_n(\tau)\right\}.
\end{equation*}
Define
\begin{equation*}
	V_{n,1}(\tau)=\frac{1}{P_{\tau}^{2}}\sum_{i\neq j}^{\tau}(\X_{i}^{\top}\X_{j})^2+\frac{1}{P_{n-\tau}^{2}}\sum_{i,j=\tau+1,i\neq j}^{n}(\X_{i}^{\top}\X_{j})^2-\frac{2}{\tau(n-\tau)}\sum_{i=1}^{\tau}\sum_{j=\tau+1}^{n}(\X_{i}^{\top}\X_{j})^2,
\end{equation*}
for $\tau=2,\ldots,n-2$.
Further, define
\begin{equation*}
	V_{n,1}=\sum_{\tau=2}^{n-2}{\frac{\tau(n-\tau)}{n}V_{n,1}(\tau)}.
\end{equation*}
We first show the asymptotic normality of $V_{n,1}$.
Rewrite $V_{n,1}$ as 
\begin{align*}
	V_{n,1}=\sum_{k=2}^{n}\sum_{i=1}^{k-1}a_{i,k}(\X_{i}^{\top}\X_{k})^2,
\end{align*}
where $a_{i,k}=\sum_{\tau=2}^{i-1}2(n-\tau-1)^{-1}(1-n^{-1})+\sum_{\tau=k}^{n-2}2(\tau-1)^{-1}(1-n^{-1})+6n^{-1}-2$.
Let $a_{i,k}=a_{k,i}$ for $i>k$, and $a_{i,i}=0$.
We have
\begin{align*}
	\var(V_{n,1})=&\sum_{k=2}^{n}\sum_{i=1}^{k-1}a_{i,k}^2\var\left\{(\X_{1}^{\top}\X_{2})^2\right\}\\
	&+\sum_{k=2}^{n}\sum_{i=1}^{k-1} a_{i,k}\left\{\sum_{j=1}^{n}a_{i,j}+\sum_{j=1}^{n}a_{j,k}-2a_{i,k}\right\}   \cov\left\{(\X_{1}^{\top}\X_{2})^2,(\X_{1}^{\top}\X_{3})^2\right\}.
\end{align*}
Using the facts $\sum_{j=1}^{n}a_{i,j}=\sum_{j=1}^{n}a_{j,k}=0$,
\begin{equation*}
	\sum_{k=2}^{n}\sum_{i=1}^{k-1} a_{i,k}^{2}= \frac{2\pi^2-18}{3} n^2\{1+o(1)\},
\end{equation*}
and the results of Lemma~\ref{lemma_basic_results} (2)
\begin{equation*}
	\cov\left\{(\X_{1}^{\top}\X_{2})^2,(\X_{1}^{\top}\X_{3})^2\right\}=o\{\tr^2(\bms_1^2)\},
\end{equation*}
and
\begin{equation*}
	\var\left\{(\X_{1}^{\top}\X_{2})^2\right\}=2\tr^2(\bms_1^2)\{1+o(1)\},
\end{equation*}
we conclude
\begin{equation*}
	\var(V_{n,1})= \frac{4\pi^2-36}{3} n^2\tr^2(\bms_1^2)\{1+o(1)\}.
\end{equation*}
Define
\begin{equation*}
	D_{n,k}=\tr\left\{(\X_k\X_{k}^{\top}-\bms_1)\sum_{i=1}^{k-1}a_{i,k}(\X_i\X_{i}^{\top}-\bms_1) \right\}.
\end{equation*}
We have $V_{n,1}=\sum_{k=1}^{n}D_{n,k}$.
Let $\mathcal{F}_0=\{\varnothing ,\Omega\}$, $\mathcal{F}_k=\sigma\{\X_1,\ldots \X_k\}$, for $k=1,\ldots,n$.
For each $n$, $\{D_{n,k}, 1\leq k\leq n\}$ is a martingale difference sequence with respect to the $\sigma$-fields $\{\mathcal{F}_k, 1\leq k\leq n\}$. To establish the asymptotic normality of $V_{n,1}$, we only need to show
\begin{equation}\label{equation_proof_var_CLT1}
	\frac{\sum_{k=1}^{n}\E( D_{n,k}^4)}{\var^2(V_{n,1})}\rightarrow 0,
\end{equation}
and
\begin{equation}\label{equation_proof_var_CLT2}
	\frac{\sum_{k=1}^{n}\E( D_{n,k}^2 \mid \mathcal{F}_{k-1})}{\var(V_{n,1})}\cp 1,
\end{equation}
as $n,p\rightarrow \infty$.

For (\ref{equation_proof_var_CLT1}), let $\xi_{i,k}=\tr\left 
\{(\X_{k}\X_{k}^{\top}-\bms_1)(\X_{i}^{\top}\X_{i}^{\top}-\bms_1)\right\}$.
Then $D_{n,k}=\sum_{i=1}^{k-1}a_{i,k}\xi_{i,k}$,
\begin{align*}
	\E(D_{n,k}^4)=\left(\sum_{i=1}^{k-1} a_{i,k}^4\E\xi_{i,k}^4 + \sum_{i_1\neq i_2}^{k-1}3a_{i_1,k}^2a_{i_2,k}^2\E\xi_{i_1,k}^2\xi_{i_2,k}^2\right).
\end{align*}
Using the results of Lemma~\ref{lemma_basic_results} (2), $\E\xi_{k,i}^4= O\{\tr^4(\bms_1^2)\}$, $\E(\xi_{k,i_1}^2\xi_{k,i_2}^2)\leq \E\xi_{k,i_1}^4 =O\{\tr^4(\bms_1^2)\}$, and the fact $\sum_{k=2}^{n}\sum_{i_1,i_2}^{k-1}a_{i_1,k}^2a_{i_2,k}^2=O(n^3)$, we have
\begin{equation*}
	\sum_{k=1}^{n}\E( D_{n,k}^4)=O\{n^3\tr^4(\bms_1^2)\}=o\{\var^2(V_{n,1})\}.
\end{equation*}

To prove (\ref{equation_proof_var_CLT2}), it is sufficient to verify that as $n,p\rightarrow \infty$,
\begin{equation*}
	\var\left\{  \sum_{k=1}^{n}\E( D_{n,k}^2 \mid \mathcal{F}_{k-1})  \right\}/\var^2(V_{n,1})\rightarrow 0.
\end{equation*}
For $1\leq i<k\leq n$, let $\xi_{i,k}=\tr\left\{ 
(\X_k\X_k^{\top}-\bms_1)(\X_i\X_i^{\top}-\bms_1) \right\}$ and $\zeta_i=\E(\xi_{i,k}^2\mid \X_i)$.
Let $\phi_{i_1,i_2}=\E(\xi_{i_1,k}\xi_{i_2,k}\mid \X_{i_1},\X_{i_2})$ for $i_1\neq i_2$.
Then,
\begin{equation*}
	\sum_{k=1}^{n}\E( D_{n,k}^2 \mid \mathcal{F}_{k-1})=\sum_{k=1}^{n}\sum_{i=1}^{k-1}a_{i,k}^{2}\zeta_i+\sum_{k=1}^{n}\sum_{i_1\neq i_2}^{k-1}a_{i_1,k}a_{i_2,k}\phi_{i_1,i_2},
\end{equation*}
and
\begin{align*}
	&\var\left\{  \sum_{k=1}^{n}\E( D_{n,k}^2 \mid \mathcal{F}_{k-1})  \right\}\\
	=&\sum_{i=1}^{n-1}\left(\sum_{k=i+1}^{n}a_{i,k}^2\right)^2\var(\zeta_i)+4\sum_{i_1<i_2}\left(\sum_{k=i_2+1}^{n}a_{i_1,k}a_{i_2,k}\right)^2\var(\phi_{i_1,i_2})\\
	=&O\{n^3\tr^4(\bms_1^2)\}+o\{n^4\tr(\bms_1^4)\}=o\{\var^2(V_{n,1})\},
\end{align*}
where we use the facts that $\sum_{i=1}^{n-1}\left(\sum_{k=i+1}^{n}a_{i,k}^2\right)^2=O(n^3)$,
$\sum_{i_1<i_2}\left(\sum_{k=i_2+1}^{n}a_{i_1,k}a_{i_2,k}\right)^2=O(n^4)$, and the results of Lemma~\ref{lemma_basic_results} (2) that $\E\zeta_i^2=O\{\tr^4(\bms_1^2)\}$, $\E(\phi_{i_1,i_2}^2)=o\{\tr^4(\bms_1^2)\}$ for $i_1\neq i_2$ and $\E(\phi_{i_1,i_2}\phi_{i_3,i_4})=0$ if $i_1 \notin \{i_3,i_4\}$ or $i_2 \notin \{i_3,i_4\}$. 

By the facts (\ref{equation_proof_var_CLT1})  and (\ref{equation_proof_var_CLT2}), we have
$V_{n,1}/\var^{1/2}(V_{n,1})\rightarrow N(0,1)$ in distribution as $n,p\rightarrow \infty$. 
Note that 
\begin{align*}
	V_n-V_{n,1}=&-\sum_{\tau=2}^{3}\frac{\tau(n-\tau)}{n}V_{n,1}(\tau)-\sum_{\tau=n-3}^{n-2}\frac{\tau(n-\tau)}{n}V_{n,1}(\tau)\\
	&+\sum_{\tau=4}^{n-4}\frac{\tau(n-\tau)}{n}\{V_n(\tau)-V_{n,1}(\tau)\}.
\end{align*}
It is easy to check that
\begin{equation*}
	\var\left\{\sum_{\tau=2}^{3}\frac{\tau(n-\tau)}{n}V_{n,1}(\tau)+\sum_{\tau=n-3}^{n-2}\frac{\tau(n-\tau)}{n}V_{n,1}(\tau)\right\}=O\{\tr^2(\bms_1)\}=o\{\var(V_{n,1})\}.
\end{equation*}
Using the facts that
\begin{align*}
	&\var(\X_{1}^{\top}\X_2\X_{3}^{\top}\X_4)=\tr^2(\bms_1^2),\\
	&\cov(\X_{1}^{\top}\X_2\X_{2}^{\top}\X_{3},\X_{1}^{\top}\X_4\X_{4}^{\top}\X_{3})=\tr(\bms_1^4),\\
	&\var(\X_{1}^{\top}\X_2\X_{2}^{\top}\X_{3})=O\{\tr^2(\bms_1^2)\},
\end{align*}
and the Cauchy--Schwarz inequality, we have
\begin{align*}
	&\var\left(\sum_{\tau=4}^{n-4}\frac{\tau(n-\tau)}{n}\{V_n(\tau)-V_{n,1}(\tau)\}\right)\\
	\leq& n\sum_{\tau=4}^{n-4}\var\left( \frac{\tau(n-\tau)}{n}\{V_n(\tau)-V_{n,1}(\tau)\}\right)\\
	=&O\{n\log(n)\tr^2(\bms_1^2)\}+O\{n^2\tr(\bms_1^4)\}=o\{\var(V_{n,1})\}.
\end{align*}
Hence $\var^{1/2}(V_{n})=\var^{1/2}(V_{n,1})\{1+o(1)\}$ and $V_n-V_{n,1}=o_p\{\var^{1/2}(V_{n,1})\}$. 
By Slutsky’s Theorem, we have that $V_{n}/\var^{1/2}(V_{n})\rightarrow N(0,1)$ in distribution as $n,p\rightarrow \infty$. \done

\section{Efficient computation of $M_n$ and $V_n$}

In this section, we present an efficient method for computing the statistics $M_n$ and $V_n$. 
Note that $M_n(\tau)$ and $V_n(\tau)$ have the following equivalent representations:
\begin{equation*}
	M_n(\tau)=\frac{1}{P_{\tau}^{2}}\sum_{i\neq j}^{\tau}\X_{i}^{\top}\X_{j}+\frac{1}{P_{n-\tau}^{2}}\sum_{i,j=\tau+1,i\neq j}^{n}\X_{i}^{\top}\X_{j}-\frac{2}{\tau(n-\tau)}\sum_{i=1}^{\tau}\sum_{j=\tau+1}^{n}\X_{i}^{\top}\X_{j},
\end{equation*}
and 
\begin{equation*}
	V_n(\tau)=A_n(\tau)+B_n(\tau)-2C_n(\tau),
\end{equation*}
where
\begin{align*}
	A_n(\tau)=&\frac{1}{P_{\tau}^{2}}\sum_{i\neq j}^{\tau}(\X_{i}^{\top}\X_{j})^2-\frac{2}{P_{\tau}^{3}}\sum_{1\leq i,j,k\leq \tau}^{\ast}\X_{i}^{\top}\X_{j}\X_{j}^{\top}\X_k+\frac{1}{P_{\tau}^4}\sum_{1\leq i,j,k,l\leq \tau}^{\ast}\X_{i}^{\top}\X_{j}\X_{k}^{\top}\X_{l};\\
	B_n(\tau)=&\frac{1}{P_{n-\tau}^{2}}\sum_{\tau+1\leq i,j\leq n}^{\ast}(\X_{i}^{\top}\X_{j})^2-\frac{2}{P_{n-\tau}^{3}}\sum_{\tau+1\leq i,j,k\leq n}^{\ast}\X_{i}^{\top}\X_{j}\X_{j}^{\top}\X_k\\
	&+\frac{1}{P_{n-\tau}^4}\sum_{\tau+1\leq i,j,k,l\leq n}^{\ast}\X_{i}^{\top}\X_{j}\X_{k}^{\top}\X_{l};\\
	C_n(\tau)=&\frac{1}{P_{\tau}^{1}P_{n-\tau}^{1}}\sum_{j=1}^{\tau}\sum_{j=\tau+1}^{n}(\X_{i}^{\top}\X_{j})^2-\frac{1}{P_{\tau}^{2}P_{n-\tau}^{1}}\sum_{1\leq i,k\leq \tau}^{\ast}\sum_{j=\tau+1}^{n}\X_{i}^{\top}\X_{j}\X_{j}^{\top}\X_k\\
	&-\frac{1}{P_{\tau}^{1}P_{n-\tau}^{2}}\sum_{\tau+1\leq i,k\leq n}^{\ast}\sum_{j=1}^{\tau}\X_{i}^{\top}\X_{j}\X_{j}^{\top}\X_k+\frac{1}{P_{\tau}^{2}P_{n-\tau}^{2}}\sum_{1\leq i,k\leq \tau}^{\ast}\sum_{\tau+1\leq j,l\leq n}^{\ast}\X_{i}^{\top}\X_{j}\X_{k}^{\top}\X_l.
\end{align*}
Based on the equivalent representations, we first describe an efficient procedure for computing $M_n$.
For each possible $\tau$, compute the following quantities: $\sum_{i=1}^{\tau} \X_{i}^{\top} \X_i$, $\sum_{i=\tau+1}^{n} \X_{i}^{\top} \X_i$, $\sum_{i=1}^{\tau} \X_i$, and $\sum_{i=\tau+1}^{n} \X_i$. These sums can be computed with a total time complexity of $O(np)$.
Using these quantities in the subsequent formulas
\begin{equation*}
	\sum_{i\neq j}^{\tau}\X_{i}^{\top}\X_j=\left(\sum_{i=1}^{\tau}\X_i\right)^{\top}\left(\sum_{i=1}^{\tau}\X_i\right)-\sum_{i=1}^{\tau}\X_{i}^{\top}\X_i,
\end{equation*}
\begin{equation*}
	\sum_{i,j=\tau+1,i\neq j}^{n}\X_{i}^{\top}\X_j=\left(\sum_{i=\tau+1}^{n}\X_i\right)^{\top}\left(\sum_{i=\tau+1}^{n}\X_i\right)-\sum_{i=\tau+1}^{n}\X_{i}^{\top}\X_i,
\end{equation*}
\begin{equation*}
	\sum_{i=1}^{\tau}\sum_{j=\tau+1}^{n}\X_{i}^{\top}\X_j=  \left(\sum_{i=1}^{\tau}\X_i\right)^{\top}\left(\sum_{i=\tau+1}^{n}\X_i\right),
\end{equation*}
we can then compute all $M_n(\tau)$'s and obtain $M_n$ with the same time complexity of $O(np)$.

To obtain $V_n$, we first compute $(\X_{i}^{\top}\X_{j})^2$, $i=1,\ldots,n$, $j=1,\ldots, n$.
For each $\tau$, compute the sum $\sum_{i\neq j}^{\tau}(\X_{i}^{\top}\X_j)^2$.
Then, for every $\tau$ and $j$, calculate $\sum_{i=1}^{\tau}(\X_{i}^{\top}\X_j)^2$ and $\X_{j}^{\top}\sum_{i=1}^{\tau}\X_i$.
The time complexity to calculate these quantities is $O(n^2p)$.
Using the formulations
\begin{equation*}
	\sum_{i\neq j\neq k}^{\tau}\X_{i}^{\top}\X_j\X_{j}^{\top}\X_k=\left(\X_{j}^{\top}\sum_{i=1}^{\tau}\X_i-\X_{j}^{\top}\X_j \right)^2-\sum_{i=1}^{\tau}(\X_{i}^{\top}\X_j)^2+ (\X_{j}^{\top}\X_j)^2,
\end{equation*}
\begin{equation*}
	\sum_{i\neq j \neq k \neq l}^{\tau}\X_{i}^{\top}\X_{j}\X_{k}^{\top}\X_{l}=\left(\sum_{i\neq j}^{\tau}\X_i^{\top}\X_j\right)^2-2\sum_{i\neq j}^{\tau}(\X_{i}^{\top}\X_j)^2-4\sum_{i\neq j\neq k}^{\tau}\X_{i}^{\top}\X_j\X_{j}^{\top}\X_k,
\end{equation*}
we can compute all $A_n(\tau)$'s with a time complexity of $O(n^2p)$.
All $B_n(\tau)$'s and $C_n(\tau)$'s can be calculated in a similar manner.
Therefore, all $V_n(\tau)$'s and the final $V_n$ can be obtained with a total time complexity of $O(n^2p)$.

\end{document}